\documentclass[10pt]{article}
\usepackage[utf8]{inputenc}

\usepackage{amsmath,amsbsy,amsthm,amssymb}
\usepackage{multirow,multicol,tabularx,booktabs}
\usepackage{graphicx,placeins,url}
\usepackage{bbm,bm}
\usepackage[ruled]{algorithm2e}
\usepackage[shortlabels]{enumitem}
\usepackage[dvipsnames]{xcolor}
\usepackage[hidelinks, bookmarks=false]{hyperref}

\usepackage[top=2cm,bottom=2cm,right=2.5cm,left=2.5cm]{geometry}

\usepackage{natbib}
\newcommand{\blind}{1}

\newtheorem{conj}{Conjecture}
\newtheorem{thm}[conj]{\bf Theorem}

\newtheorem{cor}[conj]{\bf Corollary}
\newtheorem{prop}[conj]{\bf Proposition}
\newtheorem{lemma}[conj]{\bf Lemma}
\newtheorem{rem}[conj]{\bf Remark}

\newtheorem{assumpt}{\bf Assumption}[section]

\def\bar{\overline}
\def\to{\rightarrow}
\def\To{\longrightarrow}
\def\inprobto{{\buildrel \mathbb P \over \To \,}}
\def\sgn{{\rm\ sign}}

\def\indistrto{\buildrel {D} \over \longrightarrow}

\def\Dc{\mbox{$\mathcal D$}}

\def\Gc{\mbox{$\mathcal G$}}
\def\Hc{\mbox{$\mathcal H$}}

\def\Nc{\mbox{$\mathcal N$}}

\def\Sc{\mbox{$\mathcal S$}}
\def\Tc{\mbox{$\mathcal T$}}

\def\Vc{\mbox{$\mathcal V$}}
\def\Wc{\mbox{$\mathcal W$}}

\def\Zc{\mbox{$\mathcal Z$}}

\def\FF{{\mathbb F}}
\def\HH{{\mathbb H}}

\def\NN{{\mathbb N}}
\def\GG{{\mathbb G}}

\def\Rb{\mbox{$\mathbb R$}}

\def\Zb{\mbox{$\mathbb Z$}}

\def\EE{ {\rm I} \kern-.15em {\rm E} }
\def\PP{ {\rm I} \kern-.15em {\rm P} }

\def\t{ {\bf t}}
\def\u{ {\bf u}}
\def\v{ {\bf v}}
\def\x{ {\bf x}}

\def\z{ {\bf z}}

\def\W{ {\bf W}}
\def\X{ {\bf X}}
\def\Y{ {\bf Y}}
\def\Z{ {\bf Z}}

\def\psibm{ {\bm \psi}}

\def\mds{\medskip}
\def \1{\mathbbm{1} }

\begin{document}

\if1\blind
{
  \title{About Kendall's regression}
  \author{Alexis Derumigny\thanks{
    CREST-ENSAE, 5, avenue Henry Le Chatelier,
    91764 Palaiseau cedex, France. alexis.derumigny@ensae.fr}\hspace{.2cm}
    and
    Jean-David Fermanian\thanks{
    CREST-ENSAE, 5, avenue Henry Le Chatelier,
    91764 Palaiseau cedex, France. jean-david.fermanian@ensae.fr. \newline
    The authors are grateful for helpful discussions with Christian Francq, Johanna Neslehov\'{a}, Alexandre Tsybakov, Jean-Michel Zako\"{\i}an, the participants at the ``Copulas and their Applications'' workshop (Almeria 2017), at the Computational and Financial Econometrics 2017 congress, and at the CREST Financial Econometrics seminar (Feb. 2018). The authors have been supported by the labex Ecodec.
    }}
    
  \maketitle
} \fi

\if0\blind
{
  \bigskip
  \bigskip
  \bigskip
  \begin{center}
    {\LARGE\bf About Kendall's regression}
\end{center}
  \medskip
} \fi


\begin{abstract}
Conditional Kendall's tau is a measure of dependence between two random variables, conditionally on some covariates.
We assume a regression-type relationship between conditional Kendall's tau and some covariates, in a parametric setting with a large number of transformations of a small number of regressors.
This model may be sparse, and the underlying parameter is estimated through a penalized criterion.
We prove non-asymptotic bounds with explicit constants that hold with high probabilities.
We derive the consistency of a two-step estimator, its asymptotic law and some oracle properties.
Some simulations and applications to real data conclude the paper.
\end{abstract}


\noindent%



%
%
%

\section{Introduction}
\label{introduction}

In dependence modeling, it is common to work with scalar dependence measures which are margin-free. They can be used to quantify the positive or negative relationship between two random variables $X_1$ and $X_2$.
One of the most popular of them is Kendall's tau, a dependence measure defined by 
\begin{align*}
    \tau_{1,2}
    := \PP \big( (X_{1,1}-X_{2,1})(X_{1,2}-X_{2,2}) > 0 \big)
    - \PP \big( (X_{1,1}-X_{2,1})(X_{1,2}-X_{2,2}) < 0 \big),
\end{align*}
where $(X_{i,1}, X_{i,2})$, $i=1,2$ are i.i.d. copies of $(X_1, X_2)$, see \citet{nelsen2007introduction}.
When a covariate $\Z$ is available, it is natural to work with the conditional version of this, i.e. the conditional Kendall's tau. It is defined as
\begin{eqnarray*}
\lefteqn{    \tau_{1,2|\Z = \z}
    := \PP \big( (X_{1,1}-X_{2,1})(X_{1,2}-X_{2,2}) > 0 \big| \Z_1 = \Z_2 = \z \big) }\\
    &-& \PP \big( (X_{1,1}-X_{2,1})(X_{1,2}-X_{2,2}) < 0 \big| \Z_1 = \Z_2 = \z \big),
\end{eqnarray*}
where $(X_{i,1}, X_{i,2}, \Z_i)$, $i=1,2$ are i.i.d. copies of $(X_1, X_2, \Z)$.
In such a model, the goal is to study to what extent a $p$-dimensional covariate $\z$ can affect the dependence between the two variables of interest~$X_1$ and $X_2$.

\mds

Most often, it is difficult to have a clear intuition about the functional link between some measure of dependence and the underlying explanatory variables. Sometimes, it is even unclear whether the covariates have an influence on the dependence between the variables of interest. 
This is the so-called ``simplifying assumption'', well-known in the world of copula modeling (see~\cite{derumigny2017tests} and the references therein). This issue is particularly crucial with pair-copula constructions, as pointed out in~\cite{HobaekAasFrigessi},~\cite{AcarGenestNeslehova},~\cite{kurz2017testing}, among others. In our case, we will evaluate an explicit and flexible link between some dependence measure, the Kendall's tau, and the vector of covariates. As a sub-product of our model, we will be able to provide a test of the ``simplifying assumption''.

\mds

Given a dataset $(X_{i,1}, X_{i,2}, \Z_i)$, $i=1,\dots,n$, we will focus on the function $\z \mapsto \tau_{1,2|\Z = \z}$ for $\z \in \Zc$, where $\Zc$ denotes a compact subset of $\Rb^p$. This $\Zc$ represents a set of ``reasonable'' values for $\z$, so that the density $f_\Z$ is bounded from below on $\Zc$.
In order to simplify notations, the reference to the conditioning event $\Z \in \Zc$ will be omitted.
A first natural choice would be to invoke a nonparametric estimator of $\tau_{1,2|\Z = \z}$ as in~\cite{gijbels2011conditional},~\cite{veraverbeke2011ScandinJ} 
and~\cite{derumigny2018kernelBased}. Here, we prefer to obtain parameters that can be interpreted and that would sum up the information about the conditional Kendall's tau. Moreover, kernel-based estimation can be very costly under a computational point of view: for $m$ values of $\z$, the prediction of all these conditional Kendall's taus has a total cost of $O(m n^2)$, that can be large if a large number $m$ is required.
Other estimators of the conditional Kendall's tau, based on classification methods, are proposed in~\cite{derumigny2018classification}.

\mds

In this paper, our idea is to decompose the function $\z \mapsto \tau_{1,2|\Z = \z}$ on some functional basis $(\psi_i)_{i \geq 1}$, as any element of a space of functions from $\Zc$ to $\Rb$. First note that a Kendall's tau takes its values in the interval $[-1, 1]$, and not on the whole real line.
Nevertheless, for some known increasing and continuously differentiable function $\Lambda: [-1, 1] \to \Rb$, the function $\z \mapsto \Lambda \left( \tau_{1,2|\Z=\z} \right)$ takes values on up to the whole real line potentially, and it can be decomposed on any basis $(\psi_i)_{i \geq 1}$.
Typical transforms are $\Lambda(\tau) = \log \big(\frac{1+\tau}{1-\tau} \big)$ (the Fisher transform) or $\Lambda(\tau) = \log(-\log((1-\tau)/2))$.
We will assume that only a finite number of elements are necessary to represent this function. This means that we have
\begin{equation}
    \Lambda \left( \tau_{1,2|\Z=\z} \right)
    = \sum_{i=1}^{p'} \psi_i(\z) \beta^*_i
    = \psibm (\z)^T \beta^*,
    \label{model:lambda_cond_tau_Z}
\end{equation}
for all $\z \in \Zc$, with $p' > 0$ and a ``true'' unknown parameter $\beta^* \in \Rb^{p'}$.
The function $\psibm(\cdot) := \big(\psi_1(\cdot), \dots, \psi_{p'}(\cdot) \big)^T$ from $\Rb^p$ to $ \Rb^{p'}$ is known and corresponds to deterministic transformations of the covariates $\z$.
In practice, it is not easy to have intuition about which kind of basis to use, especially in our framework of conditional dependence measurement. Therefore, the most simple solution is the use of a lot of different functions : polynomials, exponentials, sinuses and cosinuses, indicator functions, etc... 
They allow to take into account potential non-linearities and even discontinuities of
conditional Kendall's taus with respect to $\z$.
For the sake of identifiability, we only require their linear independence, as seen in the following proposition (whose straightforward proof is omitted).
\begin{prop}
    The parameter $\beta^*$ in Model (\ref{model:lambda_cond_tau_Z}) is identifiable if and only if the functions
    $(\psi_1, \dots, \psi_{p'})$ are linearly independent $\PP_\Z$-a.e. in the sense that, for any given vector $\t = (t_1, \dots, t_{p'}) \in \Rb^{p'}$, $\PP_\Z \big( \psibm(\Z)^T \t = 0 \big) = 1$ implies $\t = 0$.
    \label{prop:identifiability_condition}
\end{prop}
%
%
With such a large choice among flexible classes of functions, it is unlikely we will be able to guess the right ones ex ante. 
Therefore, it will be necessary to consider a large number of functions $\psi_i$ under a sparsity constraint:
the cardinality of $\Sc$, the set of non-zero components of $\beta^*$, is less than some $s \in \{1, \dots, p'\}$. It is denoted by $|\Sc|=|\beta^*|_0$, where
$| \cdot |_0$ yields the number of non-zero components of any vector in $\Rb^{p'}$.
Note that, in this framework, $p'$ can be moderately large, for example $10$ or $30$ while the original dimension $p$ is small, for example $p=1$ or~$2$.
This corresponds to the decomposition of a function, defined on a small-dimension domain, in a mildly large basis.

\mds

Once an estimator $\hat \beta$ of $\beta^*$ has been computed, the prediction of all the conditional Kendall's tau's for $m$ values of $\z$, which is just the computation of
$\Lambda^{(-1)} \big( \psibm (\z)^T \hat \beta \big)$
can be done in $O(m s)$, that is much faster than what was previously required with a kernel-based estimator for large $m$, as soon as $s \leq n^2$ (see Section~\ref{example:computation_time} for a discussion).

\mds

Estimating Model~(\ref{model:lambda_cond_tau_Z}) not only provides an estimator of the conditional Kendall's tau $\tau_{1,2|\Z=\z}$, but also easily provides estimators of the marginal effects of $\z$ as by-product. For example, given $\z \in \Zc$, the marginal effect of $z_1$, i.e. $\partial \tau_{1,2|\Z=\z}(\z) / \partial z_1 $, can be directly estimated by $\big( \partial_{z_1} \psibm(\z) \big)^T \hat \beta \cdot \Lambda^{(-1)}{}' \big( \psibm (\z)^T \hat \beta \big)$, assuming that $\psibm$ and $\Lambda^{(-1)}$ are differentiable respectively at $\z$ and $\psibm (\z)^T \hat \beta$. Such sensitivities can be useful in many applications.

\mds

A desirable empirical feature of Model~(\ref{model:lambda_cond_tau_Z}) would be the possibility of obtaining
very high/low levels of dependence between $X_1$ and $X_2$, for some $\Z$ values, i.e. $\Lambda^{(-1)}(\psibm (\z)^T \beta^*)$
should be close (or even equal) to $1$ or $-1$ for some $\z$. This can be the case even if $\Zc$ is compact, that is here required for theoretical reasons.
Indeed, the image of $\{ \tau_{1,2|\z} \vert \z\in \Zc \}=[\tau_{\min},\tau_{\max}]$ through $\Lambda$ is an interval $[\Lambda_{\min},\Lambda_{\max}]$.
If $\psibm (\z)^T\beta^*\geq \Lambda_{\max}$ (resp. $\psibm (\z)^T\beta^*\leq \Lambda_{\min}$), then simply set $\tau_{1,2|\Z=\z}=\tau_{\max}$ or even one (resp. $\tau_{1,2|\Z=\z}=\tau_{\min}$ or even $(-1)$).

\mds

Contrary to more usual models, the ``explained variable'' - the conditional Kendall's tau $\tau_{1,2|\Z=\z}$ - is not observed in (\ref{model:lambda_cond_tau_Z}). Therefore, a direct estimation of the parameter $\beta^*$ (for example, by the ordinary least squares, or by the Lasso) is unfeasible.
In other words, even if the function $\z \mapsto \Lambda \big( \tau_{1,2|\Z=\z} \big)$ is deterministic, finding the best $\beta$ in Model (\ref{model:lambda_cond_tau_Z}) is far from being just a numerical analysis problem since the function to be decomposed is unknown.
Nevertheless, we will replace $\tau_{1,2|\Z=\z}$ by a nonparametric estimate $\hat \tau_{1,2|\Z=\z}$, and use it as an approximation of the explained variable.
More precisely, we fix a finite collection of points $\z'_1, \dots, \z'_{n'} \in \Zc^{n'}$ and we estimate $\hat \tau_{1,2|\Z=\z}$ for each of these points.\
Then, $\hat \beta$ is estimated as the minimizer of the $l_1$-penalized criteria
\begin{equation}
    \hat \beta := \arg \min_{\beta \in \Rb^{p'}}
    \Big[ \frac{1}{n'} \sum_{i=1}^{n'} \left( \Lambda(\hat \tau_{1,2|\Z=\z'_i}) - \psibm (\z'_i)^T\beta \right)^2 + \lambda |\beta|_1 \Big],
    \label{def:estimator_hat_beta}
\end{equation}
where $\lambda$ is a positive tuning parameter (that may depend on $n$ and $n'$), and $|\cdot|_q$ denotes the $l_q$ norm, for $1\leq q \leq \infty$. This procedure is summed up in the following Algorithm~\ref{algo:estimation_beta}.
Note that even if we study the general case with any $\lambda \geq 0$, the properties of the unpenalized estimator can be derived by choosing the particular case $\lambda = 0$.

\begin{algorithm}[htb]
\label{algo:estimation_beta}
\SetAlgoLined
    \vspace{0.1cm}
    \KwIn{A dataset $(X_{i,1}, X_{i,2}, \Z_i)$, $i=1,\dots,n$}
    \KwIn{A finite collection of points $\z'_1, \dots, \z'_{n'} \in \Zc^{n'}$}
    \For{$j\leftarrow 1$ \KwTo $n'$} {
        Compute the estimator $\hat \tau_{1,2|\Z=\z'_j}$ using the sample $(X_{i,1}, X_{i,2}, \Z_i)$, $i=1,\dots,n$ \;
    }
    Compute the minimizer $\hat \beta$ of (\ref{def:estimator_hat_beta}) using the $\hat \tau_{1,2|\Z=\z'_j},$ $j=1,\dots,n',$ estimated in the above step \;
    \KwOut{An estimator $\hat \beta$.}
\caption{Two-step estimation of $\beta$}
\end{algorithm}

\mds

Several nonparametric estimators of $\hat \tau_{1,2|\Z=\z'_j}$ can potentially be used. We refer to~\cite{derumigny2018kernelBased} for a detailed analysis of their statistical properties. 
They are of the form
\begin{equation}
    \hat \tau_{1,2|\Z=\z}
    := \sum_{i=1}^n \sum_{j=1}^n w_{i,n}(\z) w_{j,n}(\z)
    g^*(\X_i, \X_j),
    \label{def:hat_conditional_tau}
\end{equation}
where $g^*$ is a bounded function, $\X_i := (X_{i,1} , X_{i,2})$ for $i=1, \dots, n$ and $w_{i,n}(\z) := K_h(\Z_i-\z) / \sum_{j=1}^n K_h(\Z_j-\z)$, $h = h(n) > 0$ denoting the bandwidth sequence.
In the same way, the conditional Kendall's tau can be rewritten as $\tau_{1,2|\Z=\z} = \EE[g^*(\X_1, \X_2) | \Z_1 = \Z_2 = \z]$ for the same choices of $g^*$.
Possible choices of $g^*$ are given in Section~\ref{section:choice_g}.

\mds

In Section~\ref{section:finite_distance_bounds}, we state non-asymptotic results for the our estimator $\hat \beta$ that hold with high probability. In Section~\ref{section:asymptotic_case}, its asymptotic properties are stated.
In particular, we will study the cases when $n'$ is fixed and $n\rightarrow \infty$, and when both indices tend to the infinity.
We also give some oracle properties and suggest a related adaptive estimator.
Sections~\ref{section:simulations} and~\ref{section:real_data} illustrate respectively the numerical performances of $\hat \beta$ on simulated and real data. 
All proofs and two supplementary figures have been postponed into the supplementary material.

\begin{rem}
	At first sight, in Model (\ref{model:lambda_cond_tau_Z}), there seems to be no noise perturbing the variable of interest. In fact, this is a simple consequence of our formulation of the model. In the same way, a classical linear model $Y=\X^T \beta^* + \varepsilon$ can be rewritten as $\EE[Y|\X=\x] = \x^T \beta^*$ without any explicit noise. By definition, $\EE[Y|\X=\x]$ is a deterministic function of a given $\x$. In our case, $\Lambda \big( \tau_{1,2|\Z=\z} \big)$ is a deterministic function of the variable $\z$. This means that we cannot formally write a model with noise, such as $\Lambda \big( \tau_{1,2|\Z=\z} \big) = \psibm (\z)^T \beta^* + \varepsilon$ where $\varepsilon$ is independent of the choice of $\z$. Indeed, the left-hand side of the latter equality is a $\z$-mesurable quantity, unless $\varepsilon$ is constant almost surely.
\end{rem}

\begin{rem}
    Note that the conditioning event of Model (\ref{model:lambda_cond_tau_Z}) is unusual: usual regression models consider $\EE[g(\X)|\Z=\z]$ as a function of the conditioning variable $\z$. Here, the probabilities of concordant/discordant pairs are made conditionally on $\Z_1 = \Z_2 = \z$. This unusual conditioning event will necessitate some peculiar theoretical treatments.
\end{rem}

\begin{rem}
    Instead of a fixed design setting $(\z'_i)_{i=1,\ldots,n'}$ in the optimization program, it would be possible to consider a random design: simply draw $n'$ realizations of $\Z$, independently of the $n$-sample that has been used for the estimation of the conditional Kendall's taus.
    The differences between fixed and random designs are mainly a matter of presentation and the reader could easily rewrite our results in a random design setting.
    We have preferred the former one to study the finite distance properties and asymptotics when $n'$ is fixed (Section~\ref{Asymp_n_fixed_nprime}). When $n$ and $n'$ will tend to the infinity (Section~\ref{Asymp_n_nprime}), both designs are encompassed de facto because we will assume the weak convergence of the empirical distribution associated to the sample $(\z'_i)_{i=1,\ldots,n'}$, when $n'\rightarrow \infty$.
    \label{rem:choice_sample_z_prime}
\end{rem}

\section{Finite-distance bounds on $\hat \beta$}
\label{section:finite_distance_bounds}

Our first goal is to prove finite-distance bounds in probability for the estimator $\hat \beta$.
Let $\Zb'$ be the matrix of size $n' \times p'$ whose lines are $\psibm (\z'_i)^T$, $i=1,\ldots,n'$, and let $\Y \in \Rb^{n'}$ be the column vector whose components are
$Y_i = \Lambda(\hat \tau_{1,2|\Z=\z'_i})$, $i=1, \dots, n'$.
For a vector $\v \in \Rb^{p'}$, denote by $||\v||_{n'}:= |\v|_2/\sqrt{n'}$ its empirical norm.
We can then rewrite the criterion~(\ref{def:estimator_hat_beta}) as
$\hat \beta := \arg \min_{\beta \in \Rb^{p'}}
    \Big[ ||\Y - \Zb' \beta||_{n'}^2 + \lambda |\beta|_1 \Big],$
where $\Y$ and $\Zb'$ may be considered as ``observed'', so that the practical problem is reduced to a standard Lasso estimation procedure.
Define some ``residuals'' by
$\xi_{i,n} := \Lambda(\hat \tau_{1,2|\Z=\z'_i}) - \psibm (\z'_i)^T\beta^*
= \Lambda(\hat \tau_{1,2|\Z=\z'_i}) - \Lambda(\tau_{1,2|\Z=\z'_i}) ,$
for~$i=1,\dots,n'$.
Note that these $\xi_{i,n}$ are not ``true residuals'' in the sense that they do not depend on the estimator $\hat \beta$, but on the true parameter $\beta^*$.
We also emphasized the dependence on $n$ in the notation $\xi_{i,n}$, which is a consequence of the estimated conditional Kendall's tau.

\mds

To get non-asymptotic bounds on $\hat \beta$, assume the \emph{Restricted Eigenvalue} (RE) condition, introduced by~\cite{bickel2009simultaneous}.
For $c_0 > 0$ and $s\in \{ 1, \dots, p \}$, assume

\medskip

\noindent
$RE(s,c_0)$ \textbf{condition :}
\emph{The design matrix $\Zb'$ satisfies}
\begin{equation*}
    \kappa(s, c_0) := \min_{\small
    \begin{array}{c}
        J_0 \subset \{1, \dots, p'\} \\ Card(J_0) \leq s
    \end{array} }
    \min_{\small
    \begin{array}{c}
        \delta \neq 0 \\ |\delta_{J_0^C}|_1 \leq c_0 |\delta_{J_0}|_1
    \end{array} }
    \dfrac{|\Zb' \delta|_2}{\sqrt{n'} |\delta|_2} > 0.
\end{equation*}

Note that this condition is very mild, and is satisfied with a high probability for a large class of random matrices: see~\citet[Section 8.1]{bellec2016slope} for references and a discussion.

\begin{assumpt}
    The function $\z\mapsto \psibm(\z)$ are bounded on $\Zc$ by a constant $C_\psibm$. Moreover, $\Lambda(\cdot)$ is continuously differentiable.
    Let $\Tc$ be the range of $\z\mapsto \tau_{1,2|\Z=\z}$, from $\Zc$ towards $[-1,1]$.
    On an open neighborhood of $\Tc$, the derivative of $\Lambda(\cdot)$ is bounded by a constant $C_{\Lambda'}$.
    \label{assumpt:compact_case}
\end{assumpt}

\begin{thm}[Fixed design case]
    Suppose that Assumptions \ref{assumpt:kernel_integral}-\ref {assumpt:f_XZ_Holder} and~\ref{assumpt:compact_case}
    hold and that the design matrix $\Zb'$ satisfies the $RE(s,3)$ condition.
    Choose the tuning parameter as $\lambda = \gamma t$,
    with $\gamma \geq 4$ and $t>0$, and assume that we choose $h$ small enough such that
    \begin{align}
        h^\alpha &\leq
        \min \bigg( \frac{f_{\Z, min} \alpha !}
        {4 \,  C_{K, \alpha}} \; , \;
        \frac{f_{\Z, min}^4 \alpha! \, t}
        {8 \, C_{\psi } C_{\Lambda'} (f_{\Z, min}^2 + 8 f_{\Z, max}^2)
        C_{\X\Z, \alpha}  }
        \bigg)
        \label{cond:h_1_2} .
    \end{align}
    Then, we have
    \begin{align}
        \PP \Big( & ||\Zb'(\hat \beta - \beta^*)||_{n'}
        \leq \dfrac{4(\gamma+1) t \sqrt{s}}{\kappa(s,3)}
        \text{ and } |\hat \beta - \beta^*|_q
        \leq \dfrac{4^{2/q}(\gamma+1) t s^{1/q} }{\kappa^2(s,3)},
        \text{ for every } 1 \leq q \leq 2 \Big) \nonumber \\
        &\hspace{3cm} \geq 1 - 2 n' \exp \Big( - n h^p C_1 \Big)
        - 2 n' \exp \Big( - \frac{(n-1) h^{2p} t^2}
        {C_2 + C_3  t} \Big),
        \label{eq:bound_proba_hat_beta}
    \end{align}
    where $C_1 := f_{\Z, min}^2 / \big(32 f_{\Z, max} \int K^2 + (8/3) C_K f_{\Z, min} \big)$,
    $C_3 := (64/3) C_{\psi } C_{\Lambda'} C_K^2 (f_{\Z, min}^2 + 8 f_{\Z, max}^2)/f_{\Z, min}^4$,
    and $C_2 :=  \{16 C_{\psi} C_{\Lambda'} (f_{\Z, min}^2 + 8 f_{\Z, max}^2) f_{\Z,max}\int K^2 \}^2 / f_{\Z, min}^8$.
    \label{thm:bound_proba_hat_beta}
\end{thm}

This theorem, proved in Section \ref{proof:thm:bound_proba_hat_beta}, yields some bounds that hold in probability for the prediction error $||\Zb'(\hat \beta - \beta^*)||_{n'}$ and for the estimation error $|\hat \beta - \beta^*|_q$, $1 \leq q \leq 2$, under the specification (\ref{model:lambda_cond_tau_Z}). Note that the influence of $n'$ and $p'$ is hidden through the Restricted Eigenvalue number $\kappa(s,3)$. The result depends on three parameters $\gamma$, $t$ and $h$. Apparently, the choice of $\gamma$ seems to be easy, as a larger $\gamma$ deteriorates the upper bounds.
Nonetheless, it is a bit misleading because $\hat\beta$ implicitly depends on $\lambda$ and then on $\gamma$ (for a fixed $t$). Nonetheless, choosing $\gamma =4$ is a reasonable ``by default'' choice.
Moreover, a lower $t$ provides a smaller upper bound, but at the same time the probability of this event is lowered.
This induces a trade-off between the probability of the desired event and the size of the bound, as we want the smallest possible bound with the highest probability.
Moreover, we cannot choose a too small $t$, because of the lower bound (\ref{cond:h_1_2}): $t$ is limited by a value proportional to $h^\alpha$. The latter $h$ cannot be chosen as too small, otherwise the probability in Equation (\ref{eq:bound_proba_hat_beta}) will decrease.
To be short:
\textit{low values of $h$ and $t$ yield a sharper upper bound with a lower probability, and the opposite}.
Therefore, a trade-off has to be found, depending of the kind of result we are interested in.

\medskip

Clearly, we would like to exhibit the sharpest upper bounds in~(\ref{eq:bound_proba_hat_beta}), with the ``highest probabilities''.
Let us look for parameters of the form $t \propto n^{-a}$ and $h \propto n^{-b}$, with $a, b > 0$.
The assumptions of Theorem \ref{thm:bound_proba_hat_beta} imply $b \alpha \geq a$ (to satisfy~(\ref{cond:h_1_2}))
and $1 - 2 a - 2 p b > 0$ (so that the right-hand side of (\ref{eq:bound_proba_hat_beta}) tends to $1$ as $n \to \infty$, i.e. $nh^{p} \to \infty$ and $n t^2 h^{2p} \to \infty$).
For fixed $\alpha$ and $ p $, what are the ``optimal'' choices $a$ and $b$ under the constraints $b \alpha \geq a$ and $1 - 2 a - 2 p b > 0$ ?
The latter domain is the interior of a triangle in the plane $(a,b)\in \Rb_+^2$, whose vertices are $O:=(0,0)$, $A:=(0,1/(2p))$ and $B:=(\alpha/(2p+2\alpha),1/(2p+2\alpha))$, plus the segment $]0,B[$.
All points in such a domain would provide admissible couples $(a,b)$ and then admissible tuning parameters $(t,h)$.
In particular, choosing the neighborhood of $B$, i.e. $a=\alpha (1-\epsilon)/(2p+2\alpha)$ and $b=1/(2p+2\alpha)$ for some (small) $\epsilon >0$, will be nice because the upper bounds will be minimized.

\begin{cor}
    For $0 < \epsilon < 1$, choosing the parameters $\lambda=4t$,
    $t = (n-1)^{ - \alpha (1-\epsilon)/(2 \alpha + 2p)}$ and
    $$h =c_h (n-1)^{ - 1/(2 \alpha + 2p)} ,\; c_h:= \Big(\frac{f_{\Z, min}^4  \alpha  !}
    {2 \, C_{\psi } C_{\Lambda'} (f_{\Z, min}^2 + 16 f_{\Z, max}^2)
    C_{\X\Z, \alpha}  } \Big)^{1/\alpha},$$
    we have, if $n$ is sufficiently large so that~(\ref{cond:h_1_2}) is satisfied,
    \begin{align*}
        \PP \Big( & ||\Zb'(\hat \beta - \beta^*)||_{n'}
        \leq \dfrac{20 \sqrt{s}}{\kappa(s,3) (n-1)^{
        \alpha (1 - \epsilon)/(2 \alpha + 2p)}} \text{\; and }  \\
        & |\hat \beta - \beta^*|_q
        \leq \dfrac{5 . 4^{2/q} s^{1/q} }{\kappa^2(s,3) (n-1)^{\alpha (1 - \epsilon)/(2 \alpha + 2p)}},
        \text{ for every } 1 \leq q \leq 2 \Big) \nonumber \\
        &\geq 1 - 2 n' \exp \Big( - C_1 c_h^{p} (n-1)^{(2 \alpha + p )/(2 \alpha + 2 p) }\Big)
        - 2 n' \exp \Big( - \dfrac{c_h^{2p} (n-1)^{2\alpha\epsilon / (2p+2\alpha)} }
        {C_2 + C_3 (n-1)^{ - \alpha (1-\epsilon)/(2 \alpha + 2p)}} \Big).
    \end{align*}
\end{cor}

%

\section{Asymptotic behavior of $\hat \beta$}
\label{section:asymptotic_case}

\subsection{Asymptotic properties of $\hat \beta$ when $n\to\infty$ and for fixed $n'$}
\label{Asymp_n_fixed_nprime}

In this part, $n'$ is still supposed to be fixed and we state the consistency and the asymptotic normality of $\hat\beta$ as $n \to \infty$.
As above, we adopt a fixed design: the $\z'_i$ are arbitrarily fixed or, equivalently, our reasonings are made conditionally on the second sample.

\mds

For $n, n' > 0$, denote by $\hat \beta_{n, n'}$ the estimator (\ref{def:estimator_hat_beta}) with $h = h_{n}$ and $\lambda = \lambda_{n, n'}$.
The following lemma, proved in Section \ref{proof:lemma:beta_process_GG}, provides another representation of this estimator $\hat \beta_{n, n'}$ that will be useful hereafter.

\begin{lemma}
    We have $\hat \beta_{n, n'} = \arg \min_{\beta \in \Rb^{p'}}
    \GG_{n,n'}(\beta)$, where
    \begin{equation}
        \GG_{n,n'}(\beta)
        := \frac{2}{n'} \sum_{i=1}^{n'} \xi_{i,n} \psibm (\z'_i)^T (\beta^* - \beta)
        + \frac{1}{n'} \sum_{i=1}^{n'}
        \big\{ \psibm (\z'_i)^T (\beta^* - \beta) \big\}^2
        + \lambda_{n, n'} |\beta|_1.
        \label{def:process_GG}
    \end{equation}
    \label{lemma:beta_process_GG}
\end{lemma}

We will invoke a \textit{convexity argument}:
``Let $g_n$ and $g_\infty$ be random convex functions taking minimum values at $x_n$ and $x_\infty$, respectively.
If all finite dimensional distributions of $g_n$ converge weakly to those of $g_\infty$ and $x_\infty$ is the unique minimum point of $g_\infty$
with probability one, then $x_n$ converges weakly to $x_\infty$'' (see~\cite{kato2009asymptotics}, e.g).

\begin{thm}[Consistency of $\hat \beta$]
    Under the assumptions of Lemma~\ref{lemma:consistency_hatTau},
    if $n'$ is fixed and $\lambda = \lambda_{n, n'} \to \lambda_0 $,
    then, given $\z'_1,\ldots, \z'_{n'}$ and as $n$ tends to the infinity, $\hat \beta_{n, n'} \inprobto \beta^{**} := \inf_\beta \GG_{\infty,n'} (\beta),$ where
    $\GG_{\infty,n'} (\beta)  :=
    \sum_{i=1}^{n'} \big( \psibm (\z'_i)^T (\beta^* - \beta) \big)^2 / n'+ \lambda_0 |\beta |_1.$
    In particular, if $\lambda_0=0$ and
        $< \psibm (\z'_1), \dots, \psibm (\z'_{n'}) > \, = \Rb^{p'}$, then
    $\hat \beta_{n, n'} \inprobto \beta^*$.
\end{thm}

{\it Proof :} By Lemma \ref{lemma:consistency_hatTau}, the first term in the r.h.s. of~(\ref{def:process_GG})
converges to $0$ as $n\to\infty$.
The third term in the r.h.s. of~(\ref{def:process_GG}) converges to $\lambda_0 |\beta |_1$ by assumption.
We have just proven that $\GG_{n,n'} \to \GG_{\infty,n'}$ pointwise as $n \to \infty$.
We can now apply the convexity argument, because $\GG_{n,n'}$ and $\GG_{\infty,n'}$ are convex functions.
As a consequence, $\arg \min_\beta \GG_{n,n'}(\beta) \to \arg \min_\beta \GG_{\infty,n'}(\beta)$ in law.
Since we have adopted a fixed design setting, $\beta^{**}$ is non random, given $(\Z_1',\ldots,\Z_{n'}')$.
The convergence in law towards a deterministic quantity implies convergence in probability, which concludes the proof.
Moreover, when $\lambda_0=0$, $\beta^*$ is the minimum of $\GG_{\infty,n'}$ because the vectors $\psibm (\z'_i)$, $i=1, \dots, p'$ generate the space $\Rb^{p'}$.
Therefore, this implies the consistency of $\hat\beta_{n,n'}$. $\;\;\Box$

\mds

To evaluate the limiting behavior of $\hat\beta_{n,n'}$,
we need the joint asymptotic normality of $(\xi_{1,n},\ldots,\xi_{n',n})$, when $n\to\infty$ and given $\z'_1,\ldots,\z'_{n'}$.
By applying the Delta-method to the function $\Lambda(\, \cdot \,)$ component-wise, this is given by the following corollary of Lemma \ref{lemma:asymptNorm_hatTau}.

\begin{cor}
    Under the assumptions of Lemma \ref{lemma:asymptNorm_hatTau},
    $(n h_{n}^p)^{1/2}
    \left[ \xi_{1,n},\ldots, \xi_{n',n} \right]^T$ tends in law towards a random vector $\Nc \big(0, \tilde \HH \big)$ given $(\z'_1,\ldots,\z'_{n'})$,
    where $\tilde \HH$ is a $n' \times n'$ real matrix defined, for every integers $1 \leq i, j \leq n'$, by
    \begin{align*}
        [\tilde \HH]_{i, j} &:= \frac{ 4 \int K^2 \1_{ \{ \z'_{i} = \z'_{j} \} }
        }{ f_\Z(\z'_{i})}
        \Big( \Lambda' \big(\tau_{1,2|\Z=\z'_{i}} \big) \Big)^2 
        \times \Big\{ \EE[ \tilde g(\X_1,\X) \tilde g(\X_2,\X) | \Z = \Z_1 = \Z_2 = \z'_{i}] - \tau_{1,2|\Z=\z'_{i}}^2 \Big\},
    \end{align*}
    where $\tilde g$ is the symmetrized version $\tilde g(\x_1, \x_2) := (g^*(\x_1, \x_2) + g^*(\x_2, \x_1))/2$.
    \label{cor:asymptNorm_xi}
\end{cor}

\begin{thm}[Asymptotic law of the estimator]
\label{thm:WeakConvLasso}
    Under the assumptions of Lemma \ref{lemma:asymptNorm_hatTau}, and if
    $\lambda_{n, n'}  (n h_{n,n'}^p)^{1/2}$ tends to $ \ell $ when $n\rightarrow \infty$, we have
    $(n h_{n,n'}^p)^{1/2} (\hat \beta_{n, n'} - \beta^*)
    \indistrto \u^* := \arg \min_{\u \in \Rb^{p'}} \FF_{\infty, n'}(\u),$
    given $\z'_1,\ldots, \z'_{n'}$, where
    \begin{align*}
        \FF_{\infty, n'}(\u)
        &:= \frac{2}{n'} \sum_{i=1}^{n'} \sum_{j=1}^{p'}
        W_{i} \psi_j (\z'_i) u_{j}
        + \frac{1}{n'} \sum_{i=1}^{n'} \left( \psibm (\z'_i)^T \u \right)^2
        + \ell \sum_{i=1}^{p'} \big( |u_i| \1_{\{\beta_i^*=0\}} + u_i\sgn(\beta_i^*) \1_{\{\beta_i^* \neq 0\}} \big) ,
    \end{align*}
    with $\W = (W_1, \dots, W_{n'}) \sim \Nc \left(0, \tilde \HH \right).$
\end{thm}

This theorem is proved in Section \ref{proof:thm:WeakConvLasso}.
When $\ell=0$, we can say more about the limiting law in general. Indeed, in such a case,
$\u^*=\arg \min_{\u \in \Rb^{p'}} \FF_{\infty, n'}(\u)$ is the solution of the first order conditions $\nabla  \FF_{\infty, n'}(\u)=0$, that are written as
$\sum_{i=1}^{n'} W_{i} \psibm (\z'_i)
+ \sum_{i=1}^{n'} \psibm (\z'_i)  \psibm (\z'_i)^T \u = 0.$
Therefore,
$$ \u^* = - \Big( \sum_{i=1}^{n'} \psibm (\z'_i) \psibm (\z'_i)^T \Big)^{-1}
\sum_{i=1}^{n'} W_{i} \psibm (\z'_i),$$
when $\Sigma_{n'}:=\sum_{i=1}^{n'} \psibm (\z'_i)  \psibm (\z'_i)^T $ is invertible. Then, the limiting law of
$(n h_{n,n'}^p)^{1/2} (\hat \beta_{n, n'} - \beta^*)$ is Gaussian, and its asymptotic covariance is
$ V_{as}:=\Sigma_{n'}^{-1} \sum_{i,j=1}^{n'} [\tilde\HH]_{i,j}
\psibm(\z'_i) \psibm(\z'_j)^T\Sigma_{n'}^{-1}. $

\mds

The previous results on the asymptotic normality of $\hat\beta_{n,n'}-\beta^* $ can be used to test $\Hc_0: \beta^*=0$ against the opposite. 
As said in the introduction, this would constitute a test of the ``simplifying assumption'', i.e. the fact that the conditional copula of $(X_1,X_2)$ given 
$\Z$ does not depend on this covariate. Some tests of significance of $\beta^*$ would be significantly simpler than most of the tests of the simplifying assumption that have been proposed in the literature until now. Indeed, the latter ones have been built on nonparametric estimates of conditional copulas and, as sub-products of the weak convergence of the associated processes, the test statistics behaviors are obtained. Therefore, such statistics depend on a preliminary non-parametric estimation of conditional marginal distributions (see~\cite{veraverbeke2011ScandinJ},~\cite{derumigny2017tests}, e.g.), a source of complexities and statistical noise. 
At the opposite, some tests of $\Hc_0$ based on $\hat\beta_{n,n'}$ do not require this stage, at the cost of a (probably small) loss of power.
For instance, in the case of $\ell=0$, we propose the Wald-type test statistics 
$$\Wc_n:= n h_{n,n'}^p (\hat \beta_{n, n'} - \beta^*)^T V_n (\hat \beta_{n, n'} - \beta^*),\; V_{n}:=\Sigma_{n'}^{-1} \sum_{i,j=1}^{n'} \hat\HH_{i,j}
\psibm(\z'_i) \psibm(\z'_j)^T\Sigma_{n'}^{-1}. $$
\begin{align*}
    \hat \HH_{i, j}
    &:= \frac{ 4 \int K^2 \1_{ \{ \z'_{i} = \z'_{j} \} }
    }{ \hat f_\Z(\z'_{i})}
    \Big( \Lambda' \big(\hat\tau_{1,2|\Z=\z'_{i}} \big) \Big)^2 
    \times \Big\{ \Gc_n(\z'_i) - \hat\tau_{1,2|\Z=\z'_{i}}^2 \Big\},
\end{align*}
where $\hat{f}_\Z(\Z)$ and $\Gc_n(\z)$ denote consistent estimators of $f_\Z(\z)$ and 
$\EE[ \tilde g(\X_1,\X) \tilde g(\X_2,\X) | \Z = \Z_1 = \Z_2 = \z]$ respectively. Under $\Hc_0$, $\Wc_n$ tends to a chi-square distribution with $n'$ degrees of freedom. For instance, with the notations of Section~\ref{introduction}, we propose
$$ \Gc_n(\z)= \sum_{i,j,k=1, i\neq j \neq k}^n w_{i,n}(\z)w_{j,n}(\z) w_{k,n}(\z) \tilde g(\X_{i},\X_{k}) \tilde g(\X_j,\X_k) .$$

Note that if there is an intercept, i.e. if one of the functions in $\psibm$ (say, $\psi_1$) is constant to $1$, it should be removed in the statistics above. The corresponding coefficients of $\hat \beta$ should be removed as well. Indeed, in this case the simplifying assumption does not correspond to $\beta^* = 0$, but rather to $\beta^*_{-1} = 0$ where $\beta^*_{-i}$ denotes the vector $\beta^*$ where the $i$-th coefficient has been removed.
    
\subsection{Oracle property and a related adaptive procedure}

Let remember that $\Sc:=\{j: \beta_j^*\neq 0\}$ and assume that $|\Sc|=s<p$ so that the true model depends on a subset of predictors. In the same spirit as~\cite{fanli2001}, we say that an estimator $\hat\beta$ satisfies the oracle property if
\begin{itemize}
    \item $v_{n}(\hat \beta_{\Sc} - \beta^*_{\Sc})$ converges in law towards a continuous random vector, for some conveniently chosen rate of convergence $(v_{n})$, and
    \item  we identify the nonzero components of the true parameter $\beta^*$ with probability one when the sample size $n$ is large, i.e. the probability of the event $\big( \{j:\hat\beta_j \neq 0\} = \Sc \big)$ tends to one.
\end{itemize}
As above, let us fix $n'$ and $n$ will tend to the infinity. Then, denote $\{j:\hat\beta_j \neq 0\}$ by $\Sc_n$, that will  implicitly depend on $n'$.
It is well-known that the usual Lasso estimator does not fulfill the oracle property, see ~\cite{zou2006}. Here, this is still the case.
The following proposition is proved in Section \ref{proof:prop:not_oracle_property}.
\begin{prop}
    Under the assumptions of Theorem~\ref{thm:WeakConvLasso}, $\lim\sup_n  \PP\left( \Sc_n= \Sc \right) =c <1. $
    \label{prop:not_oracle_property}
\end{prop}
A usual way of obtaining the oracle property is to modify our estimator in an ``adaptive'' way. Following~\cite{zou2006}, consider a preliminary ``rough'' estimator of $\beta^*$, denoted by $\tilde \beta_n$, or more simply $\tilde\beta$. Moreover $\nu_n (\tilde \beta_n - \beta^*)$ is assumed to be asymptotically normal, for some deterministic sequence $(\nu_n)$ that tends to the infinity.
Now, let us consider the same optimization program as in~(\ref{def:estimator_hat_beta}) but with a random tuning parameter given by
$   \lambda_{n,n'} := \mu_{n,n'}/ |\tilde{\beta}_n|^\delta$,
for some constant $\delta >0$ and some positive deterministic sequence $(\mu_{n,n'})$. The corresponding adaptive estimator (solution of the modified Equation~(\ref{def:estimator_hat_beta})) will be denoted by $\check{\beta}_{n,n'}$, or simply $\check{\beta}$. Hereafter, we still set
$\Sc_n=\{j : \check{\beta}_j \neq 0\}$.
The following theorem is proved in Section \ref{proof:thm:WeakConvLassoAdaptive}.

\begin{thm}[Asymptotic law of the adaptive estimator of $\beta$]
    Under the assumptions of Lemma \ref{lemma:asymptNorm_hatTau}, if
    $\mu_{n, n'}  (n h_{n,n'}^p)^{1/2}  \rightarrow \ell\geq 0 $ and $\mu_{n, n'}  (n h_{n,n'}^p)^{1/2} \nu_n^\delta  \rightarrow \infty $ when $n\rightarrow \infty$, we have
    \begin{equation*}
        (n h_{n,n'}^p)^{1/2} (\check \beta_{n, n'} - \beta^*)_{\Sc}
        \indistrto \u_{\Sc}^{**} :=
        \underset{\u_{\Sc} \in \Rb^{s}}{\arg \min} \,
        \check\FF_{\infty, n'}(\u_{\Sc}),\, \text{where}
    \end{equation*}
    \begin{align*}
        \check\FF_{\infty, n'}(\u_{\Sc})
        := \frac{2}{n'} \sum_{i=1}^{n'} \sum_{j\in \Sc}
        W_{i} \psi_j (\z'_i) u_{j}
        + \frac{1}{n'} \sum_{i=1}^{n'} \Big(
        \sum_{j\in \Sc} \psi_j (\z'_i) u_{j} \Big)^2+\ell \sum_{i\in \Sc} \frac{u_i}{|\beta_i^*|^\delta} \sgn(\beta_i^*) ,
    \end{align*}
    with $\W = (W_1, \dots, W_{n'}) \sim \Nc \big(0, \tilde \HH \big).$
    Moreover, when $\ell=0$, the oracle property is fulfilled: $ \PP\left( \Sc_n= \Sc \right) \underset{n}{\rightarrow} 1$.
    \label{thm:WeakConvLassoAdaptive}
\end{thm}

\subsection{Asymptotic properties of $\hat \beta$ when $n$ and $n'$ jointly tend to $+\infty$}
\label{Asymp_n_nprime}
Now, we consider a framework in which both $n$ and $n'$ are going to the infinity, while the dimensions $p$ and $p'$ stay fixed.
To be specific, $n$ and $n'$ will not be allowed to independently go to the infinity.
In particular, for a given $n$, the other size $n'(n)$ (simply denoted as $n'$) will be constrained, as detailed in the assumptions below.
In this section, we still work conditionally on $\z'_1,\ldots,\z'_{n'},\ldots$. The latter vectors are considered as ``fixed'', inducing a deterministic sequence.
Alternatively, we could consider randomly drawn $\z'_i$ from a given law. The latter case can easily been stated from the results below but its specific statement is left to the reader.


\begin{thm}[Consistency of $\hat \beta_{n,n'}$, jointly in $(n,n')$]
    Assume that Assumptions \ref{assumpt:kernel_integral}-\ref {assumpt:f_XZ_Holder} and~\ref{assumpt:compact_case} are satisfied.
    Assume that $ \sum_{i=1}^{n'} \psibm(\z'_i)\psibm(\z'_i)^T/n'$ converges to a matrix $M_{\psi,\z'}$, as $n' \to \infty$.
    Assume that $\lambda_{n, n'} \to \lambda_{0}$ and $n' \exp ( - A n h^{2p}) \to 0$ for every $A>0$, when $(n,n') \to \infty$.
    Then $\hat \beta_{n,n'} \inprobto \arg \min_{\beta \in \Rb^{p'}} \GG_{\infty,\infty}(\beta),$ as $(n,n') \to \infty,$
    where $\GG_{\infty,\infty}(\beta):= (\beta^* - \beta) M_{\psi,\z'} (\beta^* - \beta)^T + \lambda_{0} |\beta|_1$.
    Moreover, if $\lambda_0 = 0$ and $M_{\psi,\z'}$ is invertible, then $\hat \beta_{n,n'}$ is consistent and tends to the true value~$\beta^*$.
    \label{thm:consistency_hatBeta_n_nprime}
\end{thm}
Proof of this theorem is provided in the Supplementary Material, Section \ref{proof:thm:consistency_hatBeta_n_nprime}.
Note that, since the sequence $(\z'_i)$ is deterministic, we just assume the usual convergence of $ \sum_{i=1}^{n'} \psibm(\z'_i)\psibm(\z'_i)^T/n'$ in $\Rb^{p'{}^2}$.
Moreover, if the ``second subset'' $(\z'_i)_{i=1,\ldots,n'}$ were a random sample (drawn along the law $\PP_{\Z}$), the latter convergence would be understood ``in probability''.
And if $\PP_{\Z}$ satisfies the identifiability condition (Proposition \ref{prop:identifiability_condition}), then $M_{\psi,\z'}$ would be invertible and $\hat \beta_{n,n'}\to \beta^*$ in probability.
Now, we want to go one step further and derive the asymptotic law of the estimator $\hat \beta_{n,n'}$.
\begin{assumpt}
    \begin{enumerate}[(i)]
        \item The support of the kernel $K(\cdot)$ is included into $[-1,1]^p$.
        Moreover, for all $n,n'$ and every $(i, j) \in \{1, \dots, n'\}^2$, $i \neq j$, we have
        $|\z'_{i} - \z'_{j} |_\infty > 2 h_{n,n'}$.
        \item \label{assumpt:rates_n_nprime}
        (a) $n' (n h_{n,n'}^{p+4\alpha} + h_{n,n'}^{2\alpha} + (n h_{n,n'}^p)^{-1} ) \to 0$, (b) $\lambda_{n,n'}  ( n' \, n \,  h_{n,n'}^p)^{1/2} \to 0$, \\
        (c) $ n \,  h_{n,n'}^{p+\alpha}/\ln n' \to \infty$.
        \item \label{assumpt:limit_distrib_Zprime}
        The distribution $\PP_{\z',n'} :=  \sum_{i=1}^{n'} \delta_{\z'_i}/n'$ weakly converges as $n' \to \infty$, to a distribution $\PP_{\z',\infty}$ on $\Rb^p$, with a density $f_{\z', \infty}$ with respect to the $p$-dimensional Lebesgue measure.
        \item The matrix $V_1 := \int \psibm(\z') \psibm(\z')^T f_{\z', \infty}(\z') d\z'$ is non-singular.
        \item \label{assumpt:Lambda_C2}
        $\Lambda(\cdot)$ is two times continuously differentiable. Let $\Tc$ be the range of $\z\mapsto \tau_{1,2|\Z=\z}$, from $\Zc$ towards $[-1,1]$.
        On an open neighborhood of $\Tc$, the second derivative of $\Lambda(\cdot)$ is bounded by a constant $C_{\Lambda''}$.
        %
        %
        %
    \end{enumerate}
    \label{assumpt:asymptNorm_joint}
\end{assumpt}

Part (i) of the latter assumption forbids the design points $(\z'_i)_{i \geq 1}$ from being too close to each other and too fast, with respect to the rate of convergence $(h_{n,n'})$ to $0$. This can be guaranteed by choosing an appropriate design.
For example, if $p=1$ and $\Zc = [0,1]$, choose the dyadic sequence $1/2, 1/4, 3/4, 1/8, 3/8, 5/8, 7/8, \dots$

\mds

Part (ii) can be ensured by first choosing a slowly growing sequence $n'(n)$, and then by choosing $h$ that would tend to $0$ fast enough. Note that a compromise has to be found concerning these two rates. The sequence $\lambda_{n,n'}$ should be chosen at last, so that (b) is satisfied.
Interestingly, it is always possible to choose the asymptotically optimal bandwidth, i.e. $h\propto n^{-1/(2\alpha+p)}$. In this case, we can set $n'=n^a$, with any
$a\in ] 0,2\alpha/(2\alpha + p)[$ and the constraints are satisfied.

\mds

The design points $\z'_i$ are deterministic, similarly to all results in the present paper. For a given $n'$, we can invoke the non-random measure $\PP_{\z',n'} := n'{}^{-1} \sum_{i=1}^{n'} \delta_{z'_i}.$
Equivalently, all results can be seen as given conditionally on the sample $(\z'_i)_{i \geq 1}$. In (iii), we impose the weak convergence of $\PP_{\z',n'}$ to a measure with density w.r.t. the Lebesgue measure. Intuitively, this means we do not want to observe some design points that would  be repeated infinitely often (this would result in a Dirac component in $\PP_{\z',\infty}$).
An optimal choice of the density $f_{\z', \infty}$ is not an easy task. Indeed, even if we knew exactly the true density $f_\Z$, there is no obvious reasons why
we should select the $\z'_i$ along $f_{\Z}$ (at least in the limit). If we want a small asymptotic variance $\tilde V_{as}$ (see below), the distribution of the design should concentrate the $\z'_i$ in the regions where $\Lambda' \big(\tau_{1,2|\Z = \z'} \big)^2$ is small and where $\psibm(\z') \psibm(\z')^T$ is big.

\mds

Part (iv) of the assumption is usual, and ensure that the design is somehow ``asymptotically full rank''.
This matrix $V_1$ will also appear in the asymptotic variance of $\hat\beta_{n,n'}$.

\mds

Part \ref{assumpt:Lambda_C2} allow us to control a remainder term in a Taylor expansion of $\Lambda$. Notice that this technical assumption was not necessary in the previous section, where we used the Delta-method on the vector $(\hat \tau_{1,2|\Z = \z'_i} - \tau_{1,2|\Z = \z'_i})_{i=1, \dots, n'}$. But when the number of terms $n'$ tends to infinity, we have to invoke second derivatives to control remainder terms.

\mds

The proof of the next theorem is provided in Section \ref{proof:thm:weak_conv_doubleAsympt}.
\begin{thm}[Asymptotic law of $\hat \beta_{n,n'}$, jointly in $(n,n')$]
    Under Assumptions~\ref{assumpt:asymptNorm_joint} and~\ref{assumpt:kernel_integral}-\ref{assumpt:f_XZ_Holder},
    we have
    $$ (n n' h_{n,n'}^p)^{1/2} (\hat \beta_{n, n'} - \beta^*) \indistrto \Nc(0, \tilde V_{as}),$$
    where $\tilde V_{as} := V_1^{-1} V_2 V_1^{-1}$, $V_1$ is the matrix defined in Assumption~\ref{assumpt:asymptNorm_joint}(iv), and
    \begin{align*}
        &V_2 := \int K^2 \int (\tilde g(\x_1 , \x_3) \tilde g(\x_2 , \x_3) - \tau_{1,2|\z'_1=\z'_2=\z}) \Lambda' \big(\tau_{1,2|\Z = \z} \big)^2 \psibm(\z) \psibm(\z)^T, \\
        & \hspace{2cm} \times f_{\X|\Z}(\x_1|\Z=\z) f_{\X|\Z}(\x_2|\Z=\z) f_{\X|\Z}(\x_3|\Z=\z) \frac{f_{\z',\infty}(\z)}{f_{\Z}(\z)}
        \, d\x_1 \, d\x_2 \, d\x_3 \, d\z'.
    \end{align*}
    \label{thm:weak_conv_doubleAsympt}
\end{thm}

\section{Simulations}
\label{section:simulations}

\subsection{Numerical complexity}
\label{example:computation_time}

Let us take a short numerical application to compare the complexity of our new estimator with the kernel-based ones.
Assume that the size of our dataset is $n=1.000$, with a fixed small $p$, and $p' = 100$. We want to estimate the conditional Kendall's tau on $m = 10.000$ given points $\z_1, \dots, \z_m$.
Using simple kernel-based estimation, the total number of operations is of the order of $n^2 \times m = 1.000^2 \times 10.000 = 10^{10}$.
On the contrary, using our new parametric estimators, the cost can be decomposed in the following way:
\begin{enumerate}
    \item We choose the design points $\z'_1, \dots, \z'_{n'}$ (say, equi-spaced) with $n' = 100$.
    \item We estimate the kernel-based estimator on these $n'$ points (cost: $n^2 \times n' = 1.000^2 \times 100 = 10^{8}$).
    \item We run the Lasso optimization, which is a convex program, so its computation time is linear in~$n'$ and $p'$ (cost: $n' \times p' = 100 \times 100 = 10^4 $).
    \item Finally, for each $\z_i$, we compute the prediction $\Lambda^{(-1)} \big( \hat \beta^T \z_i \big)$, and let us assume that $s = 50$ (cost: $m \times s = 10.000 \times 50 = 5 \times 10^5$).
\end{enumerate}
Summing up, the computational cost of this realistic experiment is around $10^8$, which is $100$ times faster than the kernel-based estimator. Moreover, each new point $\z_{m+1}$ will result in a marginal supplementary cost of $50$ operations, compared with a marginal cost of $n^2=1.000^2=10^6$ for the kernel-based estimator.
Such a huge difference is due to the fact that we have transformed what was previously available as U-statistic of order $2$ with a $O(n^2)$ computational cost for each prediction, into a linear parametric model with $s$~non-zero parameters, giving a cost of $O(s)$ operations for each prediction.

\subsection{Choice of tuning parameters and estimation of the components of $\beta$}
\label{subsubsection:simus:estimation_components_beta}

Now, we evaluate the numerical performance of our estimates through a simulation study.
In this subsection, we have chosen $n=3000$, $n'=100$ and $p=1$.
The univariate covariate $Z$ follows a uniform distribution between $0$ and $1$.
The marginals $X_1 | Z=z$ and $X_2 |Z=z$ follow some Gaussian distributions $\Nc(z,1)$. The conditional copula of $(X_1, X_2) | Z=z$ belongs to the Gaussian copula family. Therefore, it will be parameterized by its (conditional) Kendall's tau $\tau_{1,2|Z=z}$, and is denoted  by $C_{\tau_{1,2|Z=z}}$.
Obviously, $\tau_{1,2|Z=z}$ is given by Model (\ref{model:lambda_cond_tau_Z}).
The dependence between $X_1$ and $X_2$, given $Z=z$, is specified by $\tau_{1,2|Z=z} := 3 z (1-z) = 3/4 - (3/4)(z - 1/2)^2$.

\mds

We will choose $\Lambda$ as the identity function and the $\z'_i$ as a uniform grid on $[0.01,0.99]$. The values $0$ and $1$ for the $\z'_i$ are excluded to avoid boundaries numerical problems.
As for regressors, we will consider $p' = 12$ functions of $Z$, namely $\psi_1(z) = 1$,
$\psi_{i+1}(z) = 2^{-i} (z-0.5)^i$ for $i=1,\dots, 5$,
$\psi_{5+2i}(z) = \cos(2i \pi z )$ and
$\psi_{6+2i}(z) = \sin(2i \pi z )$ for $i=1,2$,
$\psi_{11}(z) = \1\{z \leq 0.4 \}$,
$\psi_{12}(z) = \1\{z \leq 0.6 \}$.
They cover a mix of polynomial, trigonometric and step-functions.
Then, the true parameter is $\beta^* = (3/4, 0, -3/4, {\bf 0}_9)$, where ${\bf 0}_9$ is the null vector of size $9$.

\mds

Our reference value of the tuning parameter $h$ is given by the usual rule-of-thumb, i.e. $h = \hat \sigma(Z) n^{-1/5}$, where $\hat \sigma$ is the estimated standard deviation of $Z$.
Data-driven choices of the bandwidth $h$ of the first estimator are presented in~\cite{derumigny2018kernelBased}.
Moreover, we designed a cross validation procedure (see Algorithm \ref{algo:cross_validation_lambda}) whose output is a data-driven choice for the tuning parameter $\hat \lambda^{cv}$.
Finally, we perform the convex optimization of the Lasso criterion using the R package \texttt{glmnet} by~\cite{friedman2017glmnet}.

\begin{algorithm}[htb]
\label{algo:cross_validation_lambda}
\SetAlgoLined
    Divide the dataset $\Dc = (X_{i,1}, X_{i,2}, \Z_i)_{i=1, \dots, n}$ into $N$ disjoint blocks $\Dc_1, \dots, \Dc_N$ \;
    \ForEach{$\lambda$} {
        \For{$k \leftarrow 1$ \KwTo $N$}{
            Estimate the conditional Kendall's taus $\big(\hat \tau_{1,2 |\Z = \z'_i}^{(k)}\big)_{i=1, \dots, n'}$ on the dataset~$\Dc_k$ \;
            Estimate $\hat\beta^{(-k)}$ by Equation (\ref{def:estimator_hat_beta}) on the dataset $\Dc \backslash \Dc_k$ using the tuning parameter $\lambda$ \;
            Compute $Err_k(\lambda) := \sum_{i=1, \dots, n'} \left( \hat \tau_{1,2 |\Z = \z'_i}^{(k)} - \psibm(\z'_i)^T \hat\beta^{(-k)} \right)^2$ \;
        }
    }
    Return $\hat \lambda^{cv} := \arg \min_\lambda \sum_k Err_k(\lambda)$.
\caption{Cross-validation algorithm for choosing $\lambda$.}
\end{algorithm}

\mds

In our simulations, we observed that the estimation of $\hat \beta$ is not very satisfying if the family of function $\psi_i$ is far too large. Indeed, our model will ``learn the noise'' produced by the kernel estimation, and there will be ``overfitting'' in the sense that the function $\Lambda^{(-1)} \big( \psibm (\cdot)^T \hat \beta \big)$ will be very close to $\hat \tau_{1,2|\Z = \cdot}$, but not to the target $\tau_{1,2|\Z = \cdot}$.
Therefore, we have to find a compromise between misspecification (to choose a family of $\psi_i$ that is not rich enough), and over-fitting (to choose a family of $\psi_i$ that is too rich).

\mds

\begin{figure}[tb]
    \centering
    \includegraphics[height = 8cm]{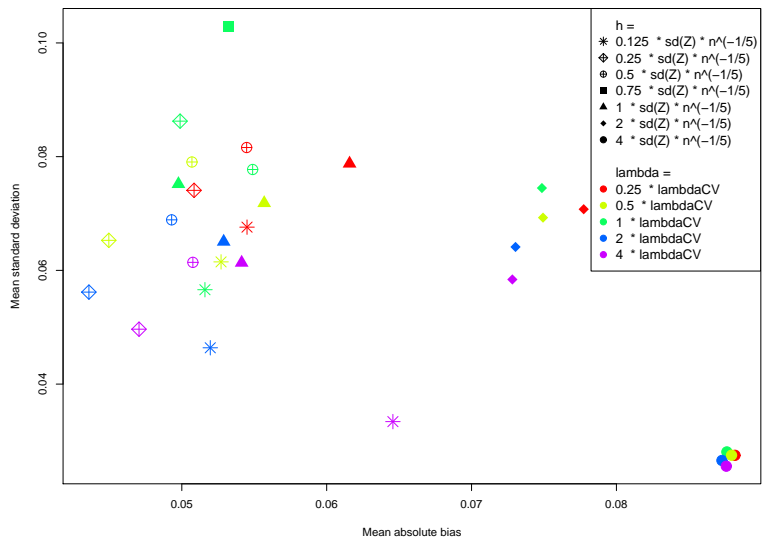}
    \caption{Mean absolute bias $ \sum_{i=1}^{12} | \EE[\hat \beta_i] - \beta_i^* |/12$ and mean standard deviation $ \sum_{i=1}^{12} \sigma (\hat \beta_i)/12$, for different data-driven choices of the tuning parameters $h$ and $\lambda$.}
    \label{fig:plot_bias_sd_fct_h_lambda}
\end{figure}

We have led 100 simulations for couples of tuning parameters $(\lambda,h)$, where $\lambda \propto \hat \lambda^{cv}$, and $h \propto \hat \sigma(Z) n^{-1/5}$. The results in term of empirical bias and standard deviation of $\hat\beta$ are displayed in Figure \ref{fig:plot_bias_sd_fct_h_lambda}.
Empirically, we find the smallest $h$ tend to perform better than the largest ones.
The influence of the tuning parameter $\lambda$ (around reasonable values) is less clear.
Finally, we selected $h = 0.25 \hat \sigma(Z) n^{-1/5}$ and $\lambda = 2 \hat \lambda^{cv}$.
With the latter choice, the coefficient by coefficient results are provided in Table \ref{tab:result_estimation_model1}.
The empirical results are relatively satisfying, despite a small amount of over-fitting.
In particular, the estimation procedure is able to identify the non-zero coefficients almost systematically.
To give a complete picture, for one particular simulated sample, we show the results of the estimation procedure, as displayed in Figures \ref{fig:plot_poly2_coeff_lambda} and \ref{fig:plot_poly2_comparison} in the supplementary material ``Supplementary figures on a simulated sample''.

\begin{table}[tb]
    \centering
    \resizebox{\textwidth}{!}{%
    \begin{tabular}{|c|cc cc cc cc cc cc|}
        \hline
        & $\hat \beta_1$ & $\hat \beta_2$ & $\hat \beta_3$ &
        $\hat \beta_4$ & $\hat \beta_5$ & $\hat \beta_6$ &
        $\hat \beta_7$ & $\hat \beta_8$ & $\hat \beta_9$ &
        $\hat \beta_{10}$ & $\hat \beta_{11}$ & $\hat \beta_{12}$ \\
        \hline
        True value  & 0.75 & 0 & -0.75 & 0 & 0 & 0 & 0 & 0 & 0 & 0 & 0 & 0 \\
        Bias  & -0.13 & 3.6e-05 & 0.26 & 0.0033 & -0.045 & -0.0051 & -0.011 & -2e-04 & -3.2e-05 & 0.073 & -0.0013 & 0.00021 \\
        Std. dev. & 0.15 & 0.00041 & 0.18 & 0.035 & 0.078 & 0.041 & 0.022 & 0.0051 & 0.00037 & 0.15 & 0.007 & 0.0041 \\
        Prob. & 1 & 0.015 & 0.96 & 0.015 & 0.4 & 0.069 & 0.36 & 0.076 & 0.0076 & 0.33 & 0.038 & 0.023 \\
        \hline
    \end{tabular}
    }
    \caption{Estimated bias, standard deviation and probability of being non-null for each estimated component of $\beta$ ($h = 0.25 \,\hat \sigma(Z) n^{-1/5}$ and $\lambda = 2 \hat \lambda^{cv}$).}
    \label{tab:result_estimation_model1}
\end{table}


\subsection{Comparison between parametric and nonparametric estimators of the conditional Kendall's tau}
\label{subsubsection:comparison_two_estimators}

We will now compare our estimator of the conditional Kendall's tau, i.e. $\z \mapsto \Lambda^{(-1)} \big( \psibm (\z)^T \hat \beta \big)$ with the kernel-based estimator, i.e. the first-step estimator.
For this, we will consider six different settings:
\begin{enumerate}
    \item as previously, a Gaussian copula parameterized by its conditional Kendall's tau, given by
    $\tau_{1,2|Z=z} := 3 z (1-z) = 3/4 - (3/4)(z - 1/2)^2$ (well-specified model) ;
    \item a badly-specified model, with a Frank copula whose parameter is given by $\theta(z) = \tan(\pi z / 2)$. Note that the parameter $\theta$ of the Frank family belongs to $\Rb \backslash \{0\}$ and that its Kendall's tau is not written in terms of standard functions of its parameter $\theta$, see~\cite[p.171]{nelsen2007introduction} ;
    \item an intermediate model with a Frank copula calibrated to have the same conditional Kendall's tau as in the first setting ;
    \item another intermediate model with a Gaussian copula calibrated to have the same conditional Kendall's tau as in the second setting ;
    \item a Gaussian copula with a conditional Kendall's tau constant equal to $0.5$ ;
    \item a Frank copula with a conditional Kendall's tau constant equal to $0.5$.
\end{enumerate}
This setting will allows to see the effect of good/bad specifications and of changes in terms of copula families.
In Table~\ref{tab:result_comparisonEstimators}, for each setting, we provide five numerical measures of performance of a given estimator:
\begin{itemize}
    \item the integrated bias: $IBias := \int_z \big( \EE[\hat \tau_{1,2| Z=z}] - \tau_{1,2| Z=z} \big) dz$ ;
    \item the integrated variance: $IVar := \int_z \EE \Big[  \big( \hat \tau_{1,2| Z=z} - \EE[\hat \tau_{1,2| Z=z}] \big)^2 \Big] dz$ ;
    \item the integrated standard deviation: $ISd := \int_z \EE \Big[ \big( \hat \tau_{1,2| Z=z} - \EE[\hat \tau_{1,2| Z=z}] \big)^2 \Big]^{1/2} dz$ ;
    \item the integrated mean square-error: $IMSE := \int_z \EE \Big[ \big( \hat \tau_{1,2| Z=z} - \tau_{1,2| Z=z} \big)^2 \Big] dz$ ;
    \item the CPU time used for the computation.
\end{itemize}
Note that integrals have been approximately computed using a discrete grid $\{ 0.0005 \times i, i = 0, \dots, 2000 \}$.
Globally, in terms of IMSE, the parametric estimator of $\tau_{1,2|\z}$ is doing a better work than a kernel estimator almost systematically (with the single exception of setting 3) and not only in terms of computation time. Surprisingly, even under mis-specification, this conclusion applies whatever the sample size. The differences are particularly striking when the conditional Kendall's tau is a constant function (i.e. under the simplifying assumption).

\begin{table}[p]
    \centering
    \resizebox{\textwidth}{!}{%
    \begin{tabular}{c | cccccc | ccccccc}
        \multicolumn{1}{c}{ } &
        \multicolumn{6}{c}{Kernel-based estimator} &
        \multicolumn{6}{c}{Two-step estimator with $n'=100$ points}  \\
        \cmidrule(lr){2-7} \cmidrule(lr){8-13}
        \multicolumn{1}{c}{Setting} & 1 & 2 & 3 & 4 & 5 & \multicolumn{1}{c}{6}
        & 1 & 2 & 3 & 4 & 5 & 6 \\
    \hline \multicolumn{1}{c}{ } & \multicolumn{12}{c}{ $ n = 500 $} \\ 
    \hline $IBias$ & -29.3 & -14.9 & -31.5 & -6.35 & -32.2 & -29.9 & -23.9 & -19.5 &  -26 & -10.5 & -31.6 & -29.9 \\ 
    $IVar$ & 17.4 & 26.4 & 16.9 & 26.2 & 18.5 & 16.8 &   27 & 17.1 &   28 & 16.8 &  1.9 & 1.65 \\ 
    $ISd$ &  123 &  158 &  120 &  157 &  132 &  126 & 43.3 & 62.5 & 43.8 & 56.4 & 29.7 & 26.6 \\ 
    $IMSE$ & 17.4 & 26.5 & 16.9 & 26.4 & 18.5 & 16.8 &   27 & 17.1 &   28 & 16.9 & 1.91 & 1.65 \\ 
    CPU time (s) & 4.63 & 5.83 & 4.62 & 4.85 & 4.74 &  4.9 & 1.47 & 1.72 & 1.42 & 1.45 & 1.52 & 1.54 \\ 
    \hline \multicolumn{1}{c}{ } & \multicolumn{12}{c}{ $ n = 1000 $} \\ 
    \hline $IBias$ & -16.6 & -11.6 & -15.8 & -2.97 & -16.6 & -17.7 & -12.6 & -12.3 & -12.3 & -5.42 & -16.6 & -17.6 \\ 
    $IVar$ & 8.92 & 17.3 & 8.23 & 13.8 & 8.82 & 8.52 & 8.06 & 7.59 & 9.03 & 6.31 & 0.622 & 0.659 \\ 
    $ISd$ & 89.2 &  116 & 84.5 &  115 & 92.2 & 90.5 & 30.2 & 47.8 & 35.5 & 43.1 & 18.2 & 18.6 \\ 
    $IMSE$ & 9.01 & 17.4 & 8.31 &   14 & 8.88 & 8.57 & 8.07 & 7.61 & 9.04 & 6.34 & 0.624 & 0.661 \\ 
    CPU time (s) &   13 & 12.5 & 12.8 & 12.3 & 12.3 & 12.7 & 3.44 & 3.58 & 3.73 & 3.59 & 3.63 & 3.68 \\ 
    \hline \multicolumn{1}{c}{ } & \multicolumn{12}{c}{ $ n = 2000 $} \\ 
    \hline $IBias$ & -9.94 & -4.96 &  -10 & -4.47 & -10.7 & -10.5 & -6.99 & -6.55 & -7.27 & -5.81 & -10.6 & -10.5 \\ 
    $IVar$ & 4.76 & 7.62 & 4.49 & 7.81 & 4.94 & 4.65 & 3.09 & 2.49 &  3.3 & 2.44 & 0.345 & 0.351 \\ 
    $ISd$ & 65.2 &   85 & 62.6 & 86.4 & 69.4 & 67.3 & 22.7 & 31.4 & 22.3 & 32.3 & 14.7 & 15.2 \\ 
    $IMSE$ & 4.77 & 7.63 &  4.5 & 7.83 & 4.95 & 4.66 & 3.09 & 2.49 &  3.3 & 2.44 & 0.345 & 0.352 \\ 
    CPU time (s) & 67.7 & 68.6 & 67.2 & 73.4 & 72.3 & 59.2 & 15.1 & 15.1 & 15.1 & 16.4 & 17.9 & 14.8 \\ 
    \hline
    \end{tabular}
    }
    \caption{Comparison of the performance between the two estimators. Integrated measures have been multiplied by $10^3$, for readability.}
    \label{tab:result_comparisonEstimators}
\end{table}

\subsection{Comparison with the tests of the simplifying assumption}

Now, under the six previous settings, we compare the test of the simplifying assumption $\Hc_0$ developed in Section~\ref{Asymp_n_fixed_nprime} with some of the bootstrapped-based tests of the latter assumption that has been introduced in~\cite{derumigny2017tests}.
In particular, they propose a nonparametric test, using the statistic $\Tc^0_{CvM}$ defined by
\begin{align*}
    \Tc^0_{CvM}
    := \int_{[0,1]^3} \Big( \hat C_{1,2| Z = \hat F_Z^{-1}(u_3)}(u_1, u_2)
    - \hat C_{s,1,2 | Z}(u_1, u_2) \Big)^2 du_1 du_2 du_3,
\end{align*}
where $\hat C_{1,2| Z= z}$ is a kernel-based nonparametric estimator of the conditional copula of $(X_1,X_2)|Z=z$ and $\hat C_{s,1,2 | Z}(u_1, u_2) := n^{-1} \sum_{i=1}^n \hat C_{1,2 | Z = Z_i}(u_1, u_2).$
We will also invoke their parametric test statistic
\begin{align*}
    \Tc_2^c := \int_0^1 \Big(
    \hat \theta \big( \hat F_Z^{-1}(u) \big) - \hat \theta \Big)^2 du,
\end{align*}
where $\hat\theta(z)$ estimates the parameter of the Gaussian (resp. Frank) copula given $Z = z$, assuming we know the right family of conditional copula, and $\hat\theta$ consistently estimates the parameter of the corresponding simplified copula (under the null). Moreover, $\hat F_Z^{-1}$ denotes the empirical quantile function that is associated to the $Z$-sample.
The latter test statistics depends on an a priori chosen parametric copula family. To evaluate the risk of mis-specification, we also include in our table the parametric test $\Tc_2^c$ assuming that the data come from a Clayton copula, whereas the true copula is Gaussian or Frank.
For these three tests, p-values are computed by the usual nonparametric bootstrap, with $100$ resampling: see Table~\ref{tab:result_comparisonTests}. 
Globally, the test based on $\Wc_n$ performs very well under all settings, compared to the alternative nonparametric test. It is only beaten by $\Tc_2^c$ that is obtained by choosing the right copula family, a not very realistic situation. When it is not the case, $\Wc_n$ does a better work.  

\begin{table}[htb]
    \centering
    \begin{tabular}{l cccccc}
        &\multicolumn{4}{c}{Not under $\Hc_0$}
        & \multicolumn{2}{c}{Under $\Hc_0$} \\
        \cmidrule(lr){2-5} \cmidrule(lr){6-7}
        & 1 & 2 & 3 & 4 & 5 & 6 \\
        \hline 
        $\Wc_n$ & 88.7 & 99.8 & 87.3 &  100 &   12 & 12.1 \\ 
        $\Tc^0_{CvM}$ & 59.5 &   52 & 64.7 & 37.5 &    0 &    0 \\ 
        $\Tc_2^c$ &  100 &  100 &  100 &  100 &  0.2 &  2.6 \\ 
        $\Tc_2^c$ (Clayton) &   68 &   13 &  100 &  100 &  1.8 &  1.8 \\
        \hline
    \end{tabular}
    \caption{Comparison of the performance between different tests of the simplifying assumption under the six settings of Section~\ref{subsubsection:comparison_two_estimators}, with $n=500$.}
    \label{tab:result_comparisonTests}
\end{table}

\subsection{Dimension $2$ and choice of $\psibm$}

In this section, we will fix the sample size $n=3000$ and the dimension $p=2$.
The random vector $\Z$ will follow a uniform distribution on $[0,1]^2$, $X_1 | \Z = \z \sim \Nc(0, z_1)$, $X_2 | \Z = \z \sim \Nc(0, z_1)$. Given $\Z = \z$, the conditional copula of $X_1$ and $X_2$ is Gaussian.
We consider three different choices for the functional form of its conditional Kendall's tau :

Setting 1. $\tau_{1,2|\Z=\z} = (3/4) \times (z_1 - z_2)$ ;

Setting 2. $\tau_{1,2|\Z=\z} = (4/8) \times \cos(2 \pi z_1)
+ (2/8) \times \sin(2 \pi z_2)$ ;

Setting 3. $\tau_{1,2|\Z=\z} = (3/4) \times \tanh(z_1 / z_2)$,

\noindent
where $\z = (z_1, z_2)$.
We try different choices of dictionaries $\psibm$. For convenience, define $p_0(x) := 1$, $p_i(x) := 2^{-i} (x - 0.5)^i$, $trig_0(x) := 1$,
and $trig_i(x) := \big(\cos(2 i \pi x), \sin(2 i \pi x) \big)$, for $x \in \Rb$ and $i \in \NN^*$.
We will use the notation $(g_1, g_2) \otimes (g_3, g_4) := (g_1 g_3, g_1 g_4, g_2 g_3, g_2 g_4)$.
We are interested in the following functions $\psibm$, that are defined for every $\z \in \Rb^p$ by
\begingroup \allowdisplaybreaks
\begin{align*}
    \psibm^{(1)}(\z)
    &:= \Big( 1, \big( p_i(z_1) \big)_{i=1, \dots, 5}, 
    \big( p_i(z_2) \big)_{i=1, \dots, 5} \Big)
    = \Big( p_i(z_1) \times p_j(z_2)
    \Big)_{\min(i,j) = 0, \, \max(i,j) \leq 5} 
    \in \Rb^{11}, \\
    \psibm^{(2)}(\z)
    &:= \Big( p_i(z_1) \times p_j(z_2)
    \Big)_{\min(i,j) \leq 1, \, \max(i,j) \leq 5} 
    \in \Rb^{20}, \\
    \psibm^{(3)}(\z)
    &:= \Big( p_i(z_1) \times p_j(z_2)
    \Big)_{\min(i,j) \leq 2, \, \max(i,j) \leq 5}
    \in \Rb^{27}, \\
    \psibm^{(4)}(\z)
    &:= \Big( p_i(z_1) \times p_j(z_2)
    \Big)_{\max(i,j) \leq 5}
    \in \Rb^{36}, \\
    \psibm^{(5)}(\z)
    &:= \Big( 1, \big( trig_i(z_1)
    \big)_{i=1, \dots, 5}, 
    \big( trig_i(z_2) \big)_{i=1, \dots, 5} \Big) 
    \in \Rb^{21}, \\
    \psibm^{(6)}(\z)
    &:= \Big( trig_i(z_1) \otimes trig_j(z_2)
    \Big)_{\min(i,j) \leq 1, \, \max(i,j) \leq 5} 
    \in \Rb^{57}, \\
    \psibm^{(7)}(\z)
    &:= \Big( trig_i(z_1) \otimes trig_j(z_2)
    \Big)_{\min(i,j) \leq 2, \, \max(i,j) \leq 5} 
    \in \Rb^{85}, \\
    \psibm^{(8)}(\z)
    &:= \Big( trig_i(z_1) \otimes trig_j(z_2)
    \Big)_{\max(i,j) \leq 5}
    \in \Rb^{121}, \\
    \psibm^{(9)}(\z)
    &:= \Big( \psibm^{(1)}(\z), \psibm^{(5)}(\z) \Big)
    \in \Rb^{31}, \hspace{1cm}
    \psibm^{(10)}(\z)
    := \Big( \psibm^{(2)}(\z), \psibm^{(6)}(\z) \Big)
    \in \Rb^{76}, \\
    \psibm^{(11)}(\z)
    &:= \Big( \psibm^{(3)}(\z), \psibm^{(7)}(\z) \Big)
    \in \Rb^{137}, \hspace{1cm}
    \psibm^{(12)}(\z)
    := \Big( \psibm^{(4)}(\z), \psibm^{(8)}(\z) \Big)
    \in \Rb^{156},
\end{align*}
\endgroup
where in the last 4 dictionaries, we count the function constant to $1$ only once.
We choose $n' = 400$ and the design points $\z'_i$ are chosen as an equispaced grid on $[0.1, 0.9]^2$. We consider similar measures of performance for our estimators as in Section~\ref{subsubsection:comparison_two_estimators}. The only difference is that the integration in $\z$ is now done on the unit square $[0,1]^2$. In practice, integrals are discretized, and estimated by a sum over the points
$\{ (0.01 \times i, 0.01 \times j), 0 \leq i,j \leq 100 \}$.
Results are displayed in the following Table~\ref{tab:result_comparison_psibm}.

\begin{table}[htb]
    \centering
    \resizebox{\textwidth}{!}{%
    \begin{tabular}{c | cccc | cccc | cccc |}
        \multicolumn{1}{c}{ } &
        \multicolumn{4}{c}{Setting 1} &
        \multicolumn{4}{c}{Setting 2} &
        \multicolumn{4}{c}{Setting 3} \\
        \cmidrule(lr){2-5} \cmidrule(lr){6-9} \cmidrule(lr){10-13}
        & IBias & ISd & IMSE & Time
        & IBias & ISd & IMSE & Time
        & IBias & ISd & IMSE & Time \\
        \hline
        $\psibm^{(1)}$ & 0.577 & 19.4 & 1.44 & 6.82
        & -0.632 &   24 &  1.4 & 6.75
        & -7.71 &   17 & 6.79 & 6.67 \\
        $\psibm^{(2)}$ & 0.309 & 18.9 & 1.43 & 6.77 & -0.166 & 23.7 & 1.35 & 6.66 & -7.57 & 16.9 &  6.8 & 6.66 \\
        $\psibm^{(3)}$ & 0.728 & 19.9 & 1.63 & 6.77 & -0.36 & 27.1 &  1.9 & 6.67 & -7.63 & 23.7 & 3.45 & 7.06 \\
        $\psibm^{(4)}$ & 0.513 & 18.9 & 1.81 & 6.77 & -0.245 & 26.5 & 2.22 & 6.68 & -7.29 &   25 & 2.06 & 7.52 \\
        \hline
        $\psibm^{(5)}$ &  1.5 & 25.7 & 15.7 & 6.77 & 0.0616 &   15 & 2.67 & 6.66 & -8.38 & 21.6 & 14.9 & 7.51 \\
        $\psibm^{(6)}$ & 1.64 &   26 & 15.7 & 6.79 & 0.269 &   15 & 2.61 & 6.66 & -8.23 & 21.9 & 14.9 & 7.52 \\
        $\psibm^{(7)}$ & 0.311 & 26.1 &   17 & 6.79 & 0.0167 &   15 & 3.14 & 6.69 & -7.33 & 23.1 & 15.1 & 7.26 \\
        $\psibm^{(8)}$ &  1.2 &   26 & 17.3 & 6.88 & -0.113 & 14.6 & 3.15 &  6.7 & -7.6 & 22.9 & 15.3 &  7.2 \\
        \hline
        $\psibm^{(9)}$ & 0.596 & 17.7 & 2.05 & 6.79 & 0.492 & 15.8 & 2.72 & 6.67 & -7.93 & 16.3 & 7.04 & 7.19 \\
        $\psibm^{(10)}$ & -0.0921 &   18 & 2.08 & 6.77 & -0.493 & 16.6 & 2.75 & 6.66 & -7.65 & 16.7 & 6.94 & 7.19 \\
        $\psibm^{(11)}$ & 0.529 & 17.3 & 2.57 & 6.83 & -0.165 & 15.8 & 3.08 &  6.7 & -6.87 &   23 & 4.76 & 7.21 \\
        $\psibm^{(12)}$ &  0.5 & 16.9 & 2.64 & 6.92 & -0.078 & 16.4 & 3.24 & 6.76 & -7.07 & 25.5 & 4.43 & 7.54 \\
        \hline
    \end{tabular}
    }
    \caption{Comparison of the estimation using different $\psibm$ families. All integrated measures have been multiplied by $1000$. Computation time is given in seconds.}
    \label{tab:result_comparison_psibm}
\end{table}

We note that the size of the family $\psibm$ seems to have a tiny influence on the computation time, which lies always between 6 and 8 seconds. In all settings, polynomial families ($\psibm^{(1)}$ to $\psibm^{(4)}$) give the best $IMSE$, even when the true function is trigonometric (Setting 2) or under misspecification (Setting 3). Nevertheless, using trigonometric functions can help to reduce the integrated biais and standard deviation. Indeed, in Setting 2, trigonometric families ($\psibm^{(5)}$ to $\psibm^{(8)}$) do a fair job according to these two measures of performance. Similarly, in Setting 3, mixed families ($\psibm^{(9)}$ to $\psibm^{(12)}$) achieve an acceptable performance. In Settings 1 and 2, they often yield improvement other a msispecified family, especially in terms of integrated standard deviation.

\mds

Comparisons between three indicators $IMSE$, $IBias$ and $ISd$ may be surprising at first sight, but there is no direct link between their values. Indeed, for every point $\z$, $MSE(\z) = Bias(\z)^2 + Sd(\z)^2$, while $IMSE = \int MSE(\z) d\z$,  $IBias = \int Bias(\z) d\z$ and $ISd = \int Sd(\z) d\z$. Therefore, a procedure that minimize both $Ibias$ and $ISd$ still may not minimize $IMSE$, and conversely. This is due to the non-linearity of the square function, combined with the integration.

\section{Real data application}
\label{section:real_data}

Now, we apply the model given by~(\ref{model:lambda_cond_tau_Z}) to a real dataset. From the website of the World Factbook of the Central Intelligence Agency, we have collected data of male and female life expectancy and GDP per capita for $n=206$ countries in the world.
We seek to analyze the dependence between male and female life expectancies conditionally on the GDP per capita, i.e. given the explanatory variable $Z = \log_{10}(GDP/capita)$.
This dataset and these variables are similar as those in the first example studied in~\cite{gijbels2011conditional}.

\mds

\begin{figure}[p]
    {\centering
    \includegraphics[height = 7.5cm]{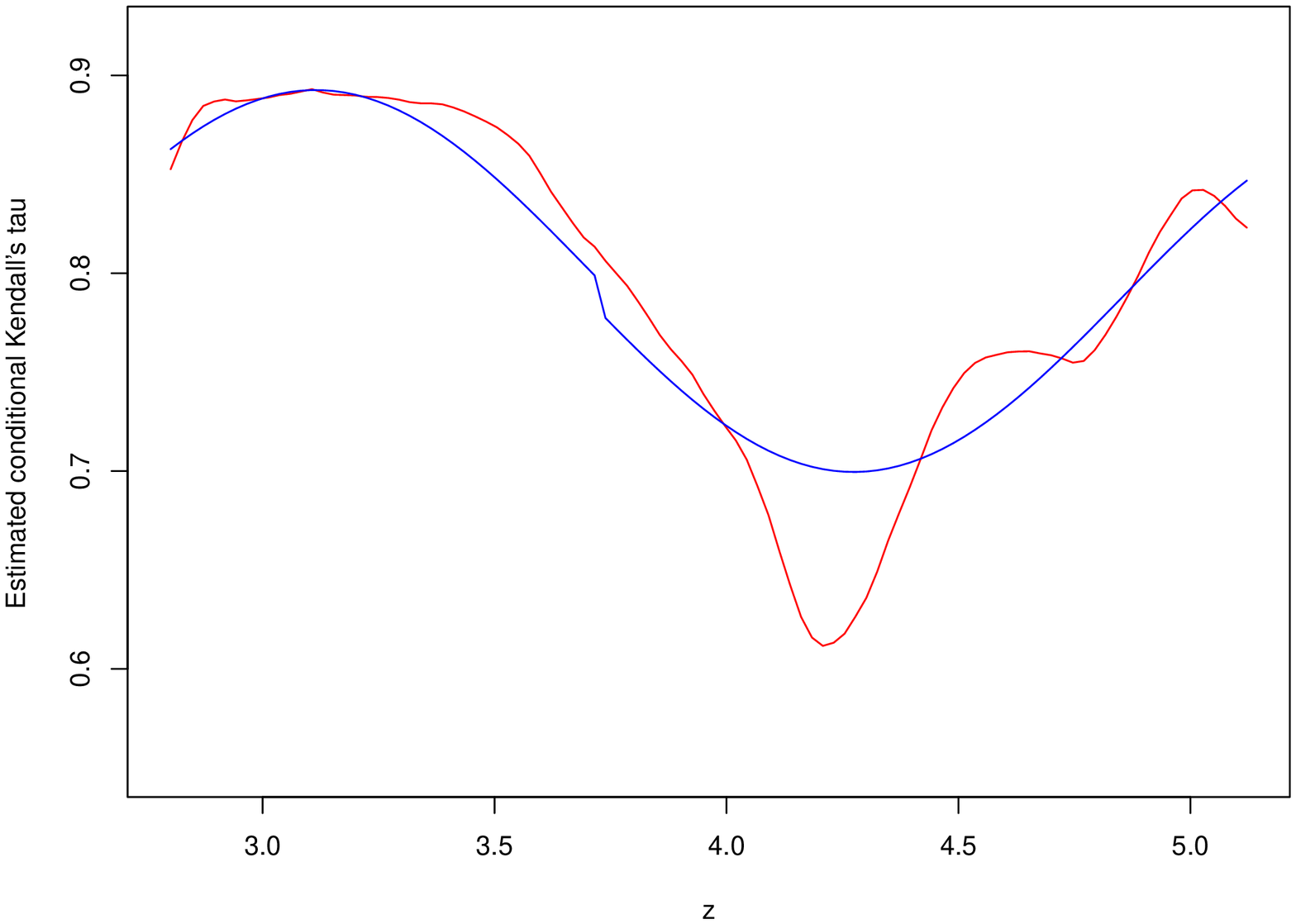}
    \vspace{-0.6cm}
    \caption{Estimated conditional Kendall's tau $\hat \tau_{1,2|\Z = \z}$ (red curve), and prediction $\Lambda^{(-1)} \big( \psibm (\z)^T \hat \beta \big)$ (blue curve) as a function of $\z$ for the application on real data, where the estimated non-zero coefficients are $\hat \beta_1 = 0.78$, $\hat \beta_7 = -0.043$, $\hat \beta_8 = 0.069$ and $\hat \beta_{11} = 0.020$.}
    \label{fig:plot_factbook_comparison}
    }
    \vspace{0.4cm}
    {
    \centering
    \includegraphics[height = 7.5cm]{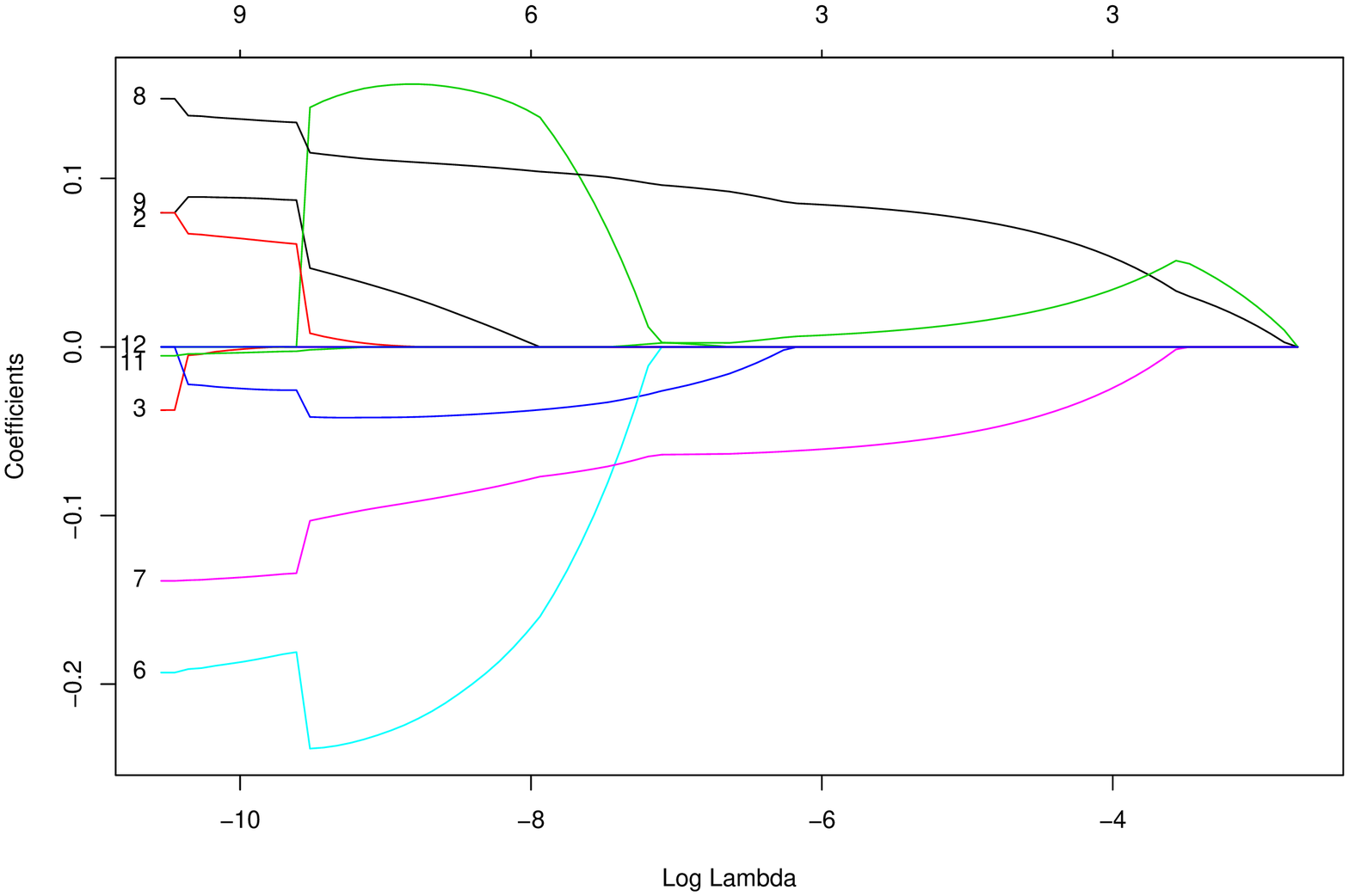}
    \vspace{-0.6cm}
    \caption{Evolution of the estimated non-zero coefficients as a function of the regularization parameter $\lambda$ for the application on real data. All the other non-displayed $\psi_i$ coefficients are zero.}
    \label{fig:plot_factbook_coeff_lambda}
    }
\end{figure}

\mds

We use $n'=100$, $h=2 \sigma(Z)  n^{-1/5}$ and the same family of functions $\psi_i$ as in Section~\ref{subsubsection:simus:estimation_components_beta} above (once composed with a linear transform
to be defined on $[\min(Z), \max(Z)]$). The results are displayed in Figure~\ref{fig:plot_factbook_comparison}.
As expected, the levels of conditional dependence between male and female expectancies are strong overall.
Many poor countries suffer from epidemics, malnutrition or even wars. In such cases, life expectancies of both genders are exposed to the same ``exogenous'' factors, inducing
high Kendall's taus.
Logically, we observe a monotonic decrease of such Kendall's taus when $Z$ is larger, up to $Z\simeq 4.5$, as already noticed by~\cite{gijbels2011conditional}.
Indeed, when countries become richer, more developed and safe, men and women less and less depend on their environment (and on its risks of death, potentially).
Nonetheless, when $Z$ become even larger (the richest countries in the world), conditional dependencies between male and female life expectancies interestingly increase again, because men and women behave similarly in terms of way of life. In particular, they can benefit from the same levels of security and health and are exposed to the same lethal risks.

\section{Supplementary material}

\begin{description}

\item[Proofs of the theoretical results in ``About Kendall's regression'':] In this supplementary material, we detail the proofs for all the results in this paper. We also recall some useful lemmas from~\cite{derumigny2018kernelBased}.

\item[Supplementary figures on a simulated sample:] To give a more precise picture of our estimators, two supplementary figures are given to illustrate their behavior on a typical sample.


\end{description}





\bibliography{biblio}

\begin{thebibliography}{}

\bibitem[Acar et~al., 2012]{AcarGenestNeslehova}
Acar, E., Genest, C., and Ne{$\check{\text{s}}$}lehov\'{a}, J. (2012).
\newblock Beyond simplified pair-copula constructions.
\newblock {\em J. Multivariate Anal.}, 110:74--90.

\bibitem[Bellec et~al., 2016]{bellec2016slope}
Bellec, P.~C., Lecu{\'e}, G., and Tsybakov, A.~B. (2016).
\newblock Slope meets lasso: improved oracle bounds and optimality.
\newblock {\em ArXiv:1605.08651}.

\bibitem[Bickel et~al., 2009]{bickel2009simultaneous}
Bickel, P.~J., Ritov, Y., and Tsybakov, A.~B. (2009).
\newblock Simultaneous analysis of lasso and dantzig selector.
\newblock {\em Ann. Statist.}, 37(4):1705--1732.

\bibitem[Derumigny and Fermanian, 2017]{derumigny2017tests}
Derumigny, A. and Fermanian, J.-D. (2017).
\newblock About tests of the “simplifying” assumption for conditional
  copulas.
\newblock {\em Depend. Model.}, 5(1):154--197.

\bibitem[Derumigny and Fermanian, 2018a]{derumigny2018kernelBased}
Derumigny, A. and Fermanian, J.-D. (2018a).
\newblock {A}bout kernel-based estimation of the conditional kendall's tau:
  finite-distance bounds and asymptotic behavior.
\newblock {\em arXiv preprint arXiv:1810.06234}.

\bibitem[Derumigny and Fermanian, 2018b]{derumigny2018classification}
Derumigny, A. and Fermanian, J.-D. (2018b).
\newblock A classification point-of-view about conditional {K}endall's tau.
\newblock {\em arXiv preprint arXiv:1806.09048}.

\bibitem[Fan and Li, 2001]{fanli2001}
Fan, J. and Li, R. (2001).
\newblock Variable selection via nonconcave penalized likelihood and its oracle
  properties.
\newblock {\em J. Amer. Statist. Assoc.}, 96(456):1348--1360.

\bibitem[Friedman et~al., 2017]{friedman2017glmnet}
Friedman, J., Hastie, T., Tibshirani, R., and Simon, N. (2017).
\newblock glmnet: Lasso and elastic-net regularized generalized linear models.
  {R} package version 2.0--2.

\bibitem[Gijbels et~al., 2011]{gijbels2011conditional}
Gijbels, I., Veraverbeke, N., and Omelka, M. (2011).
\newblock Conditional copulas, association measures and their applications.
\newblock {\em Comput. Statist. Data Anal.}, 55(5):1919--1932.

\bibitem[Hob{\ae}k~Haff et~al., 2010]{HobaekAasFrigessi}
Hob{\ae}k~Haff, I., Aas, K., and Frigessi, A. (2010).
\newblock On the simplified pair-copula construction–simply useful or too
  simplistic?
\newblock {\em J. Multivariate Anal.}, 101:1296--1310.

\bibitem[Kato, 2009]{kato2009asymptotics}
Kato, K. (2009).
\newblock Asymptotics for argmin processes: Convexity arguments.
\newblock {\em J. Multivariate Anal.}, 100(8):1816--1829.

\bibitem[Kurz and Spanhel, 2017]{kurz2017testing}
Kurz, M.~S. and Spanhel, F. (2017).
\newblock Testing the simplifying assumption in high-dimensional vine copulas.
\newblock {\em arXiv preprint arXiv:1706.02338}.

\bibitem[Nelsen, 2007]{nelsen2007introduction}
Nelsen, R.~B. (2007).
\newblock {\em An introduction to copulas}.
\newblock Springer Science \& Business Media.

\bibitem[Veraverbeke et~al., 2011]{veraverbeke2011ScandinJ}
Veraverbeke, N., Omelka, M., and Gijbels, I. (2011).
\newblock Estimation of a conditional copula and association measures.
\newblock {\em Scand. J. Stat.}, 38(4):766--780.

\bibitem[Zou, 2006]{zou2006}
Zou, H. (2006).
\newblock The adaptive lasso and its oracle properties.
\newblock {\em J. Amer. Statist. Assoc.}, 101(476):1418--1429.

\end{thebibliography}


\appendix


\pagebreak
\setcounter{equation}{0}
\setcounter{figure}{0}
\setcounter{table}{0}
\setcounter{page}{1}
\makeatletter
\renewcommand\appendixname{Supplement}
\renewcommand{\theequation}{S\arabic{equation}}
\setcounter{section}{0}
\setcounter{subsection}{0}

\if1\blind
{
\begin{center}
    {\LARGE\bf Proofs of the theoretical results in ``About Kendall's regression''}
    
    \bigskip
    
  {\large Alexis Derumigny$^1$
    and
    Jean-David Fermanian\footnote{
    CREST-ENSAE, 5, avenue Henry Le Chatelier,
    91764 Palaiseau cedex, France. \\
    Email adresses: alexis.derumigny@ensae.fr, jean-david.fermanian@ensae.fr. \\
    This research has been supported by the Labex Ecodec.
    }}
    
    \end{center}
} \fi

\if0\blind
{
  \begin{center}
    {\LARGE\bf Proofs of the theoretical results in ``About Kendall's regression''}
    \end{center}
  \medskip
} \fi

\section{Proofs of finite-distance results for $\hat \beta$}

In this section, we will use the notation $\u := \hat \beta - \beta^*$ and $\xi=[\xi_{i,n}]_{i=1,\ldots,n'}$, $\xi_{i,n}=Y_i -(\Zb' \beta)_i$.

\subsection{Technical lemmas}

\begin{lemma}
    We have $||\Zb'\u||_{n'}^2 \leq \lambda |\u|_1 + \dfrac{1}{n'}
    \left \langle \xi \; , \; \Zb' \u \right \rangle $.
    \label{lemma:beta_subdifferential}
\end{lemma}
{\it Proof :}
As $\hat \beta$ is optimal, through the Karush-Kuhn-Tucker conditions, we have
$(1/n') \Zb'{}^T(\Y-\Zb'\hat \beta) \in \partial \big( \lambda |\hat \beta|_1 \big) ,$
where $\partial \big( \lambda |\hat \beta|_1 \big)$ is the subdifferential of the norm $\lambda |\cdot|_1$ evaluated at $\hat \beta$.
The dual norm of $|\cdot|_1$ is $|\cdot|_\infty$, so there exists $\v$ such that $|\v|_\infty \leq 1$ and
$(1/n') \Zb'{}^T(\Y-\Zb'\hat \beta) + \lambda \v = 0.$
We deduce successively
$ \Zb'{}^T \Zb'(\beta^*-\hat \beta)/n' +
 \Zb'{}^T \xi /n' + \lambda \v = 0$,
\begin{gather*}
    %
    %
    \dfrac{1}{n'} |\Zb'(\beta^*-\hat \beta)|_2^2 + \dfrac{1}{n'}
    (\beta^*-\hat \beta)^T \Zb'{}^T \xi + \lambda (\beta^*-\hat \beta)^T \v = 0,\; \text{and finally} \\
    ||\Zb'(\beta^*-\hat \beta)||_{n'}^2
    \leq \dfrac{1}{n'} \left \langle \Zb'(\hat \beta - \beta^*)  \; , \; \xi \right \rangle
    + \lambda |\beta^*-\hat \beta|_1.\;\; \Box
\end{gather*}

\begin{lemma}
    We have $|\u_{\Sc^C}|_1 \leq |\u_{\Sc}|_1 + \dfrac{2}{\lambda n'}
    \left \langle \xi \; , \; \Zb' \u \right \rangle $.
    \label{lemma:first_order_consequence}
\end{lemma}
{\it Proof :}
By definition, $\hat \beta$ is a minimizer of
$||\Y-\Zb'\beta||_{n'}^2 + \lambda |\beta|_1$. Therefore, we have
\begin{equation*}
    ||\Y-\Zb'\hat \beta||_{n'}^2 + \lambda |\hat \beta|_1
    \leq ||\Y-\Zb'\beta^*||_{n'}^2 + \lambda |\beta^*|_1.
\end{equation*}
After some algebra, we derive
$||\Y-\Zb'\hat \beta||_{n'}^2 - ||\Y-\Zb'\beta^*||_{n'}^2
\leq \lambda \big( |(\beta^* - \hat \beta)_{\Sc}|_1
- |(\hat \beta - \beta^*)_{\Sc^C}|_1 \big).$
Moreover, the mapping $\beta \mapsto ||\Y-\Zb'\beta||_{n'}^2$ is convex and its gradient at $\beta^*$ is
$-2 \, \Zb'{}^T (\Y-\Zb'\beta^*)/n'
= -2 \, \Zb'{}^T \xi/n'$. So, we obtain
\begin{equation*}
    ||\Y-\Zb'\hat \beta||_{n'}^2 - ||\Y-\Zb'\beta^*||_{n'}^2
    \geq \dfrac{-2}{n'}
    \left \langle \Zb'{}^T \xi \; , \; \hat \beta - \beta^* \right \rangle.
\end{equation*}
Combining the two previous equations, we get
$$(-2)    \left \langle \Zb'{}^T \xi , \hat \beta - \beta^* \right \rangle /n'
    \leq \lambda \big( |(\beta^* - \hat \beta)_{\Sc}|_1
    - |(\hat \beta - \beta^*)_{\Sc^C}|_1 \big).\;\; \Box$$

\begin{lemma}
    Assume that
    $\max_{j=1, \dots, p'} \Big| \frac{1}{n'}
    \sum_{i=1}^{n'} Z'_{i,j} \xi_{i,n} \Big| \leq t,$
    for some $t > 0$, that the assumption $RE(s,3)$ is satisfied, and that the tuning parameter is given by $\lambda = \gamma t$, with $\gamma \geq 4$.
    Then, $||\Zb'(\hat \beta - \beta^*)||_{n'}
    \leq \dfrac{4(\gamma+1) t \sqrt{s}}{\kappa(s,3)}$
    and $|\hat \beta - \beta^*|_q
    \leq \dfrac{4^{2/q}(\gamma+1) t s^{1/q} }{\kappa^2(s,3)},$
    for every $1 \leq q \leq 2$.
    \label{lemma:lasso_bound}
\end{lemma}

{\it Proof :}
Under the first assumption, we have the upper bound
\begin{align*}
    \dfrac{1}{n'} \left \langle \Zb'{}^T \xi \; , \; \u \right \rangle
    \leq |\u|_1 \max_{j=1, \dots, p'} \Big| \frac{1}{n'}
    \sum_{i=1}^{n'} Z'_{i,j} \xi_{i,n} \Big|
    \leq |\u|_1 t.
\end{align*}
We first show that $\u$ belongs to the cone $\big\{ \delta \in \Rb^{p'}: |\delta_{\Sc^C}|_1 \leq 3 |\delta_{\Sc}|_1, Card(\Sc) \leq s \big \}$, so that we will be able to use the $RE(s, 3)$ assumption with $J_0 = \Sc$.
From Lemma \ref{lemma:first_order_consequence},
$ |\u_{\Sc^C}|_1 \leq |\u_{\Sc}|_1 + 2 t |\u|_1 / \lambda$.
With our choice of $\lambda$, we deduce
$ |\u_{\Sc^C}|_1 \leq |\u_{\Sc}|_1 + 2 |\u|_1 / \gamma$.
Using the decomposition $|\u|_1 = |\u_{\Sc^C}|_1 + |\u_{\Sc}|_1$, we get $|\u_{\Sc^C}|_1 \leq |\u_{\Sc}|_1 (\gamma+2)/(\gamma-2) \leq 3 |\u_{\Sc}|_1$.
As a consequence, we have
$$|\u|_1 = |\u_{\Sc^C}|_1 + |\u_{\Sc}|_1
\leq 4 |\u_{\Sc}|_1 \leq 4 \sqrt{s} |\u|_2
\leq 4 \sqrt{s} \, ||\Zb' \u||_{n'}/\kappa(s,3).$$
By Lemma \ref{lemma:beta_subdifferential},
\begin{align*}
    ||\Zb'\u||_{n'}^2
    \leq \lambda |\u|_1 + \dfrac{1}{n'}
    \left \langle \xi \; , \; \Zb' \u \right \rangle
    \leq \lambda |\u|_1 + |\u|_1 t
    \leq |\u|_1 (\gamma+1) t
    \leq \frac{4 \sqrt{s}}{\kappa(s,3)} \, ||\Zb' \u||_{n'}
    (\gamma+1) t
\end{align*}
We can now simplify and we get
\begin{align*}
    ||\Zb'\u||_{n'}
    &\leq \dfrac{4(\gamma+1) t}{\kappa(s,3)} \,
    \sqrt{s}, \quad
    |\u|_2 \leq \dfrac{4(\gamma+1) t}{\kappa^2(s,3)} \,
    \sqrt{s}, \text{ and }
    |\u|_1 \leq \dfrac{16 (\gamma+1) t}{\kappa^2(s,3)} \,
    s.
\end{align*}
Now, we compute a general bound for $|\u|_q$,
with $1 \leq q \leq 2$, using the H\"older norm interpolation inequality:
\begin{align*}
    |\u|_q
    \leq |\u|_1^{2/q-1} |\u|_2^{2-2/q}
    \leq \dfrac{4^{2/q}(\gamma+1) t s^{1/q}}{\kappa^2(s,3)} \cdot\;\; \Box
\end{align*}

\subsection{Proof of Theorem \ref{thm:bound_proba_hat_beta} }
\label{proof:thm:bound_proba_hat_beta}

Using Lemma \ref{lemma:exponential_bound_KendallsTau}, for every $t_1,t_2>0$ such that
$  C_{K, \alpha} h^{\alpha} /  \alpha  ! + t_1 \leq  f_{\Z, min}/2$,
with probability greater than
$1 - 2 n' \exp \Big( - \big(n h^p t_1^2 \big) \, / \, \big(2 f_{\Z, max} \int K^2 + (2/3) C_K t_1 \big)) \Big)
- 2 n' \exp \Big( -  (n-1) h^{2p} t_2^2 f^4_{\Z,min}$ $ / \, \big(4 f_{\Z,max}^2 (\int K^2)^2 + (8/3) C_K^2 f_{\Z, min}^2 t_2 \big) \Big)$, we have
\begin{align*}
    \max_{j=1, \dots, p'} \bigg| \frac{1}{n'}
    \sum_{i=1}^{n'} Z'_{i,j} \xi_{i,n} \bigg|
    &\leq C_{\psi} \max_{i=1, \dots, n'} \big| \xi_{i,n} \big|
    \leq C_{\psi } C_{\Lambda'} \max_{i=1, \dots, n'}
    \big| \hat \tau_{1,2|\Z=\z'_i} - \tau_{1,2|\Z=\z'_i} \big| \\
    &\leq 4 C_{\psi } C_{\Lambda'} \bigg(1 + \frac{16 f_{\Z, max}^2}{f_{\Z, min}^3}
    \Big( \frac{  C_{K, \alpha} h^{\alpha}}
    {\alpha  !} + t_1 \Big) \bigg)
    \bigg( \frac{C_{\X\Z, \alpha}   h^\alpha}
    {f_{\Z,min}^2  \alpha  !} + t_2 \bigg).
\end{align*}
We choose $t_1 := f_{\Z, min}/4$ so that, because of Condition~(\ref{cond:h_1_2}), we get
$  C_{K, \alpha} h^{\alpha} /  \alpha  ! + t_1 \leq f_{\Z, min}/2.$
Now we choose $t_2:= t f_{\Z, min}^2 / \{8 C_{\psi } C_{\Lambda'} (f_{\Z, min}^2 + 8 f_{\Z, max}^2)\}$.
By Condition~(\ref{cond:h_1_2}), $C_{\X\Z, \alpha}   h^\alpha /(f_{\Z,min}^2 \alpha !) \leq t_2,$ so that we have
\begin{align*}
    4 C_{\psi } C_{\Lambda'} \bigg(1 + \frac{8 f_{\Z, max}^2}{f_{\Z, min}^2} \bigg)
    \times \bigg( \frac{C_{\X\Z, \alpha}   h^\alpha}
    {f_{\Z,min}^2  \alpha !} + t_2 \bigg)
    \leq  8 t_2 C_{\psi } C_{\Lambda'} \bigg(1 + \frac{8 f_{\Z, max}^2}{f_{\Z, min}^2} \bigg)
    \leq t.
\end{align*}
As a consequence, we obtain that
\begin{align*}
     &\PP \Bigg( \max_{j=1, \dots, p'} \Big| \frac{1}{n'}
    \sum_{i=1}^{n'} Z'_{i,j} \xi_{i,n} \Big|
    >  t \Bigg)
    \leq 2 n' \exp \bigg( - \frac{n h^p f_{\Z, min}^2}{32 f_{\Z, max} \int K^2 + (8/3) C_K f_{\Z, min}} \bigg) \\
    & \hspace{5cm}
    + 2 n' \exp \bigg( - \frac{(n-1) h^{2p} t^2}
    {C_2+ C_3 t} \bigg),
\end{align*}
and we can apply Lemma \ref{lemma:lasso_bound} to get the claimed result.$\;\; \Box$







\section{Proofs of asymptotic results for $\hat\beta_{n,n'}$}

\subsection{Proof of Lemma \ref{lemma:beta_process_GG}}
\label{proof:lemma:beta_process_GG}

Using the definition (\ref{def:estimator_hat_beta}) of $\hat \beta_{n, n'}$, we get
\begin{align*}
    &\hat \beta_{n,n'} := \arg \min_{\beta \in \Rb^{p'}}
     \frac{1}{n'} \sum_{i=1}^{n'} \left( \Lambda(\hat \tau_{1,2|\Z=\z'_i}) - \psibm (\z'_i)^T\beta \right)^2 + \lambda_{n, n'} |\beta|_1  \\
    &\, = \arg \min_{\beta \in \Rb^{p'}}
     \frac{1}{n'} \sum_{i=1}^{n'} \left( \xi_{i,n} + \psibm (\z'_i)^T\beta^* - \psibm (\z'_i)^T\beta \right)^2 + \lambda_{n, n'} |\beta|_1  \displaybreak[0] \\
    &\, = \arg \min_{\beta \in \Rb^{p'}}
     \frac{1}{n'} \sum_{i=1}^{n'} \xi_{i,n}^2
    + \frac{2}{n'} \sum_{i=1}^{n'} \xi_{i,n} \psibm (\z'_i)^T (\beta^* - \beta)
    + \frac{1}{n'} \sum_{i=1}^{n'} \big( \psibm (\z'_i)^T (\beta^* - \beta) \big)^2
    + \lambda_{n, n'} |\beta|_1  \displaybreak[0] \\
    &\, = \arg \min_{\beta \in \Rb^{p'}}
     \frac{2}{n'} \sum_{i=1}^{n'} \xi_{i,n} \psibm (\z'_i)^T (\beta^* - \beta)
    + \frac{1}{n'} \sum_{i=1}^{n'}
    \big( \psibm (\z'_i)^T (\beta^* - \beta) \big)^2
    + \lambda_{n, n'} |\beta|_1.  \;\; \Box
\end{align*}

\subsection{Proof of Theorem \ref{thm:WeakConvLasso}}
\label{proof:thm:WeakConvLasso}

Let us define $r_{n,n'}:= (n h_{n,n'}^p)^{1/2}$,  $\u := r_{n,n'} (\beta - \beta^*)$ and $\hat\u_{n,n'} := r_{n,n'} (\hat \beta_{n, n'} - \beta^*)$, so that $\hat \beta_{n, n'} = \beta^* + \hat\u_{n,n'} / r_{n,n'}$.
By Lemma \ref{lemma:beta_process_GG},
$\hat \beta_{n, n'} = \arg \min_{\beta \in \Rb^{p'}} \GG_{n,n'}(\beta)$.
We have therefore
\begin{align*}
    &\hat\u_{n,n'} = \arg \min_{\u \in \Rb^{p'}}
    \Big[ \frac{-2}{n'} \sum_{i=1}^{n'}
    \xi_{i,n} \psibm (\z'_i)^T \frac{\u}{r_{n,n'}}
    + \frac{1}{n'} \sum_{i=1}^{n'}
    \big( \psibm (\z'_i)^T \frac{\u}{r_{n,n'}} \big)^2
    + \lambda_{n, n'} \big|\beta^* + \frac{\u}{r_{n,n'}} \big|_1 \Big],
\end{align*}
or $\hat\u_{n,n'} = \arg \min_{\u \in \Rb^{p'}} \FF_{n,n'}(\u)$,
where, for every $\u \in \Rb^{p'},$
$$    \FF_{n,n'}(\u) := \frac{-2 r_{n,n'}}{n'} \sum_{i=1}^{n'}
    \xi_{i,n} \psibm (\z'_i)^T \u
    + \frac{1}{n'} \sum_{i=1}^{n'}
    \left( \psibm (\z'_i)^T \u \right)^2 + \lambda_{n, n'} r_{n,n'}^2\Big( \big| \beta^* + \frac{\u}{r_{n,n'} } \big|_1
    - \big| \beta^* \big|_1 \Big).
$$
Note that, by Corollary \ref{cor:asymptNorm_xi}, we have
\begin{align*}
    \frac{2 r_{n,n'}}{n'} \sum_{i=1}^{n'} \xi_{i,n} \psibm (\z'_i)^T \u
    &= \frac{2}{n'} \sum_{i=1}^{n'} \sum_{j=1}^{p'} r_{n,n'} \xi_{i,n} \psi_j (\z'_i) u_{j}
    \indistrto \frac{2}{n'} \sum_{i=1}^{n'} \sum_{j=1}^{p'}
    W_{i} \psi_j (\z'_i) u_{j}.
\end{align*}
We also have, for any (fixed) $\u$ and when $n$ is large enough,
\begin{align*}
    \big| \beta^* + \frac{\u}{r_{n,n'} } \big|_1 - \big|  \beta^* \big|_1
    = \sum_{i=1}^{p'} \left( \frac{|u_i|}{r_{n,n'}} \1_{\{\beta^*_i = 0\}}
    + \frac{u_i}{r_{n,n'}}  \sgn(\beta^*_i) \1_{\{\beta^*_i \neq 0\}} \right).
\end{align*}
Therefore
$\lambda_{n, n'} r_{n,n'}^2 \Big( \big| \beta^* + \u/r_{n,n'}  \big|_1
- \big| \beta^* \big|_1 \Big)
\rightarrow \ell \sum_{i=1}^{p'} \big( |u_i| \1_{\{\beta_i^*=0\}}
+ u_i\sgn(\beta_i^*) \1_{\{\beta_i^* \neq 0\}} \big).$

\mds

We have shown that
$\FF_{n,n'}(\u) \indistrto \FF_{\infty, n'}(\u).$
Those functions are convex, hence the conclusion follows from the convexity argument.$\;\; \Box$

\subsection{Proof of Proposition \ref{prop:not_oracle_property}}
\label{proof:prop:not_oracle_property}

The proof closely follows Proposition 1 in~\cite{zou2006}.
It starts by noting that
$ \PP\left(\Sc_n= \Sc     \right) \leq \PP\left(\hat\beta_j = 0, \; \forall j\not\in \Sc \right)  .$
Because of the weak limit of $\hat\beta$ (Theorem~\ref{thm:WeakConvLasso} and the notations therein), this implies
$$ \lim\sup_n \PP\left(\hat\beta_j = 0, \; \forall j\not\in \Sc \right)
\leq \PP\left(u_j^* = 0, \; \forall j\not \in \Sc      \right) .$$
If $\ell=0$, then $\u^*$ is asymptotically normal, and the latter probability is zero.
Otherwise, $\ell\neq 0$ and define the Gaussian random vector $ \vec W_\psibm :=2\sum_{i=1}^{n'}
        W_{i} \psibm (\z'_i) /n'.$
The KKT conditions applied to $\FF_{\infty, n'}$ provide
$$ \vec W_\psibm  +  \frac{2}{n'}\sum_{i=1}^{n'} \psibm (\z'_i)  \psibm (\z'_i)^T \u^* + \ell  \v^* =0,$$
for some vector $\v^*\in \Rb^p$ whose components $v^*_j$ are less than one in absolute value when $j\not\in \Sc$, and $v^*_j=\sgn(\beta_j^*)$ when $j\in\Sc$.
If $u_j^* = 0$ for all $j\not \in \Sc$, we deduce
\begin{equation}
    (\vec W_\psibm)_{\Sc}
    + \Big[\frac{2}{n'}\sum_{i=1}^{n'} \psibm (\z'_i) \psibm (\z'_i)^T \Big]_{\Sc,\Sc} \u^*_{\Sc} + \ell \sgn (\beta_{\Sc}^*)
    = 0,\; \text{and}
\end{equation}
\begin{equation}
    \bigg| (\vec W_\psibm)_{\Sc^c} +  \Big[ \frac{2}{n'} \sum_{i=1}^{n'} \psibm (\z'_i) \psibm (\z'_i)^T \Big]_{\Sc^c,\Sc}  \u^*_{\Sc} \bigg| \leq \ell,
\end{equation}
componentwise and with obvious notations. Combining the two latter equations provides
\begin{equation}
    \bigg| (\vec W_\psibm)_{\Sc^c}  -  \Big[\sum_{i=1}^{n'} \psibm (\z'_i)  \psibm (\z'_i)^T \Big]_{\Sc^c,\Sc}
    \Big[\sum_{i=1}^{n'} \psibm (\z'_i)  \psibm (\z'_i)^T \Big]_{\Sc,\Sc}^{-1} \Big( \vec W_\psibm)_{\Sc}
    + \ell \sgn(\beta_{\Sc}^*) \Big)
    \bigg|\leq \ell,
\end{equation}
componentwise. Since the latter event is of probability strictly lower than one, this is still the case for the event
$\left\{ u_j^* = 0, \; \forall j\not \in \Sc \right\}$.$\;\; \Box$

\subsection{Proof of Theorem \ref{thm:WeakConvLassoAdaptive} }
\label{proof:thm:WeakConvLassoAdaptive}

The beginning of the proof is similar to the proof of Theorem~\ref{thm:WeakConvLasso}. With obvious notations,
$\check\u_{n,n'} = \arg \min_{\u \in \Rb^{p'}} \check\FF_{n,n'}(\u)$,
where for every $\u \in \Rb^{p'},$
\begin{eqnarray*}
\lefteqn{ \check\FF_{n,n'}(\u) := \frac{-2 r_{n,n'}}{n'} \sum_{i=1}^{n'}
    \xi_{i,n} \psibm (\z'_i)^T \u
    + \frac{1}{n'} \sum_{i=1}^{n'}
    \left( \psibm (\z'_i)^T \u \right)^2 }\\
    &+& \mu_{n, n'} r_{n,n'}^2
    \sum_{i=1}^{p'} \frac{1}{|\tilde\beta_i|^\delta}
    \left( |\beta_i^* + \frac{u_i}{r_{n,n'}}| - |\beta_i^*|    \right).
\end{eqnarray*}
If $\beta_i^*\neq 0$, then
$$  \frac{\mu_{n, n'} r_{n,n'}^2}{|\tilde\beta_i|^\delta}
    \left( |\beta_i^* + \frac{u_i}{r_{n,n'}}| - |\beta_i^*|    \right) =
    \frac{\mu_{n, n'} r_{n,n'}}{|\tilde\beta_i|^\delta} u_i \sgn(\beta_i^*)= \frac{\ell}{|\beta_i^*|^\delta} u_i \sgn(\beta_i^*) + o_P(1).
    $$
If $\beta_i^*=0$, then
$$  \frac{\mu_{n, n'} r_{n,n'}^2}{|\tilde\beta_i|^\delta}
    \left( |\beta_i^* + \frac{u_i}{r_{n,n'}}| - |\beta_i^*|    \right) =
  \frac{\mu_{n, n'} r_{n,n'}\nu_n^\delta}{|\nu_n\tilde\beta_i|^\delta}
    |u_i|.
$$
By assumption $\nu_n\tilde\beta_i=O_p(1)$, and the latter term tends to the infinity in probability iff $u_i\neq 0$.
As a consequence, if there exists some $i\not\in \Sc$ s.t. $u_i\neq 0$, then $\check\FF_{n,n'}(\u)$ tends to the infinity. Otherwise, $u_i=0$ when $i\not\in \Sc$ and
$\check\FF_{n,n'}(\u) \rightarrow \check\FF_{\infty,n'}(\u_{\Sc}) .$
Since $\check\FF_{\infty,n'}$ is convex, we deduce~\citep{kato2009asymptotics} that
$ \check\u_{\Sc} \rightarrow  \u^*_{\Sc}$, and $\check\u_{\Sc^c} \rightarrow  0_{\Sc^c},$
proving the asymptotic normality of $\check\beta_{n,n',\Sc}$.

\mds

Now, let us prove the oracle property. If $j\in \Sc$, then $\check\beta_j$ tends to $\beta_j$ in probability and $\PP(j\in \Sc_n)\rightarrow 1$. It suffices to show that $\PP(j\in \Sc_n) \rightarrow 0$ when $j\not\in \Sc$. If $j\not\in \Sc$ and
$j\in \Sc_n$, the KKT conditions on $\check\FF_{n,n'}$ provide
$$  \frac{-2 r_{n,n'}}{n'} \sum_{i=1}^{n'}
    \xi_{i,n} \psi_j (\z'_i)
    + \frac{2}{n'} \sum_{i=1}^{n'}
     \psi_j (\z'_i) \psibm (\z'_i)^T \check\u_{n,n'}  =-
    \frac{\mu_{n, n'} r_{n,n'} \nu_n^\delta}{|\nu_n\tilde\beta_j|^\delta} sign(\check u_j)\cdot$$
Due to the asymptotic normality of $\check\beta$ (that implies the one of $\check\u_{n,n'}$), the left hand side of the previous equation is asymptotically normal, when $\ell=0$. On the other side, the r.h.s. tends to the infinity in probability because $\nu_n\tilde\beta_j=O_P(1)$. Therefore, the probability of the latter event tends to zero when $n\rightarrow \infty$.
$\;\; \Box$

\subsection{Proof of Theorem \ref{thm:consistency_hatBeta_n_nprime}}
\label{proof:thm:consistency_hatBeta_n_nprime}

By Lemma \ref{lemma:beta_process_GG}, we have
$\hat \beta_{n, n'}
= \arg \min_{\beta \in \Rb^{p'}} \GG_{n,n'}(\beta)$, where
\begin{equation*}
    \GG_{n,n'}(\beta)
    := \frac{2}{n'} \sum_{i=1}^{n'} \xi_{i,n} \psibm (\z'_i)^T (\beta^* - \beta)
    + \frac{1}{n'} \sum_{i=1}^{n'}
    \big( \psibm (\z'_i)^T (\beta^* - \beta) \big)^2
    + \lambda_{n, n'} |\beta|_1.
\end{equation*}
Define also
$\GG_{\infty,n'}(\beta)
    :=  \sum_{i=1}^{n'}
    \big( \psibm (\z'_i)^T (\beta^* - \beta) \big)^2/n'
    + \lambda_{0} |\beta|_1$.
We have
\begin{align*}
    \big| \GG_{n,n'}(\beta) - \GG_{\infty,n'}(\beta) \big|
    \leq \bigg| \frac{2}{n'} \sum_{i=1}^{n'} \xi_{i,n} \psibm (\z'_i)^T (\beta^* - \beta) \bigg| + | \lambda_{n, n'} - \lambda_{0}|
    \times |\beta|_1.
\end{align*}
By assumption, the second term on the r.h.s. converges to $0$.
We now show that the first term on the r.h.s. is negligible.
Indeed, for every $\epsilon>0$,
\begin{align*}
    \PP\Big( \big\| \frac{1}{n'} \sum_{i=1}^{n'} \xi_{i,n} \psibm (\z'_i)  \big\| > \epsilon \Big)
    &\leq \PP\Big( \frac{ \|C_{\Lambda'} \|}{n'} \sum_{i=1}^{n'}  | \hat\tau_{\z'_i} - \tau_{\z'_i} | \times \|\psibm (\z'_i)  \big\| > \epsilon \Big) \\
    & \leq \sum_{i=1}^{n'} \PP\left(
    | \hat\tau_{\z'_i} - \tau_{\z'_i} | > Cst\epsilon  \right),
\end{align*}
where $Cst$ is the constant $(\|C_{\Lambda'} \|\times \|C_\psi\|)^{-1}$.
Apply Lemma~\ref{lemma:exponential_bound_KendallsTau} with the $t=f_{\Z,min}/4$ and $t'/\epsilon$ is a sufficiently small constant.
When $n$ is sufficiently large, we get
$$\PP \Big( |\hat \tau_{1,2|\Z=\z} - \tau_{1,2|\Z=\z} |
    > Cst \epsilon \Big)
    \leq 4 \exp \bigg( - n h^{2p} Cst' \bigg), $$
for some constant $Cst'>0$.
Thus, $ \sum_{i=1}^{n'} \xi_{i,n} \psibm (\z'_i)/n' =o_{\PP}(1)$, and
$ \GG_{n,n'}(\beta) = \GG_{\infty,n'}(\beta) + o_{\PP}(1)$ for every $\beta$.

\mds

Since $\sum_{i=1}^{n'} \psi(\z'_i)\psi(\z'_i)^T/n'$ tends towards a matrix $M_{\psi,\z'}$, deduce that
$\GG_{\infty,n'}(\beta)$ tends to $ \GG_{\infty,\infty}(\beta)$ when $n'\to\infty$.
Therefore, for all $\beta \in \Rb^{p'}$, $\GG_{n,n'}(\beta)$ weakly tends to $\GG_{\infty,\infty}(\beta)$.
By the convexity argument, we deduce that $\arg\min_\beta\GG_{n,n'}(\beta)$ weakly converges to $\arg\min_\beta\GG_{\infty,\infty}(\beta)$. Since the latter minimizer is non random,
the same convergence is true in probability.$\;\; \Box$

\section{Proof of Theorem \ref{thm:weak_conv_doubleAsympt}}
\label{proof:thm:weak_conv_doubleAsympt}

We start as in the proof of Theorem \ref{thm:WeakConvLasso}.
Define $\tilde r_{n,n'} := (n n' h_{n,n'}^p)^{1/2}$,
$\u := \tilde r_{n,n'} (\beta - \beta^*)$ and
$\hat \u_{n,n'} := \tilde r_{n,n'} (\hat \beta_{n, n'} - \beta^*)$,
so that $\hat \beta_{n, n'} = \beta^* + \hat\u_{n,n'} / \tilde r_{n,n'}$.
We define for every $\u \in \Rb^{p'}$,
\begin{align}
    \FF_{n,n'}(\u)
    &:= \frac{-2 \tilde r_{n,n'}}{n'} \sum_{i=1}^{n'}
    \xi_{i,n} \psibm (\z'_i)^T \u
    + \frac{1}{n'} \sum_{i=1}^{n'} \left( \psibm (\z'_i)^T \u \right)^2
    + \lambda_{n, n'} \tilde r_{n,n'}^2
    \Big( \big| \beta^* + \frac{\u}{\tilde r_{n,n'} } \big|_1 - \big| \beta^* \big|_1 \Big),
    \label{eq:def_F_n_nprime}
\end{align}
and we obtain
$\hat\u_{n,n'} = \arg \min_{\u \in \Rb^{p'}} \FF_{n,n'}(\u)$.

\begin{lemma}
    Under the same assumptions as in Theorem \ref{thm:weak_conv_doubleAsympt}, $T_1 := (\tilde r_{n,n'} / n')
    \sum_{i=1}^{n'} \xi_{i,n} \psibm(\z'_i)$ tends in law towards a Gaussian random vector $\Nc(0, V_2)$.
    \label{lemma:limit_T_1}
\end{lemma}
This lemma is proved in Section \ref{proof:lemma:limit_T_1}.
It will help to control the first term of Equation (\ref{eq:def_F_n_nprime}), which is simply $-2 T_1^T \u$.

\mds

Concerning the second term of Equation (\ref{eq:def_F_n_nprime}),
using Assumption \ref{assumpt:asymptNorm_joint}(iii), we have for every $\u \in \Rb^{p'}$
\begin{align}
    \frac{1}{n'} \sum_{i=1}^{n'} \left( \psibm (\z'_i)^T \u \right)^2
    \to \int \left( \psibm (\z')^T \u \right)^2 f_{\z', \infty}\, d\z'.
    \label{eq:limit_2rdterm_F}
\end{align}
This has to be read as a convergence of a sequence of real numbers indexed by $\u$, because the design points $\z'_i$ are deterministic.
We also have, for any $\u \in \Rb^{p'}$ and when $n$ is large enough,
\begin{align*}
    \big| \beta^* + \frac{\u}{\tilde r_{n,n'} } \big|_1 - \big|  \beta^* \big|_1
    = \sum_{i=1}^{p'} \Big( \frac{|u_i|}{\tilde r_{n,n'}} \1_{\{\beta^*_i = 0\}}
    + \frac{u_i}{\tilde r_{n,n'}}  \sgn(\beta^*_i) \1_{\{\beta^*_i \neq 0\}} \Big).
\end{align*}

Therefore, by Assumption \ref{assumpt:asymptNorm_joint}(ii)(b), for every $\u \in \Rb^{p'}$,
\begin{align}
    \lambda_{n, n'} \tilde r_{n,n'}^2
    \Big( \big| \beta^* + \frac{\u}{\tilde r_{n,n'} } \big|_1 - \big| \beta^* \big|_1 \Big) \to 0,
    \label{eq:limit_3rdterm_F}
\end{align}
when $(n,n')$ tends to the infinity.
Combining Lemma \ref{lemma:limit_T_1} and Equations (\ref{eq:def_F_n_nprime}-\ref{eq:limit_3rdterm_F}), and defining the function $\FF_{\infty, \infty}$ by
\begin{align*}
    \FF_{\infty, \infty}(\u) := 2 \tilde \W^T \u + \int \left( \psibm (\z')^T \u \right)^2 f_{\z', \infty}(\z') d\z', \, \u \in \Rb^{p'},
\end{align*}
where $\tilde W \sim \Nc(0, V_2)$, we obtain that every finite-dimensional margin of $\FF_{n,n'}$ converges weakly to the corresponding margin of $\FF_{\infty, \infty}$.
Now, applying the convexity lemma, we get
\begin{align*}
    \hat \u_{n,n'} \indistrto \u_{\infty, \infty}, \text{ where }
    \u_{\infty, \infty} := \arg \min_{\u \in \Rb^{p'}} \FF_{\infty, \infty}(\u).
\end{align*}
Since $\FF_{\infty, \infty}(\u)$ is a continuously differentiable convex function, we apply the first-order condition $\nabla \FF_{\infty, \infty}(\u) = 0$, which yields $2 \tilde \W  + 2 \int \psibm (\z') \psibm (\z')^T \u_{\infty, \infty} f_{\z', \infty}(\z') d\z' = 0$. As a consequence $\u_{\infty, \infty} = - V_1^{-1} \tilde \W \sim \Nc(0, \tilde V_{as}),$ using Assumption \ref{assumpt:asymptNorm_joint}(iv).
We finally obtain $\tilde r_{n,n'}  \big( \hat \beta_{n, n'} - \beta^* \big) \indistrto \Nc \big(0, \tilde V_{as} \big)$,
as claimed.$\;\; \Box$

\subsection{Proof of Lemma \ref{lemma:limit_T_1} : convergence of $T_1$}
\label{proof:lemma:limit_T_1}

Using a Taylor expansion, we have
\begin{align*}
    T_1 := \frac{\tilde r_{n,n'}}{n'} \sum_{i=1}^{n'} \xi_{i,n} \psibm (\z'_i)
    = \frac{\tilde r_{n,n'}}{n'} \sum_{i=1}^{n'}
    \Big( \Lambda \big(\hat \tau_{1,2|\Z = \z'_i} \big) - \Lambda \big(\tau_{1,2|\Z = \z'_i} \big) \Big) \psibm (\z'_i)
    = T_2 + T_3,
\end{align*}
where the main term is
\begin{align*}
    T_2 := \frac{\tilde r_{n,n'}}{n'} \sum_{i=1}^{n'} \Lambda' \big(\tau_{1,2|\Z = \z'_i} \big)
    \big(\hat \tau_{1,2|\Z = \z'_i} - \tau_{1,2|\Z = \z'_i} \big) \psibm(\z'_i),
\end{align*}
and the remainder is
\begin{align*}
    T_3 := \frac{\tilde r_{n,n'}}{n'} \sum_{i=1}^{n'} \alpha_{3,i}
    \big(\hat \tau_{1,2|\Z = \z'_i} - \tau_{1,2|\Z = \z'_i} \big)^2 \psibm(\z'_i),
\end{align*}
with $\forall i=1, \dots, n'$, $|\alpha_{3,i}| \leq C_{\Lambda''}/2$, by Assumption \ref{assumpt:asymptNorm_joint}\ref{assumpt:Lambda_C2}.

\medskip

Using the definition (\ref{def:hat_conditional_tau}) of $\hat \tau_{1,2|\Z = \z}$, the definition of the weights $w_{i,n}(\z)$ and the notation
$\bar \psibm(\z) := \Lambda' \big(\tau_{1,2|\Z = \z} \big) \psibm(\z),$ we rewrite $T_2=: T_4 + T_5$, where
\begin{align}
    &T_4 := \frac{\tilde r_{n,n'}}{n' n^2} \sum_{i=1}^{n'}
    \sum_{j_1 = 1}^n \sum_{j_2 = 1}^n
    \frac{K_h(\z'_i - \Z_{j_1}) K_h(\z'_i - \Z_{j_2})} {f_\Z^2(\z'_i)} \nonumber \\
    & \hspace{2cm} \times \, \Big(g^*(\X_{j_1}, \X_{j_2})
    - \EE \big[g^*(\X_{1}, \X_{2}) | \Z_{1} = \Z_{2} = \z'_i \big] \Big) \bar \psibm(\z'_i),
    \label{def:term_A_4} \displaybreak[0] \\
    &T_5 := \frac{\tilde r_{n,n'}}{n' n^2} \sum_{i=1}^{n'} \sum_{j_1 = 1}^n \sum_{j_2 = 1}^n K_h(\z'_i - \Z_{j_1}) K_h(\z'_i - \Z_{j_2})
    \bigg( \frac{1}{\hat f_\Z(\z'_i)^2}
    - \frac{1}{f_\Z(\z'_i)^2} \bigg) \nonumber \\
    & \hspace{2cm} \times  \,
    \Big( g^*(\X_{j_1}, \X_{j_2}) - \EE \big[g^*(\X_{1}, \X_{2}) | \Z_{1} = \Z_{2} = \z'_i \big] \Big) \bar \psibm(\z'_i).
    \label{def:term_A_5}
\end{align}
Note that we can put together the terms $(j_1, j_2)$ and $(j_2, j_1)$. This corresponds to the substitution of $g^*$ by its symmetrized version $\tilde g$. In the following, we will therefore assume that $g^*$ has been symmetrized without loss of generality.
The random variable $T_4$ can be seen (see Equation~(\ref{def:term_A_4})) as a sum of (indexed by $i$) U-statistics of order 2.
Its Hájek projection will yield the asymptotically normal dominant term of $T_2$.

\mds

To lighten notations, we denote $\tau_i:= \tau_{1,2|\Z_1=\Z_2=\z'_i}$, $f(\cdot,\cdot)=f_{\X,\Z}(\cdot,\cdot)$ and
$$  g_{i,j_1,j_2}:= g^*(\X_{j_1}, \X_{j_2})
    - \EE \big[g^*(\X_{1}, \X_{2}) | \z'_i \big]= g^*(\X_{j_1}, \X_{j_2}) - \tau_i.$$
Implicitly, all the expectations we will consider are expectations conditionally on the sequence of $\z'_i$, $i\geq 1$.

\mds

First note that, by usual $\alpha$-order limited expansions, we have
\begin{align*}
    \EE[T_4]
    & = \frac{\tilde r_{n,n'}}{n' n^2} \sum_{i=1}^{n'} n(n-1)
    \int \frac{K_h(\z'_i - \z_{1}) K_h(\z'_i - \z_{2})} {f_\Z^2(\z'_i)}  \big(g^*(\x_{1}, \x_{2})
    - \tau_i \big) \\
    & \times \bar \psibm(\z'_i)f(\x_1,\z_1)f(\x_2,\z_2)\, d\x_1 \, d\x_2\, d\z_1\, d\z_2  \\
    & - \frac{\tilde r_{n,n'}}{n' n} \sum_{i=1}^{n'} \tau_{i}
    \bar \psibm(\z'_i) \int \frac{K^2_h(\z'_i - \z)}{f_\Z^2(\z'_i)} 
    f(\x,\z) d\x\, d\z
    \displaybreak[0] \\
    & = \frac{(n-1)\tilde r_{n,n'}}{n' n} \sum_{i=1}^{n'} \int \frac{K(\t_1) K(\t_2)}{f_\Z^2(\z'_i)}  \big(g^*(\x_{1}, \x_{2}) - \tau_i \big) \\
    & \times \bar \psibm(\z'_i)f(\x_1,\z'_i - h\t_{1})f(\x_2,\z'_i - h\t_{1})\, d\x_1 \, d\x_2\, d\t_1\, d\t_2  \\
    & - \frac{\tilde r_{n,n'}}{n' n h^p} \sum_{i=1}^{n'}  \tau_{i} \bar \psibm(\z'_i)
    \int \frac{K^2(\t)}{f_\Z^2(\z'_i)}  f_{\X,\Z}(\x,\z'_i - h\t) d\x\, d\t \displaybreak[0] \\
    & = \frac{(n-1)\tilde r_{n,n'} h^{2\alpha }}{n' n} \sum_{i=1}^{n'} \int \frac{K(\t_1) K(\t_2)}{f_\Z^2(\z'_i)}  \big(g^*(\x_{1}, \x_{2}) - \tau_i \big) \\
    & \times \bar \psibm(\z'_i) d_\Z^{(\alpha)}f(\x_1,\z^*_i)\cdot \t_1^{(\alpha)} d_\Z^{(\alpha)}f(\x_2,\z^*_i)\cdot \t_2^{(\alpha)}
    \, d\x_1 \, d\x_2\, d\t_1\, d\t_2  \\
    & - \frac{\tilde r_{n,n'}}{n' n h^p} \sum_{i=1}^{n'}  \tau_{i}
    \int K^2\int \frac{\bar \psibm(\z'_i)}{f_\Z^2(\z'_i)}  f(\x,\z^*_i)\, d\x \\
    & = O\Big( \tilde r_{n,n'} h^{2\alpha} + \tilde r_{n,n'} /(nh^p)\Big)
    = O\Big( \sqrt{n n' h^{p + 4 \alpha} } + \sqrt{n' / (n h^p) } \Big)
    = o(1),
\end{align*}
under Assumption~\ref{assumpt:asymptNorm_joint}~\ref{assumpt:rates_n_nprime}.
Above, we have denoted by $\z^*_i$ some vectors in $\Rb^p$ s.t. $\| \z'_i - \Z_i^*\|_\infty < 1$. 
They depend on $\z'_i$, $\x_1$, $\x_2$ or $\x$, respectively.

\mds

Moreover, set
\begin{equation}
T_4 - \EE[T_4]=\frac{\tilde r_{n,n'}}{n' n^2} \sum_{i=1}^{n'} \sum_{j_1,j_2=1}^n \zeta_{i,j_1,j_2},
\label{T4centered}
\end{equation}
$$ \zeta_{i,j_1,j_2} =   \Big(
K_h(\z'_i - \Z_{j_1}) K_h(\z'_i - \Z_{j_2}) g_{i,j_1,j_2}
-  \EE\big[ K_h(\z'_i - \Z_{j_1}) K_h(\z'_i - \Z_{j_2}) g_{i,j_1,j_2} \big] \Big)\frac{\bar \psibm(\z'_i)} {f_\Z^2(\z'_i)} \cdot$$
Note that $Var(T_4)= \EE[T_4 T_4^T]+ o(1)$ and
$$ \EE[T_4 T_4^T]=\frac{\tilde r^2_{n,n'}}{(n')^2 n^4} \sum_{i_1,i_2=1}^{n'} \sum_{j_1,j_2=1}^n \sum_{j_3,j_4=1}^n \EE[ \zeta_{i_1,j_1,j_2} \zeta^T_{i_2,j_3,j_4} ].$$
By independence, $\EE[ \zeta_{i,j_1,j_2} \zeta^T_{i,j_3,j_4} ]=0$ when $\{j_1,j_2\}\cap \{j_3,j_4\}=\emptyset$.

\mds

Otherwise, assume that $j_1=j_3=j$ and there are no other identities among the four indices $(j_1,j_2,j_3,j_4)$.
Set
\begin{equation}
    \bar \zeta_{i} :=
    \EE\big[ K_h(\z'_i - \Z_{1}) K_h(\z'_i - \Z_{2}) g_{i,1,2} \big] \frac{\bar \psibm(\z'_i)} {f_\Z^2(\z'_i)}.
\label{zeta_i}
\end{equation}
Then,
$$ \EE[ \zeta_{i_1,j,j_2} \zeta^T_{i_2,j,j_4} ] = \zeta_{i_1,j,j_2,i_2,j,j_4} - \bar \zeta_{i_1} \bar \zeta_{i_2}^T, $$
where
\begin{eqnarray*}
    \lefteqn{\zeta_{i_1,j,j_2,i_2,j,j_4}:= \EE\Big[
    K_h(\z'_{i_1} - \Z_{j}) K_h(\z'_{i_1} - \Z_{j_2})
    K_h(\z'_{i_2} - \Z_{j}) K_h(\z'_{i_2} - \Z_{j_4})
    g_{i_1,j,j_2} g_{i_2,j,j_4}^T \Big] }\\
    &\times &  \frac{\bar \psibm(\z'_{i_1}) \bar \psibm(\z'_{i_2})^T}{f_\Z^2(\z'_{i_1})f_\Z^2(\z'_{i_2})}  \\
    & = &
    \frac{  \bar \psibm(\z'_{i_1}) \bar \psibm(\z'_{i_2})^T}{h^p f_\Z^2(\z'_{i_1})f_\Z^2(\z'_{i_2})}
    \int K(\t_1)K(\t_2) K\big(\frac{\z'_{i_2}-\z'_{i_1}}{h} + \t_1 \big) K(\t_4)
    \big( g^*(\x_1,\x_2)-\tau_{i_1} \big) \\
    &\times & \big( g^*(\x_1,\x_4)-\tau_{i_2} \big)
    f(\x_1,\z'_{i_1}- h\t_1) f(\x_2,\z'_{i_1}- h\t_2) f(\x_4,\z'_{i_4}- h\t_4)\, d\x_1\,d\x_2\,d\x_4\,d\t_1\,d\t_2\,d\t_4.
\end{eqnarray*}
By assumption, $\zeta_{i_1,j,j_2,i_2,j,j_4}$ is zero when $i_1\neq i_2$. Otherwise, when $i_1=i_2=i$,
\begin{eqnarray*}
\lefteqn{\zeta_{i,j,j_2,i,j,j_4} \simeq  \frac{  \bar \psibm(\z'_{i}) \bar \psibm(\z'_{i})^T}{h^p f_\Z(\z'_{i})} \int K^2
\int \big( g^*(\x_1,\x_2)-\tau_{i} \big)
 \big( g^*(\x_1,\x_4)-\tau_{i} \big) }\\
 &\times & f_{\X|\Z}(\x_1 | \z'_i) f_{\X|\Z}(\x_2|\z'_i) f_{\X|\Z}(\x_4| \z'_i)\, d\x_1\,d\x_2\,d\x_4 := C_{i,1,2,4}/h^p.
\end{eqnarray*}

It is easy to check that the terms with other identities among the four indices $j_k$,
as $\zeta_{i,j,j_2,i,j,j_2}$ or $\zeta_{i,j,j_2,i,j,j}$ will induce negligible remainder terms. Therefore, we get
$$ \frac{\tilde r^2_{n,n'}}{(n')^2 n^4} \sum_{i_1,i_2=1}^{n'} \sum_{j,j_2,j_4=1}^n \zeta_{i_1,j,j_2,i_2,j,j_4}
\simeq \frac{1}{n'} \sum_{i=1}^{n'} C_{i,1,2,4}.$$

Concerning the terms induced by the product of two $\bar \zeta_{i}$, note that, by limited expansions,
\begin{eqnarray*}
\lefteqn{\bar \zeta_{i} = \frac{\bar \psibm(\z'_i)} {f_\Z^2(\z'_i)}\int K_h(\z'_i - \z_1) K_h(\z'_i - \z_{2}) \big(g^*(\x_1,\x_2) - \tau_i \big) f(\x_1,\z_1) f(\x_2,\z_2)
\, d\x_1\,d\z_1 \, d\x_2\,d\z_2       }\\
&=& \frac{\bar \psibm(\z'_i)} {f_\Z^2(\z'_i)}\int K(\t_1) K(\t_{2}) \big(g^*(\x_1,\x_2) - \tau_i \big) f(\x_1,\z'_i- h\t_1) f(\x_2,\z'_i- h\t_2)
\, d\x_1\,d\t_1 \, d\x_2\,d\t_2 \\
&=& \frac{h^{2\alpha}\bar \psibm(\z'_i)} {f_\Z^2(\z'_i)}\int K(\t_1) K(\t_{2}) \big(g^*(\x_1,\x_2) - \tau_i \big) d^{(\alpha)}_\Z f(\x_1,\z^*_i)\cdot \t_1^{(\alpha)}
d^{(\alpha)}_\Z f(\x_2,\z^*_i)\cdot \t_2^{(\alpha)}
\, d\x_1\,d\t_1 \, d\x_2\,d\t_2 ,
\end{eqnarray*}
with the same notations as above. As a consequence, $\sup_i \bar \zeta_{i}= O(h^{2\alpha})$ and
$$\frac{\tilde r^2_{n,n'}}{(n')^2 n^4} \sum_{i_1,i_2=1}^{n'} \sum_{j,j_2,j_4=1}^n \bar \zeta_{i_1} \bar \zeta_{i_2}
\simeq   \frac{\tilde r^2_{n,n'}}{n} \Big( \frac{1}{n'}\sum_{i=1}^{n'} \bar \zeta_{i,1,2} \Big)^2
=O\Big(\frac{h^{4\alpha}\tilde r^2_{n,n'}}{n}\Big)=O(n' h^{4\alpha+p})=o(1).$$

Therefore, we obtain
$$ \frac{\tilde r^2_{n,n'}}{(n')^2 n^4} \sum_{i_1,i_2=1}^{n'} \sum_{j,j_2,j_4=1}^n \EE[ \zeta_{i_1,j,j_2} \zeta^T_{i_2,j,j_4} ]
\simeq \frac{1}{n'} \sum_{i=1}^{n'} C_{i,1,2,4}.$$
To calculate $\EE[T_4 T_4^T]$, there are three other similar terms, that respectively correspond to the cases $j_1=j_4$, $j_2=j_3$ or $j_2=j_4$.
Therefore, we deduce
\begin{eqnarray*}
\lefteqn{
Var(T_4) \simeq \EE[T_4 T_4^T] \simeq \frac{4}{n'} \sum_{i=1}^{n'} C_{i,1,2,4}   }\\
&\simeq &  4 \int K^2
\int \frac{  \bar \psibm(\z) \bar \psibm(\z)^T}{ f_\Z(\z)}
\int \big( g^*(\x_1,\x_2)-\tau_{1,2|\Z_1=\Z_2=\z} \big)
 \big( g^*(\x_1,\x_4)-\tau_{1,2|\Z_1=\Z_2=\z} \big) \\
 &\times & f_{\X|\Z}(\x_1 | \z) f_{\X|\Z}(\x_2|\z) f_{\X|\Z}(\x_4| \z)
 f_{\z',\infty}(\z)  \, d\x_1\,d\x_2\,d\x_4\, d\z,
\end{eqnarray*}
that is equal to the so-called variance-covariance matrix $V_2$.
Now assume that $T_4-\EE[T_4]$ is asymptotically normal, i.e. 
$T_{4}-\EE[T_4] \indistrto \Nc(0, V_2)$. This result will be proved in Subsection~\ref{proof:lemma:limit_T_4}.

\mds

Let us decompose the term $T_5$, as defined in Equation (\ref{def:term_A_5}). For every $i=1, \dots, n'$, a usual Taylor expansion yields
\begin{align*}
    \frac{1}{\hat f_{\Z}^2(\z'_i)} - \frac{1}{f_{\Z}^2(\z'_i)}
    = \frac{1}{f_{\Z}^2(\z'_i)} \Big\{
    \dfrac{1}{\Big( 1 + \dfrac{\hat f_{\Z}(\z'_i) - f_{\Z}(\z'_i)}{f_{\Z}(\z'_i)} \Big)^2 } - 1 \Big\}
    = -2 \frac{\hat f_{\Z}(\z'_i) - f_{\Z}(\z'_i)}{f_{\Z}^3(\z'_i)} + T_{7,i},
\end{align*}
where
\begin{align*}
    T_{7,i} = \frac{3}{f_{\Z}^2(\z'_i)} (1 + \alpha_{7,i})^{-4}
    \Big( \frac{\hat f_{\Z}(\z'_i) - f_{\Z}(\z'_i)}{f_{\Z}(\z'_i)} \Big)^2,
    \text{ for some } |\alpha_{7,i}| \leq \Big| \frac{\hat f_{\Z}(\z'_i) - f_{\Z}(\z'_i)}{f_{\Z}(\z'_i)} \Big|.
\end{align*}
Therefore, we obtain the decomposition $T_5 = -2 \, T_{6} + T_{7}$, where
\begin{align*}
    &T_{6} := \frac{\tilde r_{n,n'}}{n' n^2} \sum_{i=1}^{n'} \sum_{j_1 = 1}^n \sum_{j_2 = 1}^n K_h(\z'_i - \Z_{j_1}) K_h(\z'_i - \Z_{j_2}) \Big(
    \frac{\hat f_{\Z}(\z'_i) - f_{\Z}(\z'_i)}{f_{\Z}^3(\z'_i)} \Big) \nonumber \\
    & \hspace{4cm} \times \,
    \Big( g^*(\X_{j_1}, \X_{j_2}) - \EE \big[g^*(\X_{1}, \X_{2}) | \Z_{1} = \Z_{2} = \z'_i \big] \Big) \bar \psibm(\z'_i), \displaybreak[0] \\
    &T_{7} := \frac{\tilde r_{n,n'}}{n' n^2} \sum_{i=1}^{n'} \sum_{j_1 = 1}^n \sum_{j_2 = 1}^n K_h(\z'_i - \Z_{j_1}) K_h(\z'_i - \Z_{j_2})     T_{7,i} \nonumber \\
    & \hspace{4cm} \times \,
    \Big( g^*(\X_{j_1}, \X_{j_2}) - \EE \big[g^*(\X_{1}, \X_{2}) | \Z_{1} = \Z_{2} = \z'_i \big] \Big) \bar \psibm(\z'_i).
\end{align*}

Summing up all the previous equations, we get
\begin{align}
    T_1 = \big(T_{4} - \EE[T_4] \big)  - 2 \, T_{6} + T_{7} + T_3 + o(1).
    \label{eq:decomposition_T_1}
\end{align}
Afterwards, we will prove that all the remainders terms $T_6$, $T_7$ and $T_3$ are negligible, i.e. they tend to zero in probability.
These results are respectively proved in Subsections \ref{proof:lemma:limit_T_6}, \ref{proof:lemma:limit_T_7} and \ref{proof:lemma:limit_T_3}.
Combining all these elements with the asymptotic normality of $T_4$ (proved in Subsection~\ref{proof:lemma:limit_T_4}), we get $T_1 \indistrto \Nc(0, V_2)$, as claimed.$\;\; \Box$

\subsection{Proof of the asymptotic normality of $T_4$}
\label{proof:lemma:limit_T_4}
We will lead the usual H\'ajek projection of $T_4$.
To weaken notations, denote $\EE[ \zeta_{i,j_1,j_2} | \X_{j_1},\Z_{j_1}]:=\EE[ \zeta_{i,j_1,j_2} | j_1]$.
Then, recalling~(\ref{T4centered}), we can write
$$T_4 - \EE[T_4]= T_{4,1} +  T_{4,2} +  T_{4,3} ,\;\; \text{with}$$
$$  T_{4,1} := \frac{2\tilde r_{n,n'}}{n' n^2} \sum_{i=1}^{n'} \sum_{j_1, j_2=1}^n \1(j_1\neq j_2) \EE[ \zeta_{i,j_1,j_2} | j_1] ,$$
$$  T_{4,2} := \frac{2\tilde r_{n,n'}}{n' n^2} \sum_{i=1}^{n'} \sum_{j=1}^n \EE[ \zeta_{i,j,j} | j] ,\;\;\text{and}$$
$$  T_{4,3} := \frac{\tilde r_{n,n'}}{n' n^2} \sum_{i=1}^{n'} \sum_{j_1, j_2=1}^n \Big( \zeta_{i,j_1,j_2} - \EE[ \zeta_{i,j_1,j_2}| j_1]- \EE[ \zeta_{i,j_1,j_2}| j_2]\Big) .$$
We will prove that $T_{4,2}$ and $T_{4,3}$ are $o_P(1)$. Therefore, the asymptotic normality of $T_4$ reduces to the one of $T_{4,1}$.

\mds

Note that $n T_{4,1} /2(n-1) =  \sum_{j=1}^n \beta_{j,n,n'}$, where
 $$\beta_{j,n,n'} :=  \frac{\tilde r_{n,n'}}{n' n} \sum_{i=1}^{n'} \EE[ \zeta_{i,j,0} | j],\; j=1,\ldots,n,$$
by formally considering a random vector $\Z_0$ that is independent of the other $\Z_j$, $j\geq 1$. Therefore, we get a triangular array of random vectors
$(\beta_{j,n,n'})_{j=1,\ldots,n}$, s.t., for a fixed $n$, the variables $\beta_{j,n,n'}$ are mutually independent given the vectors $\z'_i$, $i\geq 1$.
Let us check Lyapunov's sufficient condition, that will imply the asymptotic normality of $T_{4,1}$.
In other words, it is sufficient to prove that
\begin{equation}
\sum_{j=1}^n \| \beta_{j,n,n'} \|^3_{\infty} \longrightarrow 0,
\end{equation}
when $n$ and $n'$ tend to the infinity. Recalling~(\ref{zeta_i}), we can rewrite
$$\beta_{j,n,n'} = \frac{\tilde r_{n,n'}}{n' n} \sum_{i=1}^{n'} \Big\{
 K_h(\z'_i - \Z_j) \frac{\bar \psibm(\z'_i)} {f_\Z^2(\z'_i)}  \int K_h(\z'_i - \z) \big( g^*(\x,\X_j) - \tau_i \big) f(\x,\z) \, d\x\, d\z - \bar \zeta_{i}
         \Big\} := \frac{\tilde r_{n,n'}}{n' n} \sum_{i=1}^{n'} \gamma_{i,j},$$
where $\sup_i\bar \zeta_{i} := O(h^{2\alpha})$.
Note that
$$ \| \beta_{j,n,n'}\|^3_\infty  \leq p^3\frac{\tilde r^3_{n,n'}}{(n')^3 n^3} \sum_{i_1,i_2,i_3=1}^{n'}
\| \gamma_{i_1,j}\|_\infty \| \gamma_{i_2,j}\|_\infty\| \gamma_{i_3,j}\|_\infty. $$
The terms that that involve some products by the means $\bar \zeta_{i_k}$, $k=1,2,3$, are negligible and they may be forgotten here.
For some constants $Cst$, this provides
\begin{eqnarray*}
\lefteqn{
\sum_{j=1}^n \EE\Big[ \| \beta_{j,n,n'} \|^3_{\infty} \Big] \leq
\frac{Cst \, \tilde r^3_{n,n'}}{(n')^3 n^3} \sum_{j=1}^n \sum_{i_1,i_2,i_3=1}^{n'} \frac{\| \bar \psibm \|_\infty (\z'_{i_1}) \| \bar \psibm \|_\infty (\z'_{i_2}) \| \bar \psibm \|_\infty (\z'_{i_3}) } {f_\Z^2(\z'_{i_1}) f_\Z^2(\z'_{i_2}) f_\Z^2(\z'_{i_3})}    }\\
&\times &  \EE\Big[
 \big| K_h(\z'_{i_1} - \Z_j)  \int K_h(\z'_{i_1} - \z_1) \big( g^*(\x_1,\X_j) - \tau_{i_1} \big) f(\x_1,\z_1) \, d\x_1\, d\z_1 \big| \\
&\times & \big| K_h(\z'_{i_2} - \Z_j)  \int K_h(\z'_{i_2} - \z_2) \big( g^*(\x_2,\X_j) - \tau_{i_2} \big) f(\x_2,\z_2) \, d\x_2\, d\z_2 \big| \\
&\times & \big| K_h(\z'_{i_3} - \Z_j)  \int K_h(\z'_{i_3} - \z_3) \big( g^*(\x_3,\X_j) - \tau_{i_3} \big) f(\x_3,\z_3) \, d\x_3\, d\z_3 \big|
\Big].
\end{eqnarray*}
By some now usual changes of variables, the latter expectations are zero when one of the three indices $i_1,i_2$ and $i_3$ is different from the others.
Thus, the non-zero expectations are obtained when $i_1=i_2=i_3$. In the latter case, we get
\begin{eqnarray*}
\lefteqn{
\sum_{j=1}^n \| \beta_{j,n,n'} \|^3_{\infty}  \leq
\frac{Cst\, \tilde r^3_{n,n'}}{(n')^3 n^2}  \sum_{i=1}^{n'}
\frac{\| \bar \psibm \|^3_\infty (\z'_{i}) } {f_\Z^6(\z'_{i}) }    }\\
&\times &  \int
 | K |^3_h(\z'_{i} - \z)  \Big|\int K_h(\z'_{i} - \z_1) \big( g^*(\x_1,\x) - \tau_{i} \big) f(\x_1,\z_1) \, d\x_1\, d\z_1 \Big|^3
 f(\x,\z) \, d\x\, d\z  \\
&\leq  &
\frac{Cst\, \tilde r^3_{n,n'}}{(n')^3 n^2 h^{2p}}  \sum_{i=1}^{n'}
\frac{\| \bar \psibm \|^3_\infty (\z'_{i}) } {f_\Z^6(\z'_{i}) }    \int
  |K|^3(\t)  \big| \int K(\t_1) \big( g^*(\x_1,\x) - \tau_{i} \big) f(\x_1,\z'_i-h\t_1) \, d\x_1\, d\t_1 \big|^3  \\
 &\times & f(\x,\z'_i-h\t) \, d\x\, d\t
 = O\Big( \frac{ \tilde r^3_{n,n'}}{(n')^2 n^2 h^{2p}} \Big)= O\Big( \frac{ 1}{(nn' h^p)^{1/2}}\Big)=o(1).
\end{eqnarray*}

Concerning the remainder terms $T_{4,2}$ and $T_{4,3}$,
note that $\EE[T_{4,2}]=\EE[T_{4,3}]=0$. Moreover, since $\EE[\zeta_{i,j,j} | j]$ is centered,
$$ \EE[ T_{4,2}T^T_{4,2}] = \frac{4\tilde r^2_{n,n'}}{(n')^2 n^4} \sum_{i_1,i_2=1}^{n'} \sum_{j=1}^n
\EE\Big[ \EE[ \zeta_{i_1,j,j} | j] \EE[\zeta^T_{i_2,j,j} | j]\Big] .$$

When $i_1\neq i_2$, some usual changes of variables yield
\begin{eqnarray*}
\lefteqn{
\EE\Big[ \EE[ \zeta_{i_1,j,j} | j] \EE[\zeta^T_{i_2,j,j} | j]\Big]=
\frac{\bar \psibm(\z'_{i_1})\bar \psibm(\z'_{i_2})^T} {f_\Z^2(\z'_{i_1})f_\Z^2(\z'_{i_2})} \tau_{i_1}\tau_{i_2} }\\
&\times & \bigg(
\EE\big[ K^2_h(\z'_{i_1} - \Z_j) K^2_h(\z'_{i_2} - \Z_j) \big] 
- \EE\big[ K^2_h(\z'_{i_1} - \Z_j) \big] \EE\big[K^2_h(\z'_{i_2} - \Z_j)    \big]
\bigg) =  O(h^{-2p}),
\end{eqnarray*}
uniformly w.r.t. $i$. By a similar reasoning, we can prove that
$$ \sup_i \EE\Big[ \EE[ \zeta_{i,j,j} | j] \EE[\zeta^T_{i,j,j} | j]\Big]=O(h^{-3p}). $$
Therefore,
$$ \EE[ T_{4,2}T^T_{4,2}] = O\Big(\frac{\tilde r^2_{n,n'}}{(n')^2  n^4 } \big((n')^2 n h^{-2p} + n' n h^{-3p} \big)\Big)
= O\big( \frac{n'}{n^2 h^{p} } +  \frac{1}{n^2 h^{2p} }   \big)= o(1).$$

Concerning $T_{4,3}$, this remainder term is centered and
\begin{eqnarray}
  \lefteqn{ \EE[ T_{4,3} T^T_{4,3} ]=
   \frac{\tilde r^2_{n,n'}}{(n')^2 n^4} \sum_{i_1,i_2=1}^{n'} \sum_{j_1, j_2=1}^n \sum_{j_3, j_4=1}^n
   \EE\Big[ \{ \zeta_{i_1,j_1,j_2} - \EE[ \zeta_{i_1,j_1,j_2}| j_1]- \EE[ \zeta_{i_1,j_1,j_2}| j_2] \} \nonumber }\\
    &\times &     \{ \zeta_{i_2,j_3,j_4} - \EE[ \zeta_{i_2,j_3,j_4}| j_3] - \EE[ \zeta_{i_2,j_3,j_4}| j_4] \}^T \Big]. \hspace{6cm}
\label{T43}
\end{eqnarray}
The expectations on the latter r.h.s. are zero when $\{j_1,j_2\} \cap \{j_3,j_4\}=\emptyset$ due to independence and the fact that the terms $\zeta_{i,j,j'}$ are centered.
Otherwise, there is at least an identity among the indices $j_k$, $k=1,\ldots,4$.
For instance, assume $j_1=j_3=j$ and $j\neq j_2 \neq j_4$. Then,
\begin{eqnarray*}
\lefteqn{ \EE\Big[ \{ \zeta_{i_1,j,j_2} - \EE[ \zeta_{i_1,j,j_2}| j] - \EE[ \zeta_{i_1,j,j_2}| j_2]  \}
    \{ \zeta_{i_2,j,j_4} - \EE[ \zeta_{i_2,j,j_4}| j]- \EE[ \zeta_{i_2,j,j_4}| j_4] \}^T \Big] }\\
    &=&
    \EE\Big[ \{ \zeta_{i_1,j,j_2} - \EE[ \zeta_{i_1,j,j_2}| j]  \}
    \{ \zeta_{i_2,j,j_4} - \EE[ \zeta_{i_2,j,j_4}| j] \}^T \Big] \\
    &=& \EE\bigg[ \EE\Big[ \{ \zeta_{i_1,j,j_2} - \EE[ \zeta_{i_1,j,j_2}| j] \}
    \{ \zeta_{i_2,j,j_4} - \EE[ \zeta_{i_2,j,j_4}| j] \}^T \Big| j \Big] \bigg] \\
    &=&
    \EE\Big[ \EE\big[  \zeta_{i_1,j,j_2} \zeta^T_{i_2,j,j_4} | j \big] \Big]
    -  \EE\Big[ \EE\big[  \zeta_{i_1,j,j_2} | j \big]  \EE\big[\zeta^T_{i_2,j,j_4} | j \big] \Big]=0.
\end{eqnarray*}
Due to the symmetry of the latter cross-products, all cases of a single identity among the $j_k$, $k=1,\ldots,4$, yield the same result. 
Therefore, we need (at least) two identities among them to obtain non zero covariances in the calculation of $\EE[ T_{4,3} T^T_{4,3} ]$.
Thus, let us assume that $j_1=j_3$ and $j_2=j_4$.
Then, the corresponding terms in~(\ref{T43}) is
\begin{align*}
   &\frac{\tilde r^2_{n,n'}}{(n')^2 n^4}
   \sum_{i_1,i_2=1}^{n'} \sum_{j_1, j_2=1}^n
   \EE\Big[ \{ \zeta_{i_1,j_1,j_2} - \EE[ \zeta_{i_1,j_1,j_2}| j_1]- \EE[ \zeta_{i_1,j_1,j_2}| j_2] \}
    \{ \zeta_{i_2,j_1,j_2} - \EE[ \zeta_{i_2,j_1,j_2}| j_1]- \EE[ \zeta_{i_2,j_1,j_2}| j_2] \}^T \Big] \\
    & \hspace{0.3cm} = \,\frac{\tilde r^2_{n,n'}}{(n')^2 n^4} \sum_{i_1,i_2=1}^{n'} \sum_{j_1, j_2=1}^n
   \bigg(\EE\Big[  \zeta_{i_1,j_1,j_2} \zeta_{i_2,j_1,j_2}^T  \Big]
    - 2\EE\Big[ \EE[ \zeta_{i_1,j_1,j_2}| j_1] \EE[ \zeta_{i_2,j_1,j_2}| j_1]^T  \Big] \bigg) \\
    & \hspace{0.3cm} =: \, v_{4,3,1} - v_{4,3,2}.
\end{align*}
By now usual techniques, we get
\begin{eqnarray*}
\lefteqn{
v_{4,3,1} =
   \frac{\tilde r^2_{n,n'}}{(n')^2 n^4} \sum_{i_1,i_2=1}^{n'} \sum_{j_1, j_2=1}^n \EE\Big[ \zeta_{i_1,j_1,j_2} \zeta^T_{i_2,j_1,j_2}  \Big] 
    \simeq 
   \frac{\tilde r^2_{n,n'}}{(n')^2 n^4} \sum_{i=1}^{n'} \sum_{j_1, j_2=1}^n \EE\Big[ \zeta_{i,j_1,j_2} \zeta^T_{i,j_1,j_2}  \Big] }\\
    &\simeq &
   \frac{\tilde r^2_{n,n'}}{(n')^2 n^2} \sum_{i=1}^{n'} \frac{\bar \psibm(\z'_i)\bar \psibm^T(\z'_i) } {f_\Z^4(\z'_i)} \int K_h^2(\z'_i - \z_1) K_h^2(\z'_i - \z_2) \big( g^*(\x_1,\x_2)-\tau_i  \big)^2 f(\x_1,\z_1) 
      \\
  &\times & f(\x_2,\z_2) \, d\x_1 \, d\z_1 \, d\x_2 \, d\z_2=
   O\Big( \frac{\tilde r^2_{n,n'}}{n' n^2h^{2p}} \Big) = O\big(  \frac{1}{nh^p}  \big)=o(1).
\end{eqnarray*}
Moreover,
\begin{eqnarray*}
\lefteqn{
v_{4,3,2} =
   \frac{\tilde r^2_{n,n'}}{(n')^2 n^4} \sum_{i_1,i_2=1}^{n'} \sum_{j_1, j_2=1}^n \EE\Big[ \EE[ \zeta_{i_1,j_1,j_2} | j_1] \EE[ \zeta^T_{i_2,j_1,j_2} | j_1]  \Big] }\\
   &\simeq &
   \frac{\tilde r^2_{n,n'}}{(n')^2 n^4} \sum_{i=1}^{n'} \sum_{j_1, j_2=1}^n \EE\Big[ \EE[ \zeta_{i,j_1,j_2} | j_1] \EE[ \zeta^T_{i,j_1,j_2}| j_1]  \Big] \\
    &\simeq &
   \frac{\tilde r^2_{n,n'}}{(n')^2 n^2} \sum_{i=1}^{n'} \frac{\bar \psibm(\z'_i)\bar \psibm^T(\z'_i) } {f_\Z^4(\z'_i)}
   \int K_h^2(\z'_i - \z_1)   K_h(\z'_i - \z_2) K_h(\z'_i - \z_3) \big( g^*(\x_1,\x_2)-\tau_i  \big)  \\
   &\times & \big( g^*(\x_1,\x_3)-\tau_i  \big) f(\x_1,\z_1) f(\x_2,\z_2) f(\x_3,\z_3) \, d\x_1 \, d\z_1 \, d\x_2 \, d\z_2 \, d\x_3 \, d\z_3  \\
  &= & O\Big( \frac{\tilde r^2_{n,n'}}{n' n^2h^{p}} \Big) = O\big(  \frac{1}{n}  \big)=o(1).
\end{eqnarray*}

Another case of two identities occurs when  $j_1=j_4$ and $j_2=j_3$, but it can be dealt similarly.
Then, we have proved that $\EE[ T_{4,3} T^T_{4,3} ]=o(1)$ and $T_{4,3}=o_P(1)$.

\subsection{Convergence of $T_{6}$ to $0$}
\label{proof:lemma:limit_T_6}

Replacing $\hat f_\Z$ in the definition of $T_6$ above by the normalized sum of the kernels, we get
\begin{align*}
    T_{6}
    &=\frac{ \tilde r_{n,n'}}{n' n^2} \sum_{i=1}^{n'} \sum_{j_1 = 1}^n \sum_{j_2 = 1}^n \frac{K_h(\z'_i - \Z_{j_1}) K_h(\z'_i - \Z_{j_2})}{f_{\Z}^3(\z'_i)}
     \big( \EE[\hat f_{\Z}(\z'_i)] - f_{\Z}(\z'_i) \big) \\
     & \times \,
     \Big( g^*(\X_{j_1}, \X_{j_2}) - \EE \big[g^*(\X_{1}, \X_{2}) | \Z_{1} = \Z_{2} = \z'_i \big] \Big) \bar \psibm(\z'_i) \displaybreak[0] \\
    &+ \frac{\tilde r_{n,n'}}{n' n^3} \sum_{i=1}^{n'} \sum_{j_1 = 1}^n \sum_{j_2 = 1}^n \sum_{j_3 = 1}^n \frac{K_h(\z'_i - \Z_{j_1}) K_h(\z'_i - \Z_{j_2})}{f_{\Z}^3(\z'_i)}
    \big( K_h(\z'_i - \Z_{j_3}) - \EE[\hat f_{\Z}(\z'_i)] \big) \\
    &  \times  \,
    \Big( g^*(\X_{j_1}, \X_{j_2}) - \EE \big[g^*(\X_{1}, \X_{2}) | \Z_{1} = \Z_{2} = \z'_i \big] \Big) \bar \psibm(\z'_i) =: T_{6,1}+T_{6,2}.
\end{align*}
The first term $T_{6,1}$ is a bias term. By Assumptions~\ref{assumpt:kernel_integral}-\ref{assumpt:f_Z_Holder},
$$ \sup_{i=1,\ldots, n'} \big| \EE[\hat f_{\Z}(\z'_i)] - f_{\Z}(\z'_i)  \big| = O(h^\alpha). $$
The sum of the diagonal terms in $T_{6,1}$ is
$$
-\frac{ \tilde r_{n,n'}}{n' n^2} \sum_{i=1}^{n'} \sum_{j = 1}^n \frac{K^2_h(\z'_i - \Z_{j})}{f_{\Z}^3(\z'_i)}
     \big( \EE[\hat f_{\Z}(\z'_i)] - f_{\Z}(\z'_i) \big)
     \EE \big[g^*(\X_{1}, \X_{2}) | \Z_{1} = \Z_{2} = \z'_i \big]  \bar \psibm(\z'_i) ,$$
that is $ O_{\PP} \big(\tilde r_{n,n'}h^\alpha / (nh^p)\big)$. The sum of the extra-diagonal terms in $T_{6,1}$ is the r.v.
\begin{align*}
    \bar T_{6,1}
    &:=\frac{ \tilde r_{n,n'}}{n' n^2} \sum_{i=1}^{n'} \sum_{1\leq j_1 \neq j_2 \leq n} \frac{K_h(\z'_i - \Z_{j_1}) K_h(\z'_i - \Z_{j_2})}{f_{\Z}^3(\z'_i)}
     \big( \EE[\hat f_{\Z}(\z'_i)] - f_{\Z}(\z'_i) \big) \\
     & \hspace{2cm} \times \,
     \Big( g^*(\X_{j_1}, \X_{j_2}) - \EE \big[g^*(\X_{1}, \X_{2}) | \Z_{1} = \Z_{2} = \z'_i \big] \Big) \bar \psibm(\z'_i) \displaybreak[0].
\end{align*}
Note that $\z\mapsto f_{\Z}(\z)$ and
$(\z_1,\z_2) \mapsto \EE \big[g^*(\X_{1}, \X_{2}) | \Z_{1} =\z_1, \Z_{2} = \z_2 \big]$ are $\alpha$-times continuously differentiable on $\Zc$ and $\Zc^2$ respectively, because of Assumptions~\ref{assumpt:f_Z_Holder} and~\ref{assumpt:f_XZ_Holder}.
By $\alpha$-order Taylor expansions of such terms, they yield some factors $h^\alpha$.
It is easy to check that the expectation of $(\bar T_{6,1})^2$ is of order $\tilde r^2_{n,n'}h^{2\alpha} / (n^2 h^{2p})$.
Therefore,
$$T_{6,1}= O_{\PP} \Big(\frac{\tilde r_{n,n'}h^\alpha}{nh^p}\Big)
=O_{\PP} \Big(\frac{(n')^{1/2} h^{\alpha} }{\sqrt{nh^p}}\Big)=o_{\PP}(1). $$

Concerning $T_{6,2}$, we can assume that the indices $j_1$, $j_2$ and $j_3$ are pairwise distinct. Indeed, the cases of one or two identities among such indices can be easily dealt. They yield an upper bound that is $O_{\PP}(\tilde r_{n,n'}h^\alpha/(nh^p))$ as above, and they are negligible.
Once we remove such terms from the triple sums (indexed by $(j_1,j_2,j_3)$) defining $T_{6,2}$, we get the centered r.v. $\bar T_{6,2}$.
Let us calculate the second moment of $\bar T_{6,2}$.
\begin{align*}
    & \EE \big[\bar T^2_{6,2} \big]
    := \frac{n n' h^p}{n'{}^2 n^6}
    \sum_{i_1=1}^{n'} \sum_{i_2=1}^{n'} \sum_{1\leq j_1 \neq j_2 \neq j_3 \leq n} \sum_{1\leq j_4 \neq j_5 \neq j_6 \leq n}
    \EE \Bigg[ \frac{K_h(\z'_{i_1} - \Z_{j_1}) K_h(\z'_{i_1} - \Z_{j_2})} {f_{\Z}^3(\z'_{i_1})} \\
    &  \times  \big( K_h(\z'_{i_1} - \Z_{j_3}) - \EE[ \hat f_{\Z}(\z'_{i_1})] \big)\,
    \Big( g^*(\X_{j_1}, \X_{j_2}) - \EE \big[g^*(\X_{1}, \X_{2}) | \Z_{1} = \Z_{2} = \z'_{i_1} \big] \Big) \bar \psibm(\z'_{i_1}) \\
    &  \times  \, \frac{K_h(\z'_{i_2} - \Z_{j_4}) K_h(\z'_{i_2} - \Z_{j_5})}{f_{\Z}^3(\z'_{i_2})}
    \big( K_h(\z'_{i_2} - \Z_{j_6}) - \EE[\hat  f_{\Z}(\z'_{i_2})] \big) \\
    & \hspace{1cm} \times  \,
    \Big( g^*(\X_{j_4}, \X_{j_5}) - \EE \big[g^*(\X_{1}, \X_{2}) | \Z_{1} = \Z_{2} = \z'_{i_2} \big] \Big) \bar \psibm(\z'_{i_2})^T \Bigg] \\
    &=: \frac{n n' h^p}{n'{}^2 n^6}     \sum_{i_1,i_2=1}^{n'} \sum_{1\leq j_1 \neq j_2 \neq j_3 \leq n} \sum_{1\leq j_4 \neq j_5 \neq j_6 \leq n} E_{i_1,i_2,j_1-j_6}.
\end{align*}
When all the indices of the latter sums are different, the latter expectation is zero.
Non zero terms above are obtained only when $j_3$ and $j_6$ are equal to some other indices.
In the case $j_3=j_6$ and no other identity among the indices, we obtain two extra factors $h^{\alpha}$ through $\alpha$-order limited expansions
of $(\z_1,\z_2) \mapsto \EE \big[g^*(\X_{1}, \X_{2}) | \Z_{1} =\z_1, \Z_{2} = \z_2 \big]$. This yields an order
$O(n n' h^{p+2\alpha} /(n h^p))$.
When $j_3$ and $j_6$ are equal to two different indices ($j_3=j_4$ and $j_6=j_2$, e.g.), we lose another factor $h^p$ but we still benefit from
the two latter factors $h^{\alpha}$.
This yields an upper bound
$O(n n' h^{p+2\alpha} /(n^2 h^{2p}))=o(1)$. The other situations can be managed similarly.
We get
$$ \EE \big[\bar T^2_{6,2} \big] = O\Big( \frac{n n' h^{p+2\alpha}}{ n h^p }\Big)=o(1).$$

\mds

Globally, we obtain $T_6 \to 0$ in probability under Assumptions \ref{assumpt:asymptNorm_joint}(ii)(a).
 $\;\;\Box$

\subsection{Convergence of $T_{7}$ to $0$}
\label{proof:lemma:limit_T_7}
Since $\sup_{i=1,\ldots,n'}   | \hat f_{\Z}(\z'_i) - f_{\Z}(\z'_i) |=o_{\PP}(1)$, note that
$$\sup_{i=1,\ldots,n'}| T_{7,i} | \leq \frac{6}{f_{\Z,min}^4}
\sup_{i=1,\ldots,n'}   | \hat f_{\Z}(\z'_i) - f_{\Z}(\z'_i) |^2,$$ with a probability arbitrarily close to one.
Apply Lemma \ref{lemma:bound_f_hat_f} with a fixed $t > 0$ and $\z = \z'_i$ for each $i=1, \dots, n'$
\begin{align*}
    \PP \bigg( \sup_{i=1,\ldots,n'}| T_{7,i} |
    \geq \frac{6}{f_{\Z,min}^4} \left( \frac{  C_{K, \alpha} h^{\alpha}}{ \alpha !} + t \right)^2 \bigg)
    \leq 2 n' \exp \bigg( - \frac{n h^p t^2}{2 f_{\Z, max} \int K^2 + (2/3) C_K t} \bigg).
\end{align*}
Set $t \propto h^{\alpha/2}$.
Deduce $\sup_{i=1,\ldots,n'}| T_{7,i} | = O_{\PP}(h^{\alpha})$
since $nh^{p+\alpha}/\ln n' \to \infty$ by assumption. Then,
\begin{align*}
    | T_{7} |
    & \leq \frac{\tilde r_{n,n'}}{n' n^2} \sup_i |T_{7,i} | \sum_{i=1}^{n'} |\bar \psibm(\z'_i)| \\
    &  \times \sum_{j_1 = 1}^n \sum_{j_2 = 1}^n
    |K|_h(\z'_i - \Z_{j_1}) |K|_h(\z'_i - \Z_{j_2}) \Big| g^*(\X_{j_1}, \X_{j_2}) - \EE \big[g^*(\X_{1}, \X_{2}) | \Z_{1} = \Z_{2} = \z'_i \big] \Big|.
\end{align*}
The expectation of the double sum is $O(h^\alpha)$, by an $\alpha$-order limited expansion
of $(\z_1,\z_2) \mapsto \EE \big[g^*(\X_{1}, \X_{2}) | \Z_{1} =\z_1, \Z_{2} = \z_2 \big]$.
Then, by Markov's inequality, we deduce
$$T_7=O_{\PP}(\tilde r_{n,n'} \sup_i |T_{7,i} |  h^\alpha)=O_{\PP}(\tilde r_{n,n'} h^{2\alpha})=O_{\PP}\big( (n'n  h^{p+4\alpha})^{1/2}\big),$$
and then $T_7= o_{\PP}(1)$ due to Assumption~\ref{assumpt:asymptNorm_joint}(ii)(a).\;\;$\Box$


\subsection{Convergence of $T_3$ to $0$}
\label{proof:lemma:limit_T_3}

For every $\epsilon>0$, by Markov's inequality,
\begin{align*}
    \PP( |T_3| >\epsilon) \leq
    & \frac{C_{\Lambda''}\tilde r_{n,n'}}{2 n' \epsilon } \sum_{i=1}^{n'} \EE\big[
    \big(\hat \tau_{1,2|\Z = \z'_i} - \tau_{1,2|\Z = \z'_i} \big)^2 \big] \psibm(\z'_i).
\end{align*}
An approximated calculation of
$\EE\big[
\big(\hat \tau_{1,2|\Z = \z'_i} - \tau_{1,2|\Z = \z'_i} \big)^2 \big]$ 
can be obtained following the steps of the proof of Lemma~\ref{lemma:asymptNorm_hatTau}.
Indeed, it can be easily seen that the order of magnitude of the latter expectation is the same as the variance of $U_{n,i}(g^*)$, and then of its Hájek projection $\hat U_{n,i} (g)$.
Since the latter variance is $O((nh^p)^{-1})$, we get
$$ \PP( |T_3| >\epsilon) \leq B \frac{\tilde r_{n,n'}}{ nh^p \epsilon },$$
for some constant $B$. Since $n'/(nh^p)\to 0$, we get $T_3 = o_{\PP}(1)$, as claimed.
$\;\; \Box$

\section{Technical results concerning the first-step estimator}
\label{section:choice_g}

Three possible choices for $g^*$ are given in~\cite{derumigny2018kernelBased}
\begin{align*}
    g_1(\X_i, \X_j)
    &:= 4 \cdot \1 \big\{ X_{i,1} < X_{j,1} , X_{i,2} < X_{j,2} \big\} - 1, \\
    g_2(\X_i, \X_j)
    &:= \1 \big\{ (X_{i,1} - X_{j,1}). (X_{i,2} - X_{j,2}) > 0 \big\}
    - \1 \big\{ (X_{i,1} - X_{j,1}). (X_{i,2} - X_{j,2}) < 0 \big\}, \\
    g_3(\X_i, \X_j)
    &:= 1 - 4 \cdot \1 \big\{ X_{i,1} < X_{j,1} , X_{i,2} > X_{j,2} \big\},
\end{align*}
where $\1$ is the indicator function. In the following, we assume that we have chosen $g^*$ as one of the $g_k$ for a fixed $k \in \{1,2,3\}$.

\begin{assumpt}
    The kernel $K$ is bounded, and set $\| K \|_{\infty} =: C_{K}$. It is symmetrical and satisfies $\int K = 1$, $\int |K| <\infty$.
    This kernel is of order $\alpha$ for some integer $\alpha > 1$: for all $j = 1, \dots, \alpha -1$ and every indices $i_1,\ldots,i_j$ in $\{1,\ldots,p\}$,
    $\int_{\Rb^p} K(\u)  u_{i_1} \dots u_{i_j} \; d\u = 0,$.
    \label{assumpt:kernel_integral}
\end{assumpt}

\begin{assumpt}
    $f_\Z$ is $\alpha$-times continuously differentiable and there exists a constant $C_{K,\alpha}>0$ s.t.,
    for all $\z \in \Zc$,
    $$\int |K|(\u)
    \sum_{i_1, \dots, i_{ \alpha } = 1}^{p}
    |u_{i_1} \dots u_{i_{\alpha }}| \,
    \sup_{t\in [0,1]}\big| \frac{ \partial^{\alpha } f_{\Z}}{ \partial z_{i_1} \dots  \partial z_{i_{\alpha }}} (\z+t\u) \big| \, d\u \leq  C_{K,\alpha}.$$
    \label{assumpt:f_Z_Holder}
\end{assumpt}
\begin{assumpt}
    There exist two positive constants $f_{\Z, min}$ and $f_{\Z, max}$ such that,
    for every $\z \in \Zc$, $f_{\Z, min} \leq f_{\Z}(\z) \leq f_{\Z, max}$.
    \label{assumpt:f_Z_max}
\end{assumpt}

\begin{lemma}
    Under Assumptions \ref{assumpt:kernel_integral},~\ref{assumpt:f_Z_Holder} and~\ref{assumpt:f_Z_max},
    we have for any $t > 0$,
    \begin{align*}
        \PP \bigg( \big| \hat f_{\Z}(\z)-f_{\Z}(\z) \big|
        \geq \frac{  C_{K, \alpha} h^{\alpha}}{ \alpha !} + t \bigg)
        \leq 2 \exp \bigg( - \frac{n h^p t^2}{2 f_{\Z, max} \int K^2 + (2/3) C_K t} \bigg).
    \end{align*}
    \label{lemma:bound_f_hat_f}
\end{lemma}

\begin{lemma}
    Under Assumptions \ref{assumpt:kernel_integral}-\ref{assumpt:f_Z_max} and if
    $ C_{K, \alpha} h^{\alpha}  / \alpha  ! \,
    < f_{\Z, min}$, the estimator $\hat f_{\Z}(\z)$ is strictly positive with a probability larger than
    $$1 - 2 \exp \Big( - n h^p \big( f_{\Z, min} - C_{K, \alpha} h^{ \alpha}/\alpha ! \big)^2
    / \, \big( 2 f_{\Z, max} \int K^2 + (2/3) C_K ( f_{\Z, min} -
    C_{K, \alpha} h^{ \alpha}/\alpha !) \big) \Big).$$
    \label{lemma:probaTau_Z_valid}
\end{lemma}

\begin{assumpt}
    For every $\x \in \Rb^2$, $\z \mapsto f_{\X, \Z}(\x, \z)$ is differentiable almost everywhere up to the order $ \alpha $, $\z\in \Zc$. For every $0 \leq k \leq  \alpha $
    and every $1 \leq i_1, \dots, i_{ \alpha } \leq p$, let
    \begin{equation*}
        \Hc_{k,\vec{\iota}}(\u,\v,\x_1,\x_2,\z):= \sup_{t \in [0,1]} \bigg|
        \frac{\partial^{k} f_{\X, \Z}}{\partial z_{i_1} \dots  \partial z_{i_k}}
        \Big( \x_1, \z + t \u \Big)
        \frac{\partial^{ \alpha -k} f_{\X, \Z}}
        {\partial z_{i_{k+1}} \dots  \partial z_{i_{ \alpha }}}
        \Big( \x_2, \z + t \v \Big)\bigg|,
    \end{equation*}
    denoting $\vec{\iota}=(i_1,\ldots,i_\alpha)$.
    Assume that $\Hc_{k,\vec{\iota}}(\u,\v,\x_1,\x_2,\z)$ is integrable and there exists a finite constant $C_{\X\Z, \alpha} > 0$, such that, for every $\z \in \Zc$,
    \begin{align*}
        \int |K|(\u) |K|(\v)
        \sum_{k=0}^{ \alpha } \binom{ \alpha }{k}
        \sum_{i_1, \dots, i_{ \alpha } = 1}^{p}
        \Hc_{k, \vec{\iota}}(\u,\v,\x_1,\x_2,\z)
        |u_{i_1} \dots u_{i_k} v_{i_{k+1}} \dots v_{i_{ \alpha }}|
        \, d\u \, d\v\, d\x_1\, d\x_2
    \end{align*}
is less than $C_{\X\Z, \alpha}$.
    \label{assumpt:f_XZ_Holder}
\end{assumpt}

\begin{lemma}[Exponential bound for the estimated conditional Kendall's tau]
    Under Assumptions \ref{assumpt:kernel_integral}-\ref{assumpt:f_XZ_Holder},
    for every $t>0$ such that
    $ C_{K, \alpha} h^{\alpha} /  \alpha ! + t \leq f_{\Z, min}/2$ and every $t'>0$, we have
    \begin{align*}
        &\PP \Bigg( |\hat \tau_{1,2|\Z=\z} - \tau_{1,2|\Z=\z} |
        > c_k \bigg(1 + \frac{16 f_{\Z, max}^2}{f_{\Z, min}^3}
        \Big( \frac{ C_{K, \alpha} h^{\alpha}}
        { \alpha  !} + t \Big) \bigg)\times
        \bigg( \frac{C_{\X\Z, \alpha}   h^\alpha}
        {f_{\Z}^2(\z)  \alpha  !} + t' \bigg)
        \Bigg) \\
        &  \leq 2 \exp \Big( - \frac{n h^p t^2}{2 f_{\Z, max} \int K^2 + (2/3) C_K t} \Big)
        + 2 \exp \Big( - \frac{(n-1) h^{2p} t'{}^2 f_{\Z,min}^4}{4 f_{\Z, max}^2 (\int K^2)^2 + (8/3) C_K^2 f_{\Z, min}^2 t'} \Big),
    \end{align*}
    with $c_1 := c_3 := 4$ and $c_2 := 2$.
    \label{lemma:exponential_bound_KendallsTau}
\end{lemma}

\begin{rem}
In Lemma~\ref{lemma:probaTau_Z_valid} and~\ref{lemma:exponential_bound_KendallsTau}, $f_{\Z,min}$ can be replaced by $f_{\Z}(\z)$. Moreover,
when the support of $K$ is included in $ [-c,c]$ for some $c>0$, $f_{\Z,max}$ can be replaced by
$\sup_{\tilde \z \in \Vc(\z,\epsilon)} f_{\Z}(\tilde\z)$, denoting by $\Vc(\z,\epsilon)$ a closed ball
of center $\z$ and any radius $\epsilon>0$, when $n\, c < \epsilon$.
\end{rem}

\begin{lemma}[Consistency]
    Under Assumption~\ref{assumpt:kernel_integral}, if
    $n h_{n}^p \to \infty$,
    $\lim K(\t) | \t |^p =0$ when $|\t | \to \infty$,
    $f_\Z$ and $\z\mapsto \tau_{1,2|\Z=\z}$ are continuous on $\Zc$, then $\hat \tau_{1,2|\Z=\z}$ tends to $\tau_{1,2|\Z=\z}$ in probability, when $n\to\infty$.
    \label{lemma:consistency_hatTau}
\end{lemma}

To derive the asymptotic law of this estimator, we will assume:
\begin{assumpt}
    (i) $n h_{n}^p \to \infty$ and $n h_{n}^{p+2\alpha} \to 0$;
    (ii) $K(\, \cdot \,)$ is compactly supported.
    \label{assumpt:asymptNorm}
\end{assumpt}

\begin{lemma}[Asymptotic normality]
    Assume~\ref{assumpt:kernel_integral},~\ref{assumpt:f_XZ_Holder},~\ref{assumpt:asymptNorm}, that the $\z'_i$ are distinct and that
    $f_\Z$ and
    $\z\mapsto f_{\X,\Z}(\x,\z)$ are continuous on $\Zc$, for every $\x$.

    Then, $(n h_{n}^p)^{1/2} \left( \hat \tau_{1,2|\Z=\z'_i} - \tau_{1,2|\Z=\z'_i} \right)_{i=1, \dots, n'}
        \indistrto \Nc (0,  \HH)$ as $ n \to \infty$,
    where $\HH$ is a $n' \times n'$ real matrix defined by
    \begin{align*}
        [\HH]_{i, j} = \frac{ 4 \int K^2 \1_{ \{ \z'_i = \z'_j \} }
        }{ f_\Z(\z'_i)} \big\{
        \EE[\tilde g(\X_1,\X) \tilde g(\X_2,\X) | \Z = \Z_1 = \Z_2 = \z'_{i}]
         - \tau_{1,2|\Z=\z'_i}^2
        \big\},
    \end{align*}
    for every $1 \leq i, j \leq n'$, and $(\X, \Z)$, $(\X_1, \Z_1)$, $(\X_2, \Z_2)$ are independent copies, where $\tilde g$ is the symmetrized version $\tilde g(\x_1, \x_2) := g^*(\x_1, \x_2) + g^*(\x_2, \x_1))/2$.
    \label{lemma:asymptNorm_hatTau}
\end{lemma}


\pagebreak
\setcounter{equation}{0}
\setcounter{figure}{0}
\setcounter{table}{0}
\setcounter{page}{1}
\makeatletter
\renewcommand\appendixname{Supplement}
\renewcommand{\theequation}{S\arabic{equation}}
\setcounter{section}{0}
\setcounter{subsection}{0}

\if1\blind
{
\begin{center}
    {\LARGE\bf Supplementary figures on a simulated sample}
    
    \bigskip
    
  {\large Alexis Derumigny$^1$
    and
    Jean-David Fermanian\footnote{
    CREST-ENSAE, 5, avenue Henry Le Chatelier,
    91764 Palaiseau cedex, France. \\
    Email adresses: alexis.derumigny@ensae.fr, jean-david.fermanian@ensae.fr. \\
    This research has been supported by the Labex Ecodec.
    }}
    
    \end{center}
} \fi

\if0\blind
{
  \begin{center}
    {\LARGE\bf Supplementary figures on a simulated sample}
    \end{center}
  \medskip
} \fi

%
%
%


\begin{figure}[ht]
    \centering
    \includegraphics[height = 10cm]{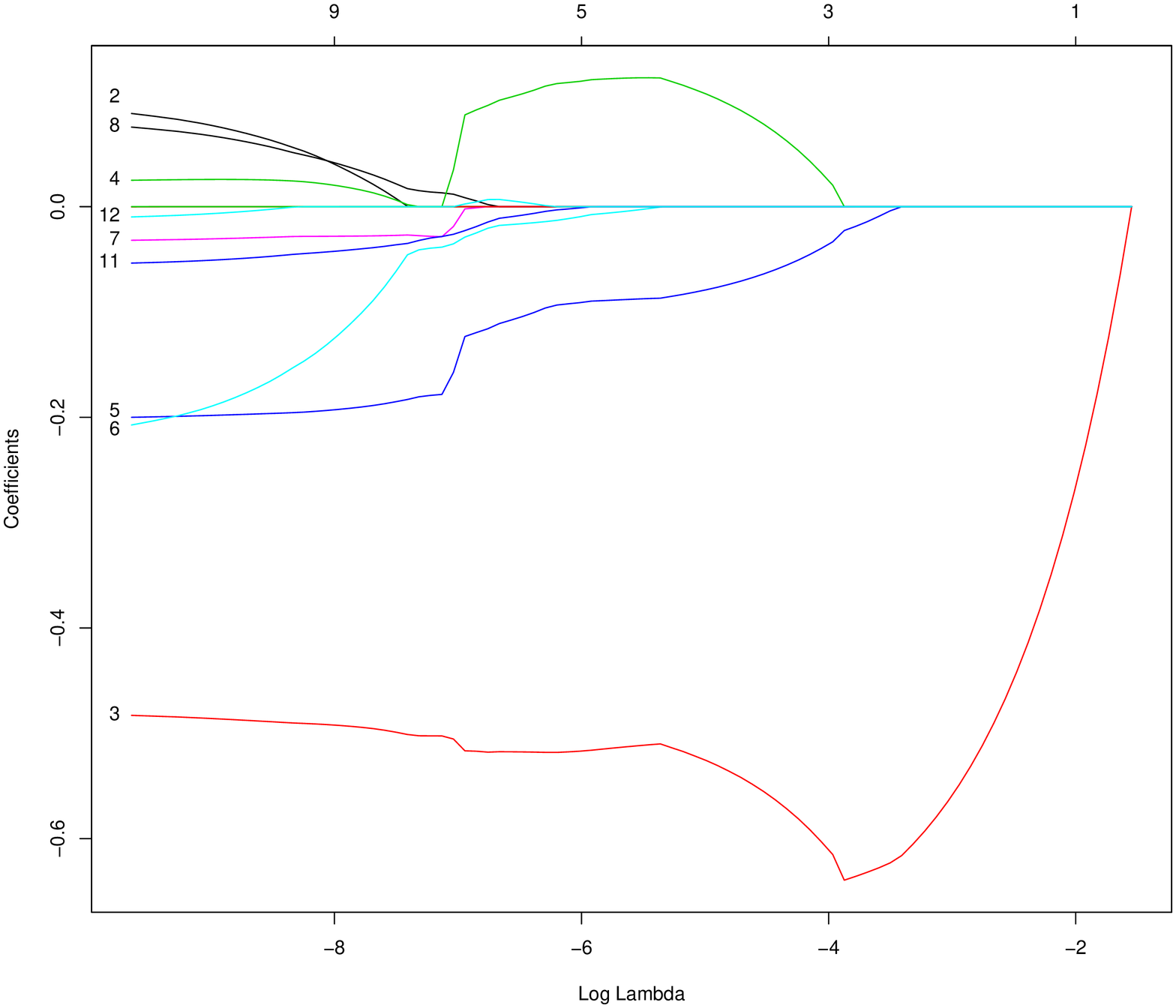}
    \caption{Evolution of the estimated non-zero coefficients as a function of the regularization parameter $\lambda$. The non-zero coefficients are $\beta_1=3/4$ and $\beta_3=3/4$. Note that the coefficients $\hat \beta_2$, $\hat \beta_5$ and $\hat \beta_9$ coefficients are always zero (and are not displayed).}
    \label{fig:plot_poly2_coeff_lambda}
\end{figure}


\begin{figure}[ht]
    \centering
    \includegraphics[height = 10cm]{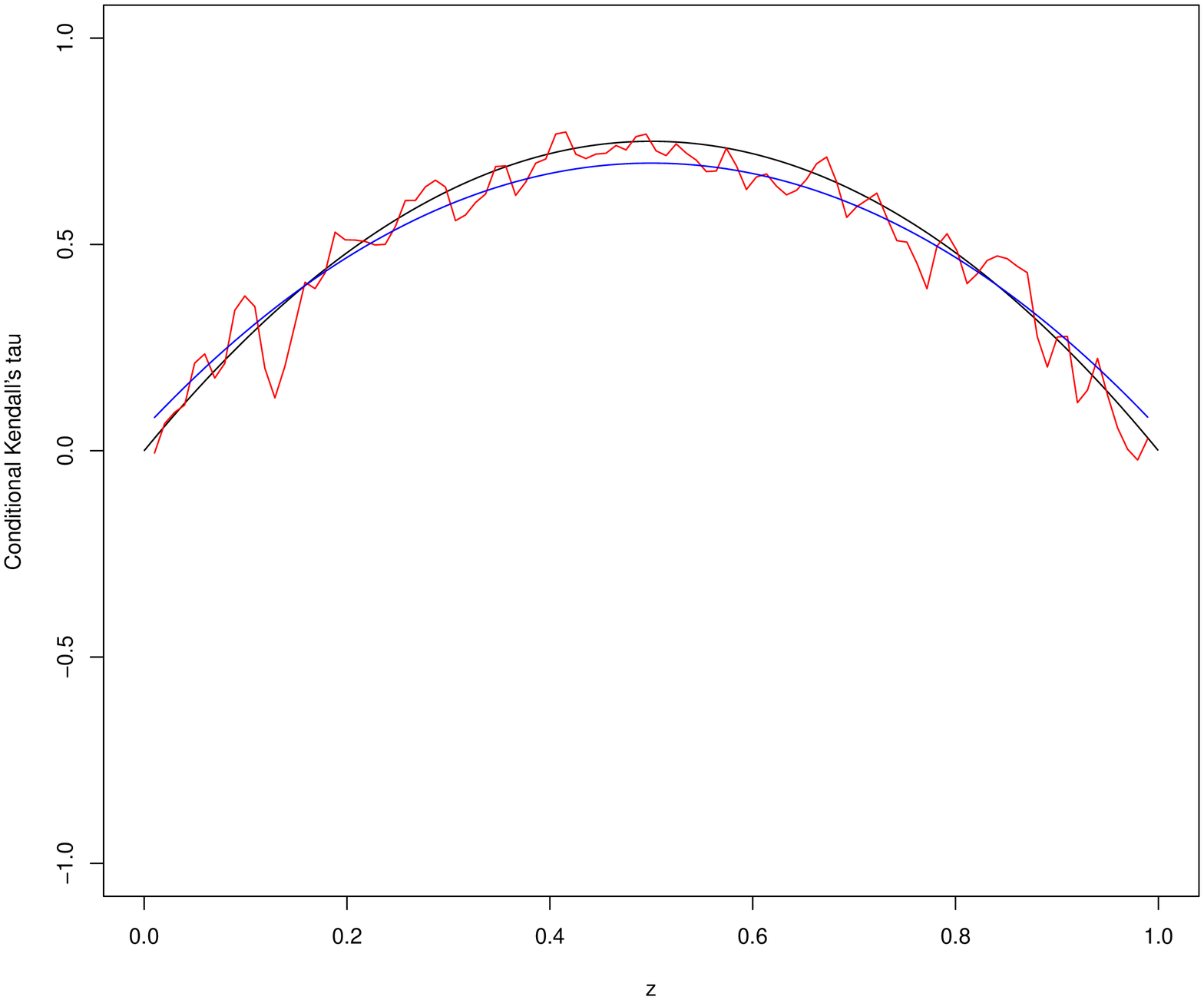}
    \caption{True conditional Kendall's tau $\tau_{1,2|Z = z}$ (black curve), estimated conditional Kendall's tau $\hat \tau_{1,2|Z = z}$ (red curve), and prediction $\Lambda^{(-1)} \big( \psibm (z)^T \hat \beta \big)$ (blue curve) as a function of $z$. For the blue curve, the regularization parameter is $2 \hat \lambda^{cv}\simeq 0.034$ where $\hat \lambda^{cv}$ is selected by Algorithm \ref{algo:cross_validation_lambda}.}
    \label{fig:plot_poly2_comparison}
\end{figure}

\FloatBarrier

\end{document}